\documentclass{amsart}
\RequirePackage[OT1]{fontenc}
\usepackage{amsfonts, amsmath, amssymb, amsthm,url,dsfont, gensymb,multicol,bbm}
\usepackage[margin=3cm,footskip=1cm]{geometry}
\usepackage[english]{babel}
\RequirePackage[colorlinks,citecolor=blue,urlcolor=blue]{hyperref}
\usepackage{booktabs}

\RequirePackage[numbers]{natbib}
\usepackage[above,below,section]{placeins}

\usepackage[shortlabels]{enumitem}
\setlist{
	listparindent=\parindent,
	parsep=0pt,
}

\theoremstyle{definition}
\theoremstyle{remark} 

\usepackage[dvipsnames]{xcolor}
\usepackage{tikz}
\usepackage{amsmath,bm,bbm,amsthm, amssymb}
\usepackage{fullpage}
\usepackage{color}
\usepackage{dsfont}
\usepackage{animate}
\usepackage{etoolbox}
\usepackage{comment}
\usepackage{pgf}
\usetikzlibrary{calc}
\usetikzlibrary{patterns}
\usetikzlibrary{arrows}
\usetikzlibrary{decorations.pathreplacing}
\usepackage[utf8]{inputenc}
\usepackage{pgfplots}

\usepackage{tcolorbox}
\usepackage{adjustbox}
\usepackage[final]{pdfpages}
\usepackage[ruled,vlined,linesnumbered,algosection,resetcount]{algorithm2e}
\usepackage{blkarray}
\usepackage{arydshln}
\usepackage{wrapfig}
\usepackage[font=scriptsize, labelfont=sc]{caption}

\theoremstyle{plain}
\newtheorem{theorem}{Theorem}
\newtheorem{proposition}[theorem]{Proposition}
\newtheorem{corollary}[theorem]{Corollary}
\newtheorem{lemma}[theorem]{Lemma}

\theoremstyle{definition}
\newtheorem{definition}[theorem]{Definition}

\newtheorem{example}[theorem]{Example}

\theoremstyle{remark}
\newtheorem{remark}[theorem]{Remark}

\def\mc{\mathcal}
\def\ms{\mathsf}
\DeclareMathOperator{\Cov}{Cov}
\DeclareMathOperator{\Var}{\mathsf{Var}}

\newcommand{\X}{\mc X}

\newcommand{\cb}[1]{{\color{black} #1}}

\renewcommand{\Im}{\textrm{Im}}

\DeclareMathOperator{\diam}{diam}
\newcommand{\E}{\mathbb{E}}
\def\P{\mathbb P}
\newcommand{\FF}{\mc F}
\newcommand{\PP}{\mc P}

\newcommand{\N}{\mathbb N}
\newcommand{\Z}{\mathbb{Z}}
\newcommand{\R}{\mathbb{R}}

\newcommand{\one}{\mathds{1}}

\newcommand{\eps}{\varepsilon}
\DeclareMathOperator{\de}{deg}

\DeclareMathOperator{\su}{sup}

\def\d{{\mathrm{d} }}

\def\xx{\boldsymbol x}
\def\dist{\mathsf{dist}}

\definecolor{wiasblue}   {cmyk}{1.0, 0.60, 0, 0}
\definecolor{mlugreen}{RGB}{172,6,52}
\definecolor{darmstadt}{RGB}{135,206,250}

\def\ni{\noindent}
\def\E{\mathbb E}

\def\P{\mathbb P}

\def\R{\mathbb R}
\def\X{\mathbb X}
\def\Z{\mathbb Z}
\def\mc{\mathcal}
\def\ms{\mathsf}

\def\lrsa{\leftrightsquigarrow}

\def\La{\Lambda}
\def\a{\alpha}

\def\s{\sigma}
\def\su{\subseteq}
\def\bs{\boldsymbol}

\def\g{\gamma}
\def\b{\beta}
\def\de{\delta}

\def\om{\omega}

\def\es{\varnothing}
\def\one{\mathbbmss{1}}

\def\De{\Delta}

\def\ff{\infty}
\def\tff{\uparrow\infty}

\def\vp{\varphi}

\def\d{{\rm d}}
\def\k{\kappa}

\def\YY{\mc Y}

\def\NN{\mc N}
\def\BB{\mc B}
\def\PP{\mc P}
\def\FF{\mc F}

\def\XX{\mc X}

\def\AA{\mc A}

\def\f{\frac}

\def\im{\item}
\def\sm{\setminus}
\def\th{\theta}
\def\bep{\begin{proof}}
	\def\enp{\end{proof}}
\def\bepr{\begin{proposition}}
	\def\enpr{\end{proposition}}
\def\bec{\begin{corollary}}
	\def\enc{\end{corollary}}
\def\bea{\begin{align}}
	\newcommand\eea{\end{align}}
\def\beas{\begin{align*}}
	\def\eeas{\end{align*}}
\def\bet{\begin{theorem}}
	\def\ent{\end{theorem}}
\def\bee{\begin{example}}
	\def\ene{\end{example}}

\def\bede{\begin{definition}}
	\def\ende{\end{definition}}
\def\ber{\begin{remark}}
	\def\enr{\end{remark}}
\def\beca{\begin{cases}}
	\def\enca{\end{cases}}
\def\bel{\begin{lemma}}
	\def\enl{\end{lemma}}
\def\been{\begin{enumerate}}
	\def\enen{\end{enumerate}}

\def\dim{\ms{dim}}
\def\beit{\begin{itemize}}
	\def\enit{\end{itemize}}
\def\befr{\begin{frame}}
	\def\enfr{\end{frame}}
\def\ti{\times}

\def\Var{\ms{Var}}
\def\Cov{\ms{Cov}}

\def\ba{\,|\,}

\def\diam{\ms{diam}}

\def\Im{\ms{Im}}

\renewcommand\le{\leqslant}
\renewcommand\ge{\geqslant}

\def\bN{{\mathbf N}}

\def\xx{{\bs x}}
\def\ww{{\bs w}}

\def\Nlf{\mathbf N}

\def\dist{\ms{dist}}

\def\becbb{\begin{center}\begin{tcolorbox}[{colback=Dandelion!20}]}
		\def\encbb{\end{tcolorbox}\end{center}}
\def\becb{\begin{center}\begin{tcbox}[{colback=Dandelion!20}]}
		\def\encb{\end{tcbox}\end{center}}

\def\bef{\begin{figure}[!h]}
	\def\enf{\end{figure}}
\def\bels{\begin{lstlisting}}
	\def\enls{\end{lstlisting}}

\def\betp{\begin{tikzpicture}}
	\def\entp{\end{tikzpicture}}
\def\co{\colon}
\def\endo{\end{document}}

\def\pa{\partial}

\def\N{\mathbb N}

\def\dtv{d_{\ms{TV}}}

\def\Tra{T^{\ms{rad}}}

\def\Pd{\PP^{}}
\def\Pds{\PP^{*}}

\def\X{\XX}

\def\dk{d_{\ms K}}

\def\tdc{T^{\ms{dc}}}

\begin{document}
	
	\title{Normal approximation for Gibbs processes\\ via disagreement couplings}
	
	\begin{abstract}
	This work improves the existing central limit theorems (CLTs) for geometric functionals of Gibbs processes in three aspects. First, we derive a CLT for weakly stabilizing functionals, thereby improving on the previously used assumption of exponential stabilization. Second, we show that this CLT holds for interaction ranges up to the percolation threshold of the dominating Poisson process. This avoids imprecise branching bounds from graphical construction. Third, by constructing simultaneous couplings of several Palm processes for Gibbs functionals, we provide a quantitative CLT in terms of Kolmogorov bounds for normal approximation. An important conceptual ingredient in these advances \cb{is the extension} of disagreement coupling adapted to unbounded windows and to the comparison at multiple spatial locations.
	\end{abstract}
\author{Christian Hirsch}
\author{Moritz Otto}
\address[Christian Hirsch, Moritz Otto]{Department of Mathematics, Aarhus University, Ny Munkegade 118, 8000 Aarhus C, Denmark}
\email{hirsch@math.au.dk, otto@math.au.dk}
\address[Christian Hirsch]{DIGIT Center, Aarhus University, Finlandsgade 22, 8200 Aarhus N, Denmark}

\author{Anne Marie Svane}
\address[Anne Marie Svane]{Department of Mathematics\\
Aalborg University \\
Aalborg, 8000, Denmark}
\email{annemarie@math.aau.dk}	

\maketitle
	
	%
	%
	\section{Introduction}
\label{sec:int}
%
%
Because of their ability to encode possibly highly complex interactions, Gibbs point processes are applied in a broad range of domains such as biology, ecology, materials science and telecommunication networks \cite{Baddeley:Rubak:Wolf:15,g1}. However, even for simple models, the use of Gibbs point processes comes at massive computational costs as the simulation requires elaborate Markov chain Monte Carlo methods \cite{Baddeley:Rubak:Wolf:15}. This issue becomes particularly pressing when devising goodness-of-fit tests on large datasets. Therefore, an attractive approach is to develop test statistics that are asymptotically normal on large windows. Then, to develop a hypothesis test only the mean and variance under the null model are needed.
This explains the need for central limit theorems (CLTs) on functionals of Gibbs point processes.

%
%
The last decade was \cb{marked} by vigorous research activities in the asymptotic theory of statistics that can be written as a sum of certain scores evaluated at the points of a point process with rapidly decaying correlations. 
Improving on the methods used in the classical setting of Poisson point processes in \cite{penrose}, CLTs for Gibbs point processes could be derived in \cite{gibbs_limit,gibbsCLT} by relying on the graphical construction of Gibbs processes \cite{ferrari}.

However, despite these recent advances, some aspects of the limit theory for Poisson processes could so far not been transferred to Gibbs point processes:
\been
\im All of \cite{CX22,gibbs_limit,gibbsCLT} are formulated under the condition of exponential stabilization. This makes it difficult to apply \cite{gibbs_limit,gibbsCLT} for delicate topological functionals, such as the persistent Betti numbers. 
\im Both \cite{gibbs_limit,gibbsCLT}  rely on the graphical construction of Gibbs processes from \cite{ferrari}, thereby imposing restrictive constraints on the interaction range of the Gibbs process for which a (partially quantitative) CLT can be established.
\im In \cite{benes}, disagreement coupling is used to establish a CLT with a Gibbs particle process as input {thus allowing a larger interaction range.} However, the investigations are restricted to U-statistics and do not provide quantitative error bounds for normal approximation. 	Very recently, \cite{CX22} derived  a quantitative CLT with optimal convergence rates (up to logarithmic corrections) and for a very general class of $\b$-mixing point processes. However, the convergence is quantified with respect to the Wasserstein distance and the methods do not extend easily to the  Kolmogorov distance.
\enen

%
%
In our work, we will address all of  the shortcomings mentioned above. The key tool for achieving the improvements will be a more refined analysis of the correlation structure of Gibbs processes using disagreement coupling. While this technique was initially suggested  for studying the spatial correlation structure of lattice-Gibbs point processes \cite{maes}, it has recently been successfully extended to derive Poisson approximation theorems for continuum Gibbs processes \cite{dp}. However, the construction in \cite{dp} is not adapted to the setting of {unbounded domains} or the coupling at several locations, which are critical for the CLT improvements derived in the present article. Therefore, a major part of our investigation is devoted to making disagreement coupling more flexible so that it can accommodate the needs from spatial limit theorems. Equipped with these conceptual advances, we improve the existing CLTs for Gibbs point processes in three aspects.
\been
\im We extend the CLT for weakly stabilizing functionals from \cite{penrose} to the setting of Gibbs point processes. In particular, this yields a CLT for persistent Betti numbers {and total edge length of minimal spanning trees.}
\im We prove a CLT for Gibbs point processes under the condition that the interaction range is smaller than the critical value of continuum percolation of an associated Poisson point process. This improves on the restrictive constraints coming from the branching bounds of the graphical construction in \cite{ferrari}. In other words, our work makes normal approximation available for point processes that were out of reach in \cite{gibbs_limit}. 
\im We show how the convergence in Wasserstein distance from \cite{CX22} can be upgraded to a convergence in Kolmogorov distance provided stricter conditions are imposed on the point process and the functional. To that end, we introduce simultaneous couplings of several Palm versions for Gibbs processes.
\enen

The rest of this manuscript is organized as follows. First, in Section \ref{sec:mod}, we introduce the model and main results, Theorems \ref{thm:1} and \ref{thm:2}. Then, in Section \ref{sec:exa}, we discuss examples satisfying the conditions of these results. {We review the definition of disagreement coupling in Section \ref{sec:const} and establish the properties needed for the proofs of the main results.} We prove Theorems \ref{thm:1} and \ref{thm:2} in Sections \ref{sec:weak} and \ref{sec:gibbs}, respectively. Finally, we provide { a detailed comparison with the existing normal approximation literature in Appendix \ref{sec:EDD} and }a technical auxiliary result in Appendix \ref{sec:iota}.

%
%
\section{Model and main results}
\label{sec:mod}
The main results of this article, Theorems \ref{thm:1} and \ref{thm:2} below, are CLTs for stabilizing functionals on Gibbs point processes. First, we discuss in greater detail the considered classes of Gibbs point processes. For a concise introduction to Gibbs point processes, we refer the reader to \cite{dereudre}.

%
%
\subsection{Conditions on the Gibbs process}
\label{ssec:gibbs}

Let $\R^d$ be Euclidean space equipped with the Borel $\sigma$-algebra $\BB^d$. 
Define
$\BB^d_0$ to be the set of bounded Borel sets.
Let $\Nlf$ be the space of all locally finite subsets of $\R^d$
and let $\NN$ denote the smallest $\s$-algebra on $N$
such that $\mu\mapsto \mu(B) :=\#(\mu \cap B)$ is measurable for
all $B\in\BB_0^d$, where $\#A$ denotes the cardinality of a set $A$. Let $\Nlf_0\su \Nlf$ denote the subspace of finite subsets and for $B\in \BB^d$, let $\Nlf_B\su \Nlf$ denote the locally finite subsets of $B$.   
A {\em point process} is a measurable map
from $(\Omega,\FF)$ to $(\Nlf,\NN)$ where 
$(\Omega,\FF,\P)$ is some probability space. 

To define Gibbs point processes, we start from an energy 
functional $E:\Nlf_0 \to (-\ff,\ff] $ which satisfies: 
\begin{itemize}
	\item {\em Non-degeneracy:} $E(\es)=0$.
	{\item {\em Hereditary:} For all $\psi\in \Nlf_0$ and $x\in \R^d$, $E(\psi)=\ff$ implies $E(\psi  \cup \{x\} )=\ff$.
		\item {\em Stability:} There is a constant $C>0$ with
		\begin{equation}\label{cond:stability}
			\inf_{\psi \in \Nlf_0} E(\psi ) \ge -C\#\psi. 
		\end{equation}
	}
\end{itemize}

The {\em local energy}
$J\colon\bN_0\times\bN_0\to (-\ff,\ff]$ is given by
\begin{align}\label{eHamilton}
	J(\vp,\psi):= E(\vp \cup \psi) - E(\psi),
\end{align}
which encodes the energy increase when adding the point configuration $\vp$ to $\psi$.  We will assume that $J$ has \emph{finite interaction range}. That is,  there exists $r_0 > 0$ such that
\begin{equation*}
	J(x, \psi) = J(x, \psi\cap B_{r_0}(x))
\end{equation*}
for all $x \in \R^d$ and $\psi \in \Nlf_0$, where  
$$B_{r_0}(x):=  \{y\in \R^d\co |y - x|\le r_0\}$$
denotes the Euclidean ball with radius $r_0 > 0$ centered at $x \in \R^d$. We write $B_{r_0}:=B_{r_0}(x)$.
Note that the assumption \cb{of finite interaction range} allows \cb{us} to define $J(\vp, \psi)$ also for $\psi \in \Nlf$.

For $B\in\BB_0^d$, the {\em partition function}
$Z_B\colon\bN\to (0,\ff)$ is defined by
\begin{align*}
	Z_B(\psi):=\E [e^{-J(\Pd_B,\psi)} ],\quad \psi\in\bN,
\end{align*}
where $\Pd_B$ denotes a unit intensity Poisson process 
on $B$. Condition \eqref{cond:stability} indeed ensures that $Z_B(\psi)<\ff$ and since $J(\es,\psi)=0$ for all $\psi\in\bN$, we also have that $ Z_B(\psi)\ge e^{-|B|}>0.$

For $B\in \BB_0^d$ and $\psi\in \Nlf$, we now define the {\em finite-volume Gibbs point process} $\XX(B,\psi)$ on $B$ with boundary conditions $\psi$ as the point process with density
\begin{align}\label{eGibbsmeasure}
	\vp \mapsto Z_B(\psi)^{-1}e^{-J(\vp ,\psi)}
\end{align}
with respect to the unit intensity Poisson process on $B$. This is well-defined since $0< Z_B(\psi)<\ff$.

The finite-volume Gibbs process is sometimes, e.g.\ in \cite{dp}, introduced in terms of the Papangelou conditional intensity (PI) {$\k\colon\R^d \times\Nlf_0\to [0,\infty)$  defined by
	\begin{equation*}
		\k(x,\psi) = e^{-J(x,\psi)}.
	\end{equation*}
}
Then, $\XX(B,\psi)$ can be characterized as the unique point process satisfying the {\em GNZ equation}
\begin{align}\label{eGNZ_finite}
	\E\Big[\sum_{X_i \in \XX(B,\psi)}  f(X_i,\XX(B,\psi))\Big]=
	\E\Big[ \int_B f(x,\XX(B,\psi) \cup\{x\}) \k(x,\XX(B,\psi) \cup \psi)\,\lambda(\d x)\Big],
\end{align}
for each $f\colon B \times\Nlf \to [0,\ff)$ that is measurable \cite{Georgii76,NgZe79}.

\cb{
	Consider the finite volume Gibbs point process $\XX(B,\psi)$ on $B\in\BB_0$ with boundary conditions $\psi\in \Nlf$. 
	Then, $\XX$ satisfies the {\em DLR equations}: For $C\su B$ and all $\vp \in \Nlf_{B \backslash C}$ outside a set of measure 0 with respect to $\XX(B, \psi)\cap (B\sm C)$, }
\begin{align}\label{eq:DLR1}
	\P(\XX(B,\psi)\cap C \in\cdot \mid \XX(B,\psi)\cap ({B\backslash C}) = \vp ) = \P(\XX(C,\psi \cup \vp)\in \cdot).
\end{align}
In words, the conditional distribution of $\XX(B,\psi)\cap C  $ given \cb{$ \XX(B,\psi)\cap ({B\backslash C})$} is again a Gibbs point process, namely the Gibbs process on $C$ with boundary conditions \cb{$\psi \cup( \XX(B,\psi)\cap ({B\backslash C}))$}.

The GNZ equations can be used to define Gibbs processes in unbounded domains. A point process $\XX$ is said to be an {\em infinite-volume Gibbs point process} on $\R^d$ if 
\begin{align}\label{eGNZ}
	\E\Big[\sum_{X_i \in \XX} f(X_i,\XX)\Big]=
	\E\Big[ \int f(x,\XX \cup \{x\}) \k(x,\XX)\,\lambda(\d x)\Big],
\end{align}
for each measurable $f\colon\R^d\times\Nlf\to [0,\ff)$. 

While the finite-volume Gibbs process was explicitly defined in \eqref{eGibbsmeasure}, \cb{the} existence and in particular \cb{the} uniqueness of infinite-volume Gibbs processes is generally far from trivial, see e.g.\ \cite{betsch,dereudre}. We therefore make some further assumptions on $\k$. \cb{For $A,B\su \R^d$, we let $A\oplus B:=\{x+y\mid x\in A, y\in B\}$ and for $x\in \R^d$, we write $x+A:=\{x\}\oplus A$. Moreover, we define $B_r(A)=A\oplus B_r$. Then, the Boolean model with intensity $\a_0$ and balls of radius $r/2$  is defined as $B_{r/2}(\Pd)$. 
	Let $r_c = r_c(\a_0)$ denote the critical threshold of continuum percolation for a Poisson point process $\Pd$ with intensity $\a_0 > 0$.  That is, the infimum of all $r > 0$ such that  $B_{r/2}(\Pd) $ has an unbounded connected component with positive probability.}

Throughout the paper, we will work under the following assumption. 
\begin{align}\label{as:(A)}
	\text{The PI is 	bounded $0 \le \k \le \alpha_0$ and has \cb{a} finite interaction range $r_0 < r_c(\a_0)$. }
\end{align}
Then, for any given PI $\k$, there exists a unique infinite-volume Gibbs point process  (see \cite[Thm. 5.4]{dereudre} or  Proposition \ref{pr:uniq} below).  We now give two well-known examples of Gibbs point process models satisfying the assumptions.

\cb{
	\bee[Pairwise interaction processes]
	In this example, the local energy is given in terms of a non-negative pair potential $v: [0,r_0]\to [0,\infty]$,
	$$J(x,\psi):= -\log(\a_0) + \sum_{y\in \psi \cap B_{r_0}(x)}v(|x-y|).$$
	The special case $\Psi = \beta \one_{[0,r_0]}$ where $\beta>0$ and $\one_{[0,r_0]}$ is the indicator for the interval $[0,r_0]$ is known as  the Strauss process. If $\beta=\ff$, one obtains the hard-sphere model.
	\ene
}

\bee[Area-interaction process \cite{bvl}]
In this example, we use for $\a_0, r_0>0$ and $\g \in[0, 1]$ the potential
$$J(x,\psi):= \cb{-}\log(\a_0) {+} V(x, \psi) \log(\g),$$
where 
\begin{align*}
	V(x,\psi):=\Big|B_{r_0/2}(x)\sm\bigcup_{y \in \psi \cap B_{r_0/2}(x)} B_{r_0/2}(y)\Big|.
\end{align*}
\ene
In both \cb{examples}, the PI  $\k$ is \emph{translation-invariant}, i.e., $\k(x, \psi) = \k(x+y, \psi +y)$
for all $x,y\in \R^d$ and $\psi\in \Nlf_0$.

We finally mention that an infinite-volume Gibbs process also satisfies the DLR equations:
It holds for all $\vp\in \Nlf_{B^c}$ \cb{outside a set of $\XX\cap B^c$-measure zero that}
\begin{equation}\label{eq:DLR2}
	\P(\XX\cap B \in\cdot \mid \XX \cap{ B^c} = \vp) = \P( \XX(B, \vp) \in \cdot). 
\end{equation}
That is, given the boundary conditions $\XX \cap{ B^c}$, the distribution of $\XX\cap B $ is a finite volume Gibbs process.

We conclude by discussing the relation between Gibbs point processes satisfying Assumption \eqref{as:(A)} and those appearing in the existing  literature \cite{benes,CX22,gibbs_limit,gibbsCLT} on normal approximation for Gibbs processes. Since we do not need the precise definitions in the rest of the paper, we keep the remarks very brief at the present point. For the interested reader, we have included a more thorough discussion in Appendix~\ref{sec:EDD}.

The class of Gibbs processes considered in \cite{CX22,gibbs_limit,gibbsCLT} were first introduced in \cite[Section 1.1]{gibbs_limit}. For finite-range interactions, the entire theory relies on a restrictive interaction-range bound formulated in \cite[Equation (1.4)]{gibbs_limit}. As observed in \cite[Remark 3]{benes}, especially for low dimensions $d$, Assumption \eqref{as:(A)} is a substantial improvement of the bound \cite[Equation (1.4)]{gibbs_limit}. In other words, our work makes normal approximation available for point processes that were out of reach in \cite{CX22,gibbs_limit,gibbsCLT}.  
Conversely, we note that our paper requires the interaction range to be finite, whereas \cite{gibbsCLT} also allows some infinite-range Gibbs processes. 

The key achievement of \cite{CX22} is to provide a quantitative CLT for point processes that are not necessarily of Gibbsian type. The main assumption here is the \emph{exponential decay of dependence (EDD)}. The authors verify this condition for the Gibbs process  from \cite{gibbs_limit}. In Section \ref{sec:EDD}, we show that the EDD property also holds under Assumption \eqref{as:(A)}. 

Finally, we mention that a qualitative normal approximation result was obtained in \cite{benes} under the assumption \eqref{as:(A)}. However, this only applies in the  case of $U$-statistics with finite stabilization radius.

%
%
\subsection{Statement of main results and conditions on the functional}
The main results of this article are two CLTs for functionals of Gibbs processes on growing domains, see Theorems \ref{thm:1} and \ref{thm:2}. These results are complementary in the following sense. While Theorem \ref{thm:1} describes a CLT for a general class of translation-invariant Gibbs functionals without convergence rates, Theorem \ref{thm:2} provides a quantitative CLT for the more restricted class of Gibbs functionals that can be written as a sum of scores. Therefore, it is natural that the two results require different sets of conditions. However, loosely speaking, in both cases the conditions \cb{consist of a moment bound and stabilization condition.}

%
%
We now discuss a CLT for weakly stabilizing functionals as the first main result of our paper.	
Let $H: \Nlf_0 \to \R$ be a functional. We say that $H$ is  \emph{translation-invariant} if  $H(\vp) = H(\vp + x)$ for all  $x\in \R^d$ and   $\vp\in \Nlf_0$.
Writing $Q_n :=[-n/2,n/2]^d$,  $n\geq 1$,  we consider functionals $H_n$ of locally finite point patterns $\vp \in \Nlf$ of the form
\begin{equation*}
	H_n (\vp):= H(\vp \cap Q_n)  
\end{equation*}

\cb{	To prove the CLT for $H_n(\XX)$, we need to impose a moment and a   stabilization condition. 	
	The moment condition requires that there are constants $c_{\ms M}, p_{\ms M}>0$ such that for any $n,m\in \N$ and $z\in \Z^d$, 
	\begin{equation}\label{eq:moment_cond}
		\E\big[|H_n( \X) - H_n( \X \sm Q_{z,m})|^5\big]\le c_{\ms M}m^{p_{\ms M}},
	\end{equation}
	where $Q_{z,m}:= z + Q_m$.
	The stabilization condition requires that, almost surely,
	\begin{equation}\label{eq:as_stab}
		H(\XX)-H(\XX\sm Q_l) := \lim_{n\to \infty} (H(\XX\cap Q_{w_n,m_n})-H((\XX\cap Q_{w_n,m_n})\sm Q_l))
	\end{equation}
	exists for all $l\ge 0$ and all sequences $m_n\in \N$ and $w_n\in \Z^d$ such that $\R^d=\bigcup_{k\ge 1} \bigcap_{n\ge k} Q_{w_n,m_n}$. Note that the limit is necessarily independent of the sequence $Q_{w_n,m_n}$.
	This is a \emph{weak stabilization condition} in the sense that there is no assumption on the rate of convergence in \eqref{eq:as_stab}. 
	
	We show in Section \ref{sss:moment}, that the conditions are satisfied if the following deterministic conditions on the functional hold. The conditions will be stated in terms of the \emph{add-one cost operator}
	$$(D_xH)(\vp) := H(\vp  \cup \{x\}) - H(\vp)$$
	for $x\in \R^d$, $\vp \in \Nlf_0$.
	Then, the moment condition \eqref{eq:moment_cond}	will be satisfied if
	there are constants $c_{\ms DM},R_{\ms DM}>0$ 
	such that for any finite point pattern $\vp \su \R^d$,
	\begin{align} 
		\label{eq:m1}
		|D_yH(\vp)| \le c_{\ms DM}\exp\big(c_{\ms DM}\vp(B_{R_{\ms DM}})\big),
	\end{align}
	see Lemma \ref{lem:(i)'}. 
	The  stabilization condition \eqref{eq:as_stab} is satisfied if for every $\vp\in \Nlf $ and any sequence  $Q_{w_n,m_n} $ with $w_n\in \Z$, $m_n\in \N$ and $\bigcup_{k\geq 1} \bigcap_{n\geq k} Q_{w_n,m_n} = \R^d$, the limit
	\begin{equation}
		\label{eq:s1}
		H(\vp \cup \{0\}) -  H(\vp) := \lim_{n\to \infty } D_0H(\vp \cap Q_{w_n,z_n}) 
	\end{equation}
	exists, see Lemma \ref{lem:s_stab}. 
	
	Let $N(0,\s^2)$ denote a normal random variable with variance $\s^2 \ge0$. Then our first main result is the following. The proof will follow from Theorem \ref{thm:gibbs_CLT} and Lemma \ref{lem:(i)'}.
}
\begin{theorem}[CLT for weakly stabilizing functionals on Gibbs processes]
	\label{thm:1}
	Let $\XX$ be an infinite-volume Gibbs point process
	with translation-invariant PI satisfying \eqref{as:(A)}. Let $H$ be a translation-invariant functional on finite point patterns satisfying conditions \cb{\eqref{eq:moment_cond} and \eqref{eq:as_stab}.} Then, 
	\begin{equation*}
		|Q_n|^{-1}\Var(H_n(\XX)  ) \to \sigma^2,
	\end{equation*}
	for some $\s^2 \ge 0$. Moreover,
	\begin{equation*}
		|Q_n|^{-1/2}(H_n(\XX) - \E[H_n(\XX)]) \to N(0,\sigma^2).
	\end{equation*}
\end{theorem}

While \cite{penrose} give\cb{s} a CLT for weakly stabilizing functionals on Poisson point processes, most CLTs for more general point process models in the literature, \cb{e.g.\ \cite{yogesh,CX22,gibbs_limit,gibbsCLT}, only apply to functionals given in terms of an exponentially stabilizing score function, and \cite{benes} even requires finite stabilization range.} The weak stabilization condition \eqref{eq:s1} closely mimics the one in \cite{penrose}, which is, however, only required to hold almost surely in \cite{penrose}. For an almost sure condition in the Gibbs setting, one would have to use \eqref{eq:as_stab}.
Condition \eqref{eq:m1} resembles the power growth condition  in \cite[Equation (1.18)]{yogesh}, except that it requires a fixed radius $R_{\ms DM}$. 

\cb{A formula for the limiting variance $\s^2$ is given in Theorem \ref{thm:gibbs_CLT}. } Positivity of $\s^2$  is not part of the statement in Theorem \ref{thm:1}. In Corollary \ref{cor:betti}, we \cb{demonstrate in an example how }\cb{the} positivity can be checked. The argument easily generalizes to a variety of other functionals.

Our second main result is a quantitative normal approximation result in the \emph{Kolmogorov distance}
\begin{align*}
	\dk(X,Y):=\sup_{u \in \R} |\P(X \le u) - \P(Y \le u)|.
\end{align*}
between real-valued random variables $X,Y$.  Concerning the setup for Theorem \ref{thm:2}, for $n \ge 1$, we consider the functional 
\begin{align}
	{H}(\vp)& :=\sum_{x \in \vp \cap {Q}} g(x, \vp), \qquad \vp \in \mathbf N, \label{eqn:poixidef}
\end{align}
for some measurable score function $g\co\R^d \ti \Nlf  \to [0, \ff)$. We follow the convention that $g(x, \vp) = g(x, \vp \cup \{x\})$ for every locally finite $\vp \su \R^d$ and \cb{$x \not\in \vp$}.  

We first impose a {bounded moment condition}. More precisely, we assume that 
\begin{align} 
	\label{eq:m2}
	\sup_{x_1, \dots, x_5 \in \R^d}	 \E [g(x_1,\XX\cup \{x_1, \dots, x_5\})^{6}]=:c_{\ms{m}}<\ff,
\end{align}
where we note that the $x_1, \dots, x_5$ do not need to be pairwise distinct.

%
%
Second, we require that $g$ is {\em exponentially stabilizing} in the spirit of \cite{yogesh}. To make this precise,  we assume that there is a \emph{stabilization radius}  $R(x,\vp)$ having exponentially decaying probabilities. That is, first, for all locally finite $\vp \su \R^d$ and $x \in \vp$, we have
$$g(x, \vp) = g\big(x, \vp \cap B_{R(x, \vp)}(x)\big).$$
We also assume that the event $\{R(x,\vp)\le r\}$ is measurable with respect to $\vp \cap B_r(x)$. Second, we impose that 
\begin{align}
	\label{eq:s2}
	-	\limsup_{r\tff }\, \sup_{x_1, \dots, x_5 \in \R^d} r^{-1}{\log\P(R(x_1,\XX\cup \{x_1,\dots,x_5\}) > r)} =:c_{\ms{es}}>0.
\end{align}

%
%
\begin{theorem}[Quantitative normal approximation of Gibbsian score sums]
	\label{thm:2}
	Let $\XX$ be an infinite-volume Gibbs point process
	satisfying \eqref{as:(A)}. Assume that \eqref{eq:m2}, \eqref{eq:s2} hold. Moreover, we assume that $g(x, \vp) = 0$ implies $g(x, \vp') = 0$  for all $x\in \R^d$ and  locally finite $\vp \su \vp'\su \R^d$. We also assume that $Q$ is a bounded Borel set with $|Q| > 1$ and $|B_{(\log |Q|)^2}(Q)| \le 2 |Q|$. Then, $H:=H(\XX)$ defined at \eqref{eqn:poixidef} satisfies 
	\begin{align*}
		\dk\Big(\f{H-\E[H]}{\sqrt{\Var(H)}}, N(0,1) \Big) \le c_{\ms{norm}}\f{|Q|(\log |Q|)^{2d}}{\Var(H)^{3/2}},
	\end{align*}
	where $c_{\ms{norm}}=c_{\ms{norm}}(d, \a_0, r_0, c_{\ms{m}}, c_{\ms{es}})$. 
\end{theorem}
The bound on the parallel set $|B_{(\log |Q|)^2}(Q)|$ is a technical condition, which for instance can be checked for sets of the form $Q = a K$ where $K$ is fixed convex and $a >0$ large enough. In particular, this condition disappears in our Corollary \ref{cor:2}. Note that Theorem \ref{thm:2} also applies to finite-volume Gibbs processes since the PI $\k(x,\vp \cup \psi) \one_{Q}(x)$ corresponds to the Gibbs process $\XX(Q,\psi)$. 
As an application, we state two variations of Theorem \ref{thm:2} { concerning the convergence rate when the observation window $Q = Q_n$ increases towards $\R^d$.} For the first version, let $\XX$ be Gibbs process with PI $\k$ satisfying \eqref{as:(A)} and let
$$
H^n:={H^n(\XX):=}\sum_{x \in \XX \cap Q_n} g(x,\XX).
$$
{Often in applications, only $\XX\cap Q_n$ can be observed, hence it is more natural to consider}
$$
H^n_{\cap}:={H^n_{\cap}(\XX):=}\sum_{x \in \XX \cap Q_n} g(x,\XX \cap Q_n)
$$
and modify the assumptions \eqref{eq:m2} and \eqref{eq:s2} by
\begin{align} 
	\sup_{n \in \N}	\sup_{x_1, \dots, x_5 \in \R^d}	 \E [g(x_1,(\XX \cap Q_n)\cup \{x_1, \dots, x_5\})^{6}]&=:c_{\ms{m}}<\ff,\label{eq:m2cap}\\
	- \sup_{n \in \N}	\limsup_{r\tff }\, \sup_{x_1, \dots, x_5 \in \R^d} r^{-1}{\log\P(R(x_1,(\XX \cap Q_n)\cup \{x_1,\dots,x_5\}) > r)} &=:c_{\ms{es}}>0.	\label{eq:s2cap}
\end{align}

\begin{corollary}
	\label{cor:2}
	Let $\XX$ be an infinite-volume Gibbs point process
	satisfying \eqref{as:(A)}. Moreover, we assume that $g(x, \vp) = 0$ implies $g(x, \vp') = 0$  for all $x\in \R^d$ and  locally finite $\vp \su \vp'\su \R^d$.
	\item[(a)] Assume that \eqref{eq:m2}, \eqref{eq:s2} hold. If $\liminf_{n \to \infty} |Q_n|^{-1} \Var(H^n) >0$, we have for $n$ large enough,
	$$
	\dk\Big(\f{H^n-\E[H^n]}{\sqrt{\Var(H^n)}}, N(0,1) \Big) \in O\Big(\f{(\log |Q_n|)^{2d}}{\sqrt{|Q_n|}}\Big).
	$$
	\item[(b)] Assume that \eqref{eq:m2cap}, \eqref{eq:s2cap} hold. If $\liminf_{n \to \infty} |Q_n|^{-1} \Var(H_{\cap}^n) >0$, we have for $n$ large enough,
	$$
	\dk\Big(\f{H_{\cap}^n-\E[H_{\cap}^n]}{\sqrt{\Var(H_{\cap}^n)}}, N(0,1) \Big) \in O\Big(\f{(\log |Q_n|)^{2d}}{\sqrt{|Q_n|}}\Big).
	$$
\end{corollary}	

{Similar convergence rates were obtained in \cite[Thm. 1.2]{gibbsCLT}, but for a smaller class of Gibbs processes.  	
	As mentioned in Section \ref{ssec:gibbs},} very recently a quantitative CLT was derived in \cite{CX22}  under the EDD conditions in  without the need for any monotonicity assumptions on the score function $g$. However, the convergence rates are {slightly worse and} given in terms of the Wasserstein distance, whereas bounds for the Kolmogorov distance are often easier to interpret. We note that also in other contexts, one frequently encounters the situation where extending Wasserstein convergence rates to Kolmogorov convergence rates requires substantial additional work \cite{kolm}. {In the special case of Poisson processes, better convergence rates were obtained in \cite{mehler}.}

{The weak stabilization condition \eqref{eq:as_stab}  is weaker than the exponential stabilization condition \eqref{eq:s2cap} when the PI is translation-invariant. We briefly explain why.} Fix $l \ge 1$. For $i \ge 2l$, we let
$$E_{i, l}^{(1)} := \Big\{\max_{x\in \X \cap Q_i}\max_{w \in \Z^d}\max_{m \in \N}R(x, \X \cap Q_{w, m} ) \le i/2 \Big\}$$
denote the event that the stabilization radius of the points in $Q_i$ is at most $i/2$ for any of the windows $Q_{w, m}$ with $w\in \Z^d$ and $m \ge 1$. We note that for fixed $x$, the number of different sets of the form $B_{i/2}(x) \cap Q_{w, m}$ is finite and of order $O(i^d)$. Hence, by exponential stabilization and the GNZ formula \eqref{eGNZ}, the probabilities $1 - \P(E_{i, l}^{(1)})$ decay to 0 exponentially in $i$. Hence, by the Borel-Cantelli lemma, there exists an almost surely finite random variable $I_0$ such that $E_{i, l}^{(1)}$ occurs for all $i \ge I_0$. We now argue that this implies the weak stabilization. Indeed, let $m_n \ge 1$ and $w_n \in \Z^d$ be such that $\R^d=\bigcup_{k\ge 1} \bigcap_{n\ge k} Q_{w_n,m_n}$. Then, $g(x, \X \cap Q_{w_n,m_n}) = g(x, \X\cap (Q_{w_n,m_n}\sm Q_l))$ holds for every $n \ge 1$   and every $x \in (\X \sm Q_{I_0}) \cap Q_{w_n, m_n}$. Therefore, whenever $n$ is such that $Q_{w_n, m_n} \supseteq Q_{2I_0}$,
\begin{align*}
	& H(\X\cap Q_{w_n,m_n}) - H(\X\cap (Q_{w_n,m_n}\sm Q_l))\\
	& = \sum_{x \in \XX \cap Q_l} g(x, \X \cap B_{I_0}(x)) + \sum_{x \in \X \cap (Q_{I_0}\sm Q_l)}(g(x, \X \cap Q_{2I_0}) -g(x, \X \cap (Q_{w_n,m_n} \sm Q_l))).
\end{align*}
Noting that the last term stabilizes when $n\to \infty$ if the stabilization radius $R(x, \X  \sm Q_l)$ is finite proves the weak stabilization condition.

To prove Theorems \ref{thm:1} and \ref{thm:2}, we extend to Gibbs processes the martingale approach from \cite{penrose} and the Palm couplings from \cite{chen}. The proofs rely on a precise control of the decay of spatial correlations through refined forms of disagreement coupling given in Section \ref{sec:const} that are tailored to Theorems \ref{thm:1} and \ref{thm:2}.

	%
	%
\section{Examples}
\label{sec:exa}
In this section, we present \cb{specific examples for Theorems \ref{thm:1} (persistent Betti numbers; minimal spanning tree) and \ref{thm:2} (total edge length in Voronoi tessellations and $k$-nearest neighbor graphs). In all examples, we consider Gibbs point processes where Assumption \eqref{as:(A)} holds.}

\subsection{Example for Theorem \ref{thm:1}: Persistent Betti numbers}
\label{sec:bet}
 The persistent Betti numbers are invariants used in topological data analysis for summarizing the persistence diagram. They are defined for $r\le s $, $q=0,1,\ldots, d-1$ and a finite set $\vp \in \Nlf$ by
	\begin{equation*}
	 \beta^{r,s}_q{(\vp)} = \dim \left(\Im (H_q( B_r(\vp);\Z/2\Z)\to H_q(B_s(\vp);\Z/2\Z) ) \right),
	\end{equation*}
 where $H_q(\cdot,\Z/2\Z)$ denotes the $q$th homology group with coefficients in $\Z/2\Z$. \cb{Homology is a fundamental quantity in algebraic topology.  Loosely speaking, the rank of $H_q$ encodes the number of $q$-dimensional holes of a topological space, see \cite{edHar} for a thorough introduction to topological data analysis.}
	 The use of persistent Betti numbers for statistics of point patterns was first suggested in \cite{robins} under the name rank functions. 

	Most existing CLTs for functionals on point processes, e.g. \cite{yogesh,gibbsCLT}, assume the functional is written in terms of an exponentially stabilizing score function. This is not easily verified for $\beta^{r,s}_q$. \cb{However, some first results on consistency and asymptotic normality for Betti numbers were given when $r=s$  in \cite{b2,b1,yogeshAdler2}. These results were extended to general $r$ and $s$ in \cite{shirai,krebs3}. All asymptotic normality results previously considered have been restricted to Poisson or Bernoulli point processes.  Existing  CLTs were derived by observing that weak stabilization in the sense of \cite{penrose} holds for persistent Betti numbers. The same observation leads to the following corollary of Theorem \ref{thm:1}. }

	\begin{corollary}\label{cor:betti}
		Let $\XX$ be an infinite-volume Gibbs point process
		with translation-invariant PI satisfying \eqref{as:(A)}. Then for all $0\leq r\leq s$
		\begin{equation*}
			|Q_n|^{-1/2}(\beta_q^{r,s}(\X \cap  Q_n) - \E\beta_q^{r,s}(\X \cap  Q_n)  ) \to N(0,\s^2)
		\end{equation*}
		for some $\s^2 \ge 0$. If $\k>0$, then $\s^2>0$.
	\end{corollary}

We show the first statement of Corollary \ref{cor:betti} below. The proof of \cb{the} positivity of the limiting variance is deferred to the end of Section \ref{sec:clt}. The proof relies on the fact that $\beta_q^{r,s}$ can be computed from the \v{C}ech complex via the Nerve Theorem, see \cite[Chap. III.2]{edHar} for details.

	\bep
		The corollary follows from Theorem \ref{thm:1} if we can show Condition \eqref{eq:m1} and  \eqref{eq:s1}. Condition \eqref{eq:s1} follows directly from \cite[Lem 5.3]{shirai}.
		For Condition \eqref{eq:m1}, recall that $k+1$ points in a finite point pattern $\vp$ form a $k$-simplex of filtration value $t$ in the \v{C}ech complex if and only if the smallest ball containing the $k+1$ points has radius $t$.  Thus, adding a point $y$ to a point pattern $\vp$ changes the $k$-simplices in the \v{C}ech complex with filtration value less than or equal $s$ only by adding and removing simplices with all $k+1$ vertices contained in a ball of radius $2s$ around $y$. By \cite[Lem. 2.11]{shirai}, this changes $\beta_q^{r,s}$  by at most the number of added and removed $q$- and $(q+1)$-simplices. This number is bounded by the total number of $(q+1)$- and $(q+2)$-tuples of points in $\vp\cap B_{2s}(y)$, which is again bounded by
		\begin{equation*}
			\big(\vp(B_{2s}(y)) + 1\big)^ {q+1} +  \big(\vp(B_{2s}(y)) + 1\big)^{q + 2} \le 2 \big(\vp(B_{2s}(y)) + 1\big)^{q+2}.
		\end{equation*}
		This shows Condition \eqref{eq:m1}. 
		\enp
\subsection{Example for Theorem \ref{thm:1}: Total edge length in minimal spanning trees}
\label{sec:mst}

 For a finite set  $\vp \su \R^d$, the minimal spanning tree (MST) is a tree defined on the vertex set $\vp$. The edge set of this tree is chosen such that it minimizes the total edge length.  We note that this tree is almost surely unique for any point process that {satisfies \eqref{eq:DLR2}. We let $H(\vp)$ denote the total edge length of the MST. For Poisson point processes, the asymptotic properties of $H$} were considered in \cite{kesten} and \cite[Theorem 2.3]{wlln}. We discuss this example because it is one of the prototypical examples where exponential stabilization fails even in the Poisson setting. Hence, we are not able to rely on Theorem \ref{thm:2}. To verify \eqref{eq:moment_cond} and \eqref{eq:as_stab}, we proceed as in \cite{kesten} with only few modifications. To make the presentation self-contained we give the details below. 
\medskip

%
%
%
%
\noindent {\bf Moment bound.} 
For any bounded Borel set $A \su \R^d$ and $\vp\in \Nlf$, we have a trivial upper bound
\begin{align}
	\label{eq:mstb}
	H(\vp \cap A) \le \ms{diam}(A)\vp(A).
\end{align}

We bound the positive and negative parts of $H_n( \X) - H_n( \X \sm Q_{z,m})$ separately.  
For the positive part, we first note that if $\X\cap Q_n \su Q_{z, m}$, then by \eqref{eq:mstb},
$$H_n(\X) - H_n(\X \sm Q_{z,m}) \le \sqrt d m \X(Q_{z, m}),$$
and by the stochastic domination \eqref{as:(A)}, the fifth moment of the right-hand side is of order $O(m^{5 + 5d})$. Thus, we may henceforth assume that $K_{z, m}\le m_0(n)$, where we write 
$$K_{z, m}:= \min\big\{ k \ge m \co \X \cap Q_{z, k} \not\su Q_{z,m}\big\}$$
and
$$m_0(n) := \min\big\{k \ge m\co  Q_n \su Q_{z, k}\big\}.$$
Now, we derive an upper bound for $H_n(\X)$. For this, we first build a minimal spanning tree on $\X \cap (Q_n \sm Q_{z, m})$, then build a spanning tree on $\X \cap Q_{z, m}$ and finally connect the two trees through a point in $\X\cap Q_{z, K_{z, m}}$ to build a spanning tree on $\XX\cap Q_n$. Therefore, 
\begin{align}
	\label{eq:mstp}
	H_n( \X) \le H_n( \X \sm Q_{z,m}) + \sqrt d K_{z, m} \X( Q_{z, m}).
\end{align}
As noted above, the tenth moment of $\X( Q_{z, m})$ is of order $O(m^{10d})$. Hence, by the Cauchy-Schwarz inequality, it suffices to show that the tenth moment of $K_{z, m}$ grows polynomially in $m$. To that end, we note that
\begin{align}
	\label{eq:zmn}
	\E[K_{z, m}^{10}] \le \sum_{k \ge m} k^{10}\P(K_{z, m} = k)
\end{align}
Now, by the Poisson-likeness derived in Proposition \ref{pr:like}, we have $\sup_{n \ge 1}\sup_{z\in \Z^d}\P(K_{z, m} = k) \le c_1(k^d - m^d) \exp\big(-(k^d - m^d)/c_1\big)$. Hence, the right-hand side of \eqref{eq:zmn} is of order $O(m^{d + 10})$, as asserted.

%
%
Finally, we bound the negative part of the difference $H_n( \X) - H_n( \X \sm Q_{z,m})$. {First, if $\X \cap Q_n \su Q_{z, m}$, then $H_n( \X) - H_n( \X \sm Q_{z,m})$ is positive. If $\X\cap Q_n \su Q_n\sm Q_{z,m}$ then $H_n( \X) = H_n( \X \sm Q_{z,m})$.} Hence, we may assume that in the MST on $\X\cap Q_n$ there are some edges crossing  $\pa Q_{z, m}$. We now bound the length of the MST on $\X \cap (Q_n \sm Q_{z, m})$ where as a starting point, we consider the  MST on $\X \cap Q_n$. Let 
$$\mc L_{z, m, n} := \Big\{\text{ edges in the MST of $\X\cap Q_n$ with exactly one endpoint in $\X \cap Q_{z, m}$}\Big\}.$$
{Then, we get a graph connecting $\X\cap (Q_n\sm Q_{z,m})$ by removing all edges with at least one endpoint in $Q_{z,m}$ and connecting  all points in $\X\cap (Q_n\sm Q_{z,m})$ that used to be connected to $\X \cap Q_{z,m}$ to any fixed  other such point. The length of any such new connection is bounded by $\sqrt d m $ plus the length of two edges from $\mc L_{z, m, n}$.
}
Henceforth, we use that there is a deterministic upper bound $d_{\ms{max}}\ge 1$ for the possible degree of any vertex in the MST.
This yields  
\begin{align}
	\label{eq:mstn}
	H_n( \X \sm Q_{z,m}) - H_n( \X)\le     \sum_{e \in \mc L_{z, m, n}} (2L_{z, m, n}  + \sqrt dm) \le d_{\ms{max}}\X( Q_{z, m})(2L_{z, m ,n} + \sqrt d m),
\end{align}
where 
$L_{z, m, n} := \max_{e \in \mc L_{z, m, n}}|e|$
denotes the maximal edge length in $\mc L_{z, m, n}$. Hence, as in the case for the positive part, by the Cauchy-Schwarz inequality, it suffices to show that $\E[L_{z,m,n}^{10}]$ is at most of polynomial order in $m$. To achieve this goal, we note that 
\begin{align}
	\label{eq:lzmn}
	\E[L_{z, m, n}^{10}] \le \sum_{\ell \ge 0}(\ell + 1) \P(L_{z, m ,n}^{10} \ge  \ell ).
\end{align}
Now, under the event $\{L_{z, m, n} \ge \ell^{1/10}\}$, there exists an edge $e\in\mc L_{z, m, n}$ of length at least $\ell^{1/10}$. {Since the MST is a subset of the Delaunay complex, there is a ball with diameter at least $|e|$ containing $e$ that does not contain any points of $\X \cap Q_n $ in its interior.} But again by the Poisson-likeness, the probability of this event can be bounded by $c_1 (m +\ell^{1/10})^d \exp(-\ell^{d/10}/c_1)$ for some constant $c_1 > 0$. Therefore, the right-hand side of \eqref{eq:lzmn} is of order $O(m^d)$, as asserted.
\bigskip

\noindent {\bf Weak stabilization.}
To verify the weak stabilization condition, we proceed as in \cite[Proposition 3]{kesten}. To make this precise, we first need the notion of separating sets. More precisely, for $\vp \in \Nlf$, a \emph{separating set of width $\de > 0$} for a box $Q_l$ is a finite set $S \su \vp \cap (Q_{l + \de  + 18r_0} \sm Q_{l + 18r_0})$ such that any line segment from a point in $\pa Q_{l + 18r_0}$ to a point in $\pa Q_{l + \de + 18r_0}$ passes within distance $3r_0$ of some point of $S$. {We note our definition is scaled by $9r_0$ compared to the original definition in \cite{kesten}, which uses 1 instead of $9r_0$ in the definition of the windows and $1/3$ instead of $3r_0$ for the distance.  After appropriate rescaling, the arguments in \cite{kesten} extend to our setting.
We shall need the following two results from  \cite{kesten}, which we restate in our setting.}

%
%
\bepr[Monotonicity {in increasing windows}; Proposition 3 of \cite{kesten}]
\label{pr:mon}
Let $l \ge 1$ {and $\vp \in \Nlf$}.  Assume that $S\su {\vp} \cap (Q_{l + \de + 18r_0} \sm Q_{l + 18r_0})$ is a separating set of width $\de > 2r_0$ for $Q_l$. Then, for any $w, w' \in \Z^d$, $m, m' \ge 1$ with $Q_{l + \de + 18r_0} \su Q_{w, m} \su Q_{w', m'}$, we have 
$$ H({\vp}\cap Q_{w,m})-H({\vp}\cap (Q_{w,m}\sm Q_l)) \le H({\vp}\cap Q_{w',m'})-H({\vp}\cap (Q_{w',m'}\sm Q_l)).$$
\enpr

The value of Proposition \ref{pr:mon} becomes apparent when combined with \eqref{eq:mstp} since the latter shows that almost surely the sequence $H(\XX\cap Q_{w_n,m_n})-H(\XX\cap (Q_{w_n,m_n}\sm Q_l))$ remains bounded as $n \tff$. Hence, the almost sure convergence in \eqref{eq:as_stab} holds as soon as we can show that almost surely $\XX$ has a $\de$-separating set for $Q_l$ for some $\de > 0$. 

%
%
To prove this final step, by the Borel-Cantelli lemma, it remains to show that the probability that there fails to be a $\de$-separating set decays at exponential speed in $\de$. To this end, we choose a set of points $\phi \su \pa Q_{l + 2r_0}$ such that every point in $\pa Q_{l + 2r_0}$ is at distance at most $r_0$ from a point in $\phi$. This is possible with a cardinality $\#\phi \in O(l^{d - 1})$. Similarly, we choose a subset $\phi'$ with the same property in $\pa Q_{l + \de + 2r_0}$, where now $\#\phi' \in O((l + \de)^{d - 1})$. In particular, if for every line segment between $\phi$ and $\phi'$ there is a point of $S$ within distance $2r_0$, then for every line segment between $\pa Q_{l + 2r_0}$ and $\pa Q_{l + \de +2r_0}$ there is a point of $S$ within distance $3r_0$. 

 Finally, the Poisson-likeness from Proposition \ref{pr:like} shows that for fixed $P \in \pa Q_{l + 2r_0}$ and $P' \in Q_{l + \de + 2r_0}$, the probability $\P\big(\X \cap B_{2r_0}([P, P']) =\es\big)$ decays exponentially in $\de$, { where $[P,P']$ denotes the line segment from $P$ to $P'$.} Since $\#\phi$ and $\#\phi'$ were of order $O(l^{d-1})$ and $O((l + \de)^{d - 1})$, a union bound concludes the proof.

%
%
\subsection{Example for Theorem \ref{thm:2}: Total edge length of Gibbs-Voronoi tessellations}
		\label{ssec:vor}
 For a locally finite $\vp \su \R^2$ we define the {\em Voronoi cell} $C(x,\vp)$ of $x \in \vp$ as the set of all points in $\R^2$ whose Euclidean distance to $x$ is less than or equal to its distance to all other points of $\vp$. The system $\{C(x,\vp)\}_{x \in \vp}$ is called the {\em Voronoi tessellation}. We define $g(x,\vp)$ as one half of the total {boundary length} of the cell $C(x,\vp)$. We use the convention that $g(x,\vp) = \ff$ if the cell $C(x, \vp)$ is unbounded. Note that this convention is only made for completeness and does not play any practical role for the functionals considered since by stationarity, all cells $C(x, \X)$ are bounded with probability 1.
 Then, the total edge length of the Voronoi tessellation induced by $\vp$ with \cb{generators} in $Q_n$ is (up to boundary effects) given by
		\begin{align*}
			H^n := \sum_{x \in \vp \cap Q_n} g(x,\vp).
		\end{align*}
As observed in Section \ref{ssec:gibbs}, Gibbs processes satisfying Assumption \eqref{as:(A)} are Poisson-like. Hence, one can proceed as in  \cite[Section 1 and 5.2]{gibbs_limit} to obtain the bounded moment condition \eqref{eq:m2} and the exponential stabilization property \eqref{eq:s2}. Moreover, \cite[Theorem 2.3]{gibbsCLT} discusses further conditions under which the variance lower bound can be verified.
For the condition $g(x,\vp\cup \{y\})\le g(x,\vp)$, we note that $C(x,\vp \cup \{y\})\su C(x,\vp)$ are convex subsets of $\R^2$. Hence, it follows from Theorem \ref{thm:2} that \cb{$H^n$} satisfies a CLT. 

This extends \cite[Theorem 4.13]{chen} to Gibbs processes satisfying \eqref{as:(A)}. { We also note that if we are not interested in rates, then this example follows from \cite[Theorem 1.13]{yogesh}. Indeed, we first recall from Corollary \ref{cor:EDD} that under Assumption \eqref{as:(A)} the Gibbs process $\XX$ is EDD. Moreover, as mentioned in \cite[Section 1]{CX22}, the EDD property is stronger than the fast decay of correlation formulated in \cite[Equation (1.10)]{yogesh}. Together with the aforementioned exponential stabilization (also  proved in \cite{benes}), this implies that the considered score function and point process form an admissible pair in the sense of \cite[Definition 1.7]{yogesh}. Finally, combining this observation with the aforementioned upper moment bound and lower variance bound, we see that \cite[Theorem 1.13]{yogesh} gives the asymptotic normality without rates.}

\subsection{Example for Theorem \ref{thm:2}: Edge lengths in $k$-nearest neighbor graph}
\label{ssec:knn}
\subsubsection{Total edge length}
We discuss the $k$-nearest neighbor (kNN) example from \cite[Theorem 4.1]{CX22}. In the kNN graph on a locally finite set $\vp \su \R^d$, we put an edge between $x, y \in \vp$ if $x$ is one of the $k$-nearest neighbors of $y$ or vice versa. Then, we define $g(x, \vp)$ as one half of the total length of all kNN edges incident to $x$. Again, we use the convention that $g(x, \vp) = \infty$ if $\vp$ contains fewer than $k + 1$ points. As in Example \ref{ssec:vor}, the quantity of interest is then given by
		\begin{align*}
			H^n := \sum_{x \in \vp \cap Q_n} g(x,\vp).
		\end{align*}
As in Example \ref{ssec:vor}, the derivation of the bounded moment condition \eqref{eq:m2} and of the exponential stabilization property \eqref{eq:s2} from the Poisson-likeness follows the existing approaches in literature, e.g. \cite[Theorem 5.2]{gibbs_limit}, \cite[Theorem 4.1]{CX22}. For examples of Gibbs point processes where the lower bound is satisfied, we refer to \cite[Theorem 2.2]{gibbsCLT}. Finally, we note that by the definition,  $g(x, \vp) = 0$ can never happen, thereby verifying the final condition from Theorem \ref{thm:2}. Hence, we deduce that \cb{$H^n$} satisfies a CLT. 

\subsubsection{Number of large $k$-nearest neighbor edges}
The extension of the total edge length example from $H^n$ to $H^n_\cap$ is more delicate, since there is a positive probability that $\X \cap Q_n$ consists of fewer than $k + 1$ points.  An example where the condition is satisfied  for $H^n_\cap$  is when $g(x, \vp)$ is the indicator that the distance to the $k$th nearest neighbor of $x$ in $\vp$ is at least some value $a > 0$.

	%
%

\section{Disagreement couplings}
\label{sec:const}

Disagreement coupling is a technique for constructing Gibbs point processes in bounded windows by thinning a dominating Poisson process. We present the algorithms we use in Section \ref{ss:fw}. The original and most important use of disagreement coupling is to control the amount of disagreement between Gibbs processes with different boundary conditions by relating them to the connected components of an associated Boolean model. Another important application in the context of normal approximation is to bound the spatial correlations of the Gibbs process. We show some  spatial homogeneity and decorrelation results in Section \ref{ss:hd}.  In Section \ref{ss:rfin}, we show how  disagreement coupling can be used to construct couplings between Gibbs processes in  bounded domains and the unbounded domain with strong control over their disagreement. This is a main ingredient in  the proof of Theorem \ref{thm:1}.
For the proof of Theorem \ref{thm:2}, we need couplings between  a whole family Gibbs processes obtained by local perturbations of the PI with bounds on the pairwise disagreement.  Such  couplings are given in  Section \ref{ss:rfin}. The result is also be essential for the proof of Theorem \ref{thm:1}, where it allows us to consider couplings in unbounded domains with locally differing boundary conditions.

The original definition of disagreement coupling in \cite{dp} applied to Gibbs processes on standard Borel spaces, thereby allowing applications to particle processes \cite{benes} and Gibbs processes with density with respect to random connection models \cite{betsch}. However, since all our applications of disagreement coupling are to processes on Euclidean space, and to keep the exposition simple, we restrict to the Euclidean case.

%
%

%
%
\subsection{Thinning algorithms}
\label{ss:fw}
%
%

In this section, we review the construction from \cite{dp} of the finite-volume Gibbs process $\XX(Q,\psi)$ on a bounded domain $Q$ with boundary conditions $\psi \in \Nlf_{Q^c}$ as a thinning of a Poisson process. In Section \ref{sss:spe}, we first review the basic version of the algorithm introduced in \cite[Sec. 5]{dp}. This is the building block for more refined algorithms in Section \ref{sss:dc}. We give a general recipe for constructing thinning algorithms with more flexible properties generalizing the one suggested in \cite[Sec. 6]{dp}. Moreover, we give two examples of such algorithms which are designed to have the properties, we need for the proofs of Theorems \ref{thm:1} and \ref{thm:2}.

Henceforth, we let $\Pds$ be a homogeneous unit-intensity Poisson point process on $\R^d \ti [0, \a_0]$ and write $\Pd$ for the projection of $\Pds$ onto the $\R^d$ coordinate, which is a Poisson process on $\R^d$ with intensity $\a_0$.  More generally, we use $\vp^*$ to denote elements of $\Nlf_{\R^d \ti [0, \a_0]}$ and $\vp$ denotes its projection to $\R^d$. We let $\vp^*_Q$ and $\vp_Q$ be the restriction of $\vp^*$ and $\vp$ to $Q\times [0,\a_0]$ and $Q$, respectively. We refer to the points of $\vp$ in a connected component of $B_r(\vp)$ as a cluster.

%
%
\subsubsection{Standard Poisson embedding}
\label{sss:spe}

Consider a bounded Borel set $Q \su \R^d$ and a set of boundary conditions $\psi \in \Nlf_{Q^c}$. In this section, we introduce \emph{the standard thinning operator}  $T_{Q, \psi}: \Nlf_{Q\times [0,\a_0]} \to \Nlf_Q$. For $\vp^* \in  \Nlf_{Q\times [0,\a_0]}$, the set $T_{Q, \psi}(\vp^*)$ is a thinning of $\vp$, that is $T_{Q, \psi}(\vp^*)\su \vp$. The $(d+1)$st coordinate of $\vp^* $ is used to decide whether or not to keep the point in the thinning. The key property of $T_{Q, \psi}$ is that the thinning $T_{Q, \psi}(\Pds_Q)\su \Pd_Q$ of the Poisson process $\Pds_Q$ has the distribution of the Gibbs process $\XX(Q, \psi)$. We call $T_{Q, \psi}(\Pds)$ the \emph{standard Poisson embedding}, since it provides an embedding $\XX(Q, \psi)\su \Pd$.

To define $T_{Q, \psi}$, we need an injective measurable map  $\iota\co Q \to \R$. Such a map is guaranteed to exist since $\R^d$ is a standard Borel space. We sometimes write $T_{Q, \psi, \iota}$ to emphasize that $T_{Q, \psi}$ depends on $\iota$. The map $\iota$ defines a total ordering $\le_{\iota}$ on $Q$. We sometimes refer to the map $\iota$ itself as an ordering.
The main idea to construct $T_{Q, \psi}(\vp^*)$ is to go through the points of $\vp$ in the order $\le_{\iota}$ and decide for each whether to keep it in the thinned process. For  $x \in  Q $, define 
\begin{align*}
	Q_{(-\ff,x)} {}:=\{ y \in Q \co y <_\iota x\} \, \text{ and }\,Q_{(x, \ff)} {}:=\{ y \in Q \co y >_\iota x\}.
\end{align*}
For $x \in  Q$ and $\psi' \su (\R^d\sm (Q_{(x, \ff)} \cup \{x\}))\times \R$, define the retention threshold
\begin{align}
	p(x, Q, \psi') := \k(x, \psi') \f{Z_{Q_{(x, \ff)}}\big(\psi' \cup \{x\}\big)}{Z_{Q_{(x, \ff)}}(\psi')}.\label{def:p}
\end{align}
It is shown in \cite[Lem. 4.1]{dp} that $0\le p(x, Q, \psi') \le \a_0$.

The algorithm for constructing $T_{Q, \psi}(\vp^*)=\{x_1,x_2,\ldots\}$ goes as follows:
\begin{itemize}
	\item[$\bullet$] To define the first point $x_1 := x_1(\vp^*)$ of $T_{Q, \psi}(\vp^*)$, we consider the set 
	$$\big\{x\in Q \co (x, u) \in \vp^*  \text{ and } u \le p(x, Q, \psi) \big\}$$
	for which $u$ is below the retention threshold. We choose $x_1$ to be the smallest element in this set with respect to the $\le_\iota$-order. If there is no such element, then the construction terminates with $T_{Q, \psi}(\vp^*) = \es$. 
	\item[$\bullet$] Otherwise, we proceed inductively. Suppose, we have defined the points $x_1,\ldots,x_k $. This will be the restriction of $T_{Q, \psi}(\vp^*)$ to $Q_{(- \ff,x_k)}$. It remains to define $T_{Q, \psi}(\vp^*)$ on $Q_{(x_k, \ff)} $. We define $x_{k+1}=x_{k + 1}(\vp^*)$ to be the $\iota$-smallest element of the set 
	\begin{equation}\label{eq:admissible}
		\big\{x\in Q_{(x_k, \ff)}  \co (x, u) \in \vp^*  \text{ and } u \le p(x, Q_{(x, \ff)}, \psi \cup \{x_1, \dots, x_k\})\big\}.
	\end{equation}
	\item[$\bullet$] The construction terminates once the  set \eqref{eq:admissible} is empty. 
\end{itemize}

We state the key result of \cite[Thm. 5.1]{dp}, namely that the thinned Poisson process $T_{Q, \psi}(\Pds) $ has the distribution of the Gibbs process $\XX(Q, \psi)$, for later reference. 

%
%
\bepr[Correctness of the Poisson embedding \cite{dp}]
\label{pr:emb}
Assume that the PI satisfies $ \k \le \a_0$. Let $Q \su \R^d$ be bounded Borel and $\psi \in \Nlf_{Q^c}$ be finite. Furthermore, let $\Pds \su \R^d \ti [0, \a_0]$ be a unit-intensity Poisson point process. Then, $T_{Q, \psi}(\Pds_Q)\su \Pd_Q$ is distributed as the Gibbs process $\XX(Q, \psi)$. 
\enpr
\bep
See \cite[Thm. 5.1]{dp}.
\enp

We remark for later that the standard thinning algorithm has the property that if $B\su Q$ such that $\iota(x)< \iota(y)$ for all $x\in B$, $y\in Q\sm B$, and $\xi = T_{Q, \psi}(\vp^*) \cap B$, then
\begin{equation}\label{eq:stand_emb_property}
	T_{Q, \psi}(\vp^*) =\xi \cup T_{Q\sm B, \psi\cup \xi}(\vp^*_{Q\sm B}).
\end{equation}

While an injective measurable map $\iota:Q \to \R$ is guaranteed to exist, it is not so intuitive what such a map looks like. While usually not strictly necessary, it is often convenient to replace $\iota$ by a more familiar map. 
Thus, let $\iota:Q \to \R$ be a measurable map such that
\begin{equation*}
	\big| \iota^{-1}(r) \big|= 0
\end{equation*}
for all $r\in \R$. We define a partial ordering of $\R^d$ by $x <_\iota y$ if $\iota(x) <\iota (y)$. If $\big| \iota^{-1}(r) \big|= 0$ for all $r\in \R$, this almost surely induces a well-defined total ordering of the Poisson process and hence $T_{Q,\psi, \iota}(\Pds)$ becomes well defined.

\begin{proposition}\label{pr:iota}
	Assume that $\iota: Q \to \R$ is a  measurable map such that for all $r\in \R$,
	\begin{equation}\label{eq:iota_cond}
		\big| \iota^{-1}(r) \big|= 0.
	\end{equation}
	Then, $T_{Q, \psi,\iota}(\Pds_Q)$ is distributed as $\XX(Q, \psi)$.  
\end{proposition}

The proof is deferred to the appendix. Using Proposition \ref{pr:iota}, one may order $Q$ e.g.\ by the distance to the origin.

%
%
\subsubsection{Disagreement couplings}\label{sss:dc}

The standard Poisson embedding defines a coupling between $\XX(Q,\psi)$ and $\Pd$ for any choice of boundary conditions $\psi$. However, it does not provide much control over how $T_{Q, \psi}(\Pds)$ are related for different choices of $\psi$. For this reason, a refined thinning algorithm called \emph{disagreement coupling} was suggested in \cite{dp} . This algorithm makes it possible to control the amount of disagreement between Gibbs processes with different boundary conditions. However, as we shall need some extra properties of the algorithms compared to \cite{dp}, we first give a general recipe for constructing disagreement couplings. We then give two examples of such algorithms. 

To define the algorithm, we need the notion of a \emph{stopping set}. A stopping set $S$ on a Borel set $B\su \R^d$ assigns to each  $\vp \in \Nlf_B$ a Borel set  $S(\vp)\su B$. We require that the map $(x,\vp) \mapsto \one_{\{x \in S(\vp)\}}$ is measurable.  Moreover, $S$ must have the property that $S(\vp)$ depends only on $\vp\cap S(\vp)$, i.e.\  $S(\vp) =  S(\vp\cap S(\vp)) $ for all $\vp \in \Nlf_B$. We write $S=S(\vp)$ when $\vp$ is clear from the context. A key property of a stopping set is the following, see e.g. \cite[Thm. A.3]{stoppingset}: If $S$ is a stopping set on $B$ and $\PP'$ is a Poisson processes on $B$, then 
\begin{equation}\label{eq:stopping}
	\text{conditionally on } \PP'_{S(\PP')}, \,	\PP'_{B\sm S(\PP')} \text{ is distributed as a Poisson process on $B\sm S(\PP')$.}
\end{equation}

%
%
The disagreement coupling algorithm needs two ingredients:
\begin{itemize}
	\item[$\bullet$] An sequence $S_0^*, S_1^*, S_2^*, \dots $ of stopping sets on $Q\times [0,\a_0]$ of the form $S_n^*=S_n\times [0,\a_0]$. The sequence should increase towards $Q\times [0,\a_0]$ in the sense that for any $\vp^* \in \Nlf_{Q\times [0,\a_0]}$,  
	$$S_0^*(\vp^*)=\es\su S_1^*(\vp^*)\su S_2^*(\vp^*) \su \dots  \quad \text{ and } \quad \bigcup_n S_n^*(\vp^*) = Q\times [0,\a_0].$$ 
	\item[$\bullet$] A family of injective ordering maps  $\iota_n:=\iota_n^{\vp^*}:Q\sm S_n(\vp^*) \to \R $ for $n\ge 0$ and $\vp^*\in \Nlf_{Q\times [0,\a_0]}$ satisfying \eqref{eq:iota_cond} such that $\iota_n^{\vp^*}$ depends only on $\vp^*$ via $\vp^* \cap S_n(\vp^*)$. Moreover, we require
	\begin{equation}\label{eq:iota_n}
		\iota_n(x) < \iota_n(y) \text{ for all } x\in S_{n+1}\sm S_n, y \in Q\sm S_{n+1}. 
	\end{equation}
	To ensure measurability, we extend $\iota_n$ to a map $\tilde{\iota}_n :Q \times \Nlf_{Q\times [0,\a_0]} \to \R\cup\{-\ff\} $ by setting $\tilde{\iota}_n^{\vp^*}(x)=-\ff$ if $x\in S_n(\vp^*)$. We then require that $\tilde{\iota}_n: Q\times \Nlf_{Q\times [0,\a_0]} \to \R \cup\{-\ff\} $, is measurable.
\end{itemize}

Let $\psi \in \Nlf_{Q^c}$ be the boundary conditions. The disagreement coupling thinning operator $\tdc_{Q,\psi}: \Nlf_{Q\times [0,\a_0]} \to \Nlf_{Q} $ constructs the thinning of a point pattern $\vp^*$ inductively in $n$ by constructing the thinned point pattern $\xi_n:=\xi_n(\vp^*)$ on the set $S_n\sm S_{n-1}$ as follows: 
\begin{itemize}
	\item[1.] First, apply the standard thinning algorithm $T_{Q,\psi, \iota_0}$ on $Q$ and set  $\chi_1:=\chi_1(\vp^*):=T_{Q,\psi, \iota_1} (\vp^*)$. Let $\xi_1:=\chi_1\cap S_1$. 
	\item[2.] Assume inductively that we have constructed $\xi_i$ on $S_{i}\sm S_{i-1}$ for $i\le n$. The thinning of $\vp^*$ on $S_n$ is $\xi_1\cup\ldots \cup \xi_n$. We next construct the thinning on $S_{n+1}\sm S_n$. For this, we add $\xi_1\cup\ldots \cup \xi_n$ to the boundary conditions $\psi$ and apply the standard thinning operator on $Q\sm S_n$ with ordering $\iota_n$ to obtain  $ \chi_{n+1}:= \chi_{n+1}(\vp^*):= T_{Q\sm S_n,\psi \cup \xi_1 \cup \ldots \cup \xi_n,\iota_n} (\vp^*)$. We set $\xi_{n+1}:=\chi_{n+1} \cap (S_{n+1}\sm S_n)$. 
	\item[3.] Finally, we define $\tdc_{Q,\psi}(\vp^*) :=\bigcup_{n\ge 1} \xi_n $. 
\end{itemize}
Note that the algorithm corresponds to applying the standard thinning operator, but with the fixed ordering replaced by an ordering $\tilde{\iota}$ that depends on  $\vp^*$, where $\tilde{\iota}$ equals $\iota_n^{\vp^*}$ on $S_{n+1}(\vp^*)\sm S_n(\vp^*)$.  

To illustrate this abstract framework, we now present a specific example, where we want to construct a disagreement coupling such that changing the boundary conditions $\psi$ locally in a set $B \su Q^c$ does not change the Gibbs process too much. Before discussing all technical details, we first describe the general intuition. Loosely speaking, the stopping sets work their way through all the connected components of  $B_{r_0/2}(\vp) $ that intersect $B_{r_0/2}(B) $. Afterwards, they proceed through the remaining components of $B_{r_0/2}(\vp) $ one at a time(note that each component $\mc C$ of $B_{r_0/2}(\vp) $ is enlarged  to $B_{r_0/2}(\mc C) $ when defining the stopping sets). 
Now, we turn to the precise construction.
%
%
\begin{example}\label{ex:dc1}[Cluster-based disagreement coupling]
	
 Define $S_1:= B_{r_0}(B)$.    This is a deterministic set consisting of all points within distance $r_0$ from $B$. Choose an ordering map $\iota_0:Q\to [0,1)$ such that $\iota_0(x)<\iota_0(y)$ whenever $x\in S_1$ and $Q\sm S_1$. Suppose $S_n$  is defined.  If $B_{r_0}(\vp \cap S_n) \sm S_n \ne \es$, we take $S_{n+1} := S_n \cup B_{r_0}(\vp\cap S_n) $. Otherwise, let $S_{n+1}=S_n\cup Q_{(-\ff,x)}$ where $x$ is the   $\iota_0$-smallest point in $(Q\sm S_n)\cap \vp$. If no such point exists, we set $S_{n+1}=Q$.  To define the ordering $\iota_n: Q \sm S_n \to [0,1)$, we let 
	$$\iota_n=\tfrac1{2} \iota_0 \one{\{B_{r_0}(\vp \cap S_n) \sm S_n\}} + ( \tfrac1{2} \iota_0 + \tfrac1{2} ) \one{\{Q\sm (S_n\cup B_{r_0}(\vp \cap S_n)  ) \}}, $$ 
	which satisfies \eqref{eq:iota_n}.
	
Note that The order in which the components not intersecting $B_{r_0/2}(B) $ are visited is determined by the $\iota_0$-smallest point in the cluster.	
	
	The advantage of this construction is that whenever $ B_{r_0}(\vp\cap S_n) = S_n$, so that we jump to a new cluster, it means that $(\vp \cap S_n) \cup (\psi \cap B)$ has distance more than $r_0$ from $Q\sm S_n$. Hence, the retention probabilities \eqref{def:p} used in the remainder of the algorithm are not influenced by $(\vp \cap S_n) \cup (\psi \cap B)$. Therefore, the decision about the thinning of each cluster does not affect the thinning of the following clusters.  Moreover, the boundary points in $B$, $\psi \cap B$, can only influence the thinning of clusters corresponding to components of $ B_{r_0/2}(\vp)$ that intersect $B_{r_0/2}(B) $. Changing $\psi \cap B$, the Gibbs processes agree on all other clusters 
\end{example}

The following proposition shows that $\tdc_{Q,\psi}(\Pds_Q)$ is indeed a Gibbs process. Note that if the maps $\iota_n$ are not injective, but satisfy \eqref{eq:iota_cond}, then $\tdc_{Q,\psi}(\Pds_Q)$ is still well defined almost surely.

%
%
\bepr[Correctness of disagreement coupling]
\label{th:dc_prop}
Suppose the PI is bounded by $\a_0$. Let $Q \su \R^d$ be a bounded Borel set and let $\psi \in \Nlf_{Q^c}$. Furthermore, let $\Pds_Q \su Q \ti [0, \a_0]$ be a unit-intensity Poisson point process. Then, $\tdc_{Q,\psi}(\Pds_Q)$ has the distribution of the Gibbs process $\XX(Q, \psi)$. This also holds if the maps $\iota_n$ are not injective, but satisfy \eqref{eq:iota_cond}.
\enpr

The proof is similar to that of \cite[Thm. 6.3]{dp}, except that in \cite{dp}, $S_{n+1}^*$ was determined by $\Pds_{S_n^*}$, while we only require $S_{n+1}^*$ to be a stopping set. This only changes the argument in \eqref{eq:iota_n-1} below.

\begin{proof}
	We start by showing that $\xi_1 \cup \dots \cup \xi_n \cup \chi_{n+1} \sim \XX(Q, \psi)$ by induction in $n$, where $\sim$ means that the point processes have the same distribution. For $n=0$, $\chi_1=T_{Q,\psi, \iota_0} (\Pds_Q)\sim \XX(Q, \psi)$ by Proposition \ref{pr:emb}.
	Assume for induction that it has been shown that $\xi_1 \cup \dots \cup \xi_{n-1} \cup \chi_n \sim \XX(Q, \psi)$.
	Taking a bounded nonnegative function $f: \Nlf_{Q} \to [0,\infty)$, we must show
	\begin{equation}\label{eq:equal_exp}
		\E[f(\xi_1 \cup \dots \cup \xi_n \cup \chi_{n+1})]= \E[f(\XX(Q, \psi))].
	\end{equation}
	Let  $\iota_{n-1}'$ be the restriction of $\iota_{n-1} $ to $Q\sm S_n$. Note that by \eqref{eq:iota_n}, both the maps $\iota_n,\iota_{n-1}':Q\sm S_n\to \R$ and $\xi_1,\ldots,\xi_n$ are completely determined by $\Pds_{ S_n^*}$. Thus, conditionally on $\Pds_{S_n^*}$, 
	\begin{equation*}
		\chi_{n+1}=T_{Q\sm S_n,\psi\cap \xi_1 \cup \dots \cup \xi_n ,\iota_n}(\Pds_{Q \sm S_n^*}) \sim \XX(Q\sm S_n, \psi \cup \xi_1 \cup \dots \cup \xi_n) \sim T_{Q\sm S_n,\psi\cup \xi_1 \cup \dots \cup \xi_n ,\iota_{n-1}'}(\Pds_{Q \sm S_n^*})
	\end{equation*}
	by definition of $\chi_{n+1}$, the property \eqref{eq:stopping} of the stopping set $S_n$, and Proposition \ref{pr:emb}.  Therefore,
	\begin{align*}
		\E[f(\xi_1 \cup \dots \cup \xi_n \cup \chi_{n+1})] {}&= \E[\E[f(\xi_1 \cup \dots \cup \xi_n \cup T_{Q\sm S_n,\psi\cup \xi_1 \cup \dots \cup \xi_n ,\iota_{n-1}'}(\Pds_{Q\sm S_n^*}) )|\Pds_{ S_n^*} ] ].
	\end{align*}
	Moreover, given $ \Pds_{S_n^*}$, \eqref{eq:iota_n} and  \eqref{eq:stand_emb_property} implies
	\begin{equation}\label{eq:iota_n-1}
		\xi_n \cup T_{Q\sm S_n,\psi\cup \xi_1 \cup \dots \cup \xi_n ,\iota_{n-1}'}(\Pds_{Q \sm S_n^*})= \chi_n(\Pds_{Q  \sm S_{n-1}^*}).
	\end{equation}
	Hence, by the induction assumption,
	\begin{align*}
		\E[f(\xi_1 \cup \dots \cup \xi_n \cup \chi_{n+1}) ]{}= \E\big[\E[f(\xi_1 \cup \dots \cup \xi_{n-1} \cup \chi_n)|\Pds_{ S_n^*} ]\big]
		= \E[f( \XX(Q, \psi))].
	\end{align*}
This shows \eqref{eq:equal_exp}.
	
	To complete the proof, take again a bounded nonnegative function $f: \Nlf_{Q} \to [0,\infty)$. Since $\Pds_Q$ is almost surely finite, there is almost surely a random $N\ge 1$ such that $\tdc_{Q,\psi}(\vp^*) := \xi_1 \cup \dots \cup \xi_n \cup \chi_{n+1}$ for all $n\ge N$. Dominated convergence shows that 
	\begin{equation*}
		\E[f(\tdc_{Q,\psi}(\Pds_Q) )]= 	\E\big[f(\bigcup_{n\ge 1} \xi_n )\big] = \lim_{n\to \infty } \E[f(\xi_1 \cup \dots \cup \xi_n \cup \chi_{n+1})] =  \E[f(\XX(Q, \psi))],
	\end{equation*}
	which proves the theorem.
\end{proof}

The algorithm in Example \ref{ex:dc1} suffices for most purposes in this paper. However, it has the disadvantage that $S_{n + 1}\sm S_n$ can be very large, thereby causing issues for deriving the total variation bounds needed for the later arguments, see Theorem \ref{th:disdis}. Therefore, we now discuss a refinement where $S_{n+1}\sm S_n $ contains at most one point from $\vp$ and $\big|S_{n+1}\sm S_n\big|\le \big|B_{r_0}\big|$. 

This is achieved by an algorithm $\Tra_{Q,B, \psi}$ that we refer to as \emph{radial coupling}. We note that the construction of this radial coupling below is rather involved. However, for most applications later, we shall not need the precise specification of the algorithm, but only its main properties. Therefore, we now state them here for easy reference. 

%
%
\bepr\label{pr:rad_prop}
Suppose the PI is bounded by $\a_0$. Let $Q \su \R^d$ be bounded Borel and $B\su  Q^c$. Furthermore, let $\Pds $ be a unit-intensity Poisson point process on $ \R^d \ti [0, \a_0]$. Then, the radial coupling has the property that if $\psi,\psi'\in \Nlf_{Q^c}$ are locally finite sets differing only on $B$, then $\Tra_{Q,B, \psi}(\Pds_Q)$ and $\Tra_{Q,B, \psi'}(\Pds_Q)$ differ only on the components of $B_{r_0/2}(\Pd_Q) $ intersecting $B_{r_0/2}(B)$. Moreover, each $S_{n+1}\sm S_n$ contains at most one point from $\Pd$, and $\big|S_{n+1}\sm S_n\big| \le \big|B_{r_0}\big|$.
\enpr

%
%
\begin{example}[Radial coupling]
	\label{ex:dc2}
	The idea of this algorithm  is similar to the one given in Example \ref{ex:dc1} to go through the clusters of $B_{r_0/2}(\vp) $ one at a time, but the description is more involved.
	We start again from a domain $Q$, boundary conditions $\psi\in \Nlf_{Q^c}$, and a set $B\su Q^c$ where we want to be able to change the boundary conditions $\psi \cap B$.  We choose a fixed ordering $\iota: Q \to [0,1)$ satisfying \eqref{eq:iota_cond} and a deterministic locally finite set $D \su Q$  with the coverage property that $Q \su B_{r_0}(D) $. We also fix a total ordering $\preceq$ on $D$. 
	
	The stopping steps $S_{i,j}$ of the algorithm are indexed by a pair of lexicographically ordered integers $(i,j)$, $i=1,\dots,I$, $j=1,\dots,N_i$, where $I:=I(\vp^*)$ and $N_i:=N_i(\vp^*)$. 
	The stopping set $S_{i,j}$ is given at the beginning of step $(i, j)$ of the algorithm. In this step, we  define three auxiliary subsets of $Q$ that depend on $\vp^*$, namely,
	\been
	\im ${v_{i,j+1}} \in D \cup \vp$, a point from which we start $(i, j)$th  step;
	\im $V_{i, j+1} := (Q\sm S_{i, j}) \cap B_{r_0}(v_{i, j+1})$, a set  used for defining $\iota_{i,j}$;
	\im ${Z_{i,j+1}} =S_{i,j+1}\sm S_{i,j} \su V_{i,j+1}$, which is used to define the stopping set $S_{i,j+1} $. 
	\enen
	Here and in the following, we use the convention $(i,N_i+1)=(i+1,1)$ and $(i,0)=(i-1,N_i)$.

	We initialize the algorithm by $S_{1,1}=\es$. For step $(i,j)$, suppose $S_{i,j}$ has been defined.
	\been
	\im We first choose $v_{i,j+1}$ either as a point in $\vp\cap S_{i,j}$ if $B_{r_0}(\vp\cap S_{i,j})\cap (Q\sm S_{i,j})\ne \es$ or, if  $B_{r_0}(\vp\cap S_{i,j}) \su S_{i,j}$, as a point in $D$. In the latter case, we let $j=N_i$ and hence $(i,j+1)$ equals $(i+1,1)$. The precise algorithm for choosing $v_{i,j+1}$ is quite involved, and will be given later. 
	
	\im Having chosen $v_{i,j+1}$, we set 
	$V_{i, j+1} := B_{r_0}(v_{i, j+1}) \cap (Q \sm S_{i, j}).$
	Define the ordering 
	$$\iota_{i,j} = \tfrac1{2}\iota \one_{V_{i,j+1}} + (\tfrac1{2}\iota +\tfrac1{2})\one_{Q\sm (S_{i,j}\cup V_{i,j+1})}.$$
	Note that the set $V_{i,j+1}$, and hence also the map $\iota_{i,j}$, is completely determined by $\vp \cap S_{i,j}$.
	\im Now,  $Z_{i, j+1}$ is defined as $Z_{i, j+1} := V_{i, j+1}$ if $\vp \cap V_{i, j+1} = \es$, and otherwise
	$$Z_{i,j+1}:=\big\{x \in V_{i,j+1}:\,\iota_{i,j}(x)\le \inf\iota_{i,j}(\vp \cap V_{i, j+1})\big\}.$$
	Finally, set $S_{i,j+1} := S_{i,j}\cup Z_{i,j+1}$. Then, $\iota_{i,j}$ satisfies \eqref{eq:iota_n}.
	\enen
	The algorithm terminates as soon as $S_{i,j}=Q$ for some $i \ge 1$, $j \ge 1$. We will see below, that this happens after finitely many steps.
	
	The regions $Z_{i,j}$ form a partition of $Q$ and satisfy $Z_{i,j} \su V_{i,j} \su V_i := \bigcup_{j=1}^{N_i} V_{i,j}$.  Moreover, the construction is such that for each $i$, $V_i \cap \vp$ is either empty or a cluster of $B_{r_0/2}(\vp ) $. Thus, the algorithm is exploring one cluster at a time just like the one in Example \ref{ex:dc1}. Moreover, each $Z_{i,j}$ is constructed to contain at most one point from $\vp$, which we call 
	$$x_{i,j} := \vp \cap Z_{i, j}.$$
	Note that $x_{i,j} $ may be the empty set.  
	Thus, each cluster is explored (at most) one point at a time.
	Figure \ref{fig:2} illustrates how the algorithm proceeds.
	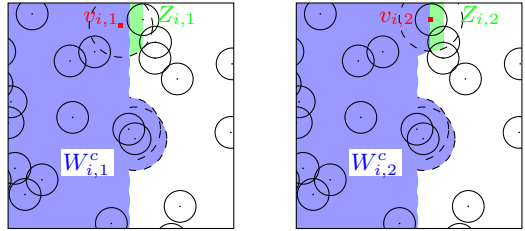
\begin{figure}[!h]
		\begin{center}
			\begin{tikzpicture}
\begin{scope};
\clip(1,2) rectangle (4,5);
\fill[green!40!white] (2.5, 4.7) circle (12.0pt) ;
\fill[white] (2.8, 0) rectangle (5,5);
\fill[blue!40!white] (0, 0) rectangle (2.3,5);
\fill[blue!40!white] (2.6041712068152414, 3.3210249874842432) circle (12.0pt) ;
\fill[blue!40!white] (2.6873888546952003, 3.193821843973229) circle (12.0pt) ;
\draw[dashed] (2.6041712068152414, 3.3210249874842432) circle (12.0pt) ;
\draw[dashed] (2.6873888546952003, 3.193821843973229) circle (12.0pt) ;
\fill[blue!40!white] (2.2, 0.0) circle (12.0pt) ;
\fill[blue!40!white] (2.2, 0.2631578947368421) circle (12.0pt) ;
\fill[blue!40!white] (2.2, 0.5263157894736842) circle (12.0pt) ;
\fill[blue!40!white] (2.2, 0.7894736842105263) circle (12.0pt) ;
\fill[blue!40!white] (2.2, 1.0526315789473684) circle (12.0pt) ;
\fill[blue!40!white] (2.2, 1.3157894736842104) circle (12.0pt) ;
\fill[blue!40!white] (2.2, 1.5789473684210527) circle (12.0pt) ;
\fill[blue!40!white] (2.2, 1.8421052631578947) circle (12.0pt) ;
\fill[blue!40!white] (2.2, 2.1052631578947367) circle (12.0pt) ;
\fill[blue!40!white] (2.2, 2.3684210526315788) circle (12.0pt) ;
\fill[blue!40!white] (2.2, 2.631578947368421) circle (12.0pt) ;
\fill[blue!40!white] (2.2, 2.894736842105263) circle (12.0pt) ;
\fill[blue!40!white] (2.2, 3.1578947368421053) circle (12.0pt) ;
\fill[blue!40!white] (2.2, 3.4210526315789473) circle (12.0pt) ;
\fill[blue!40!white] (2.2, 3.6842105263157894) circle (12.0pt) ;
\fill[blue!40!white] (2.2, 3.9473684210526314) circle (12.0pt) ;
\fill[blue!40!white] (2.2, 4.2105263157894735) circle (12.0pt) ;
\fill[blue!40!white] (2.2, 4.473684210526316) circle (12.0pt) ;
\fill[blue!40!white] (2.2, 4.7368421052631575) circle (12.0pt) ;
\fill[blue!40!white] (2.2, 5.0) circle (12.0pt) ;
\fill[red] (2.47,4.67) rectangle (2.53, 4.73);
\draw[dashed] (2.5, 4.7) circle (12.0pt) ;
\draw (2.6042, 3.3210) circle (6.0pt);
\fill (2.6042, 3.3210) circle (0.3pt);
\draw (1.2253, 2.4525) circle (6.0pt);
\fill (1.2253, 2.4525) circle (0.3pt);
\draw (2.6770, 0.2706) circle (6.0pt);
\fill (2.6770, 0.2706) circle (0.3pt);
\draw (2.8037, 1.0638) circle (6.0pt);
\fill (2.8037, 1.0638) circle (0.3pt);
\draw (4.8615, 3.9673) circle (6.0pt);
\fill (4.8615, 3.9673) circle (0.3pt);
\draw (0.8890, 4.3527) circle (6.0pt);
\fill (0.8890, 4.3527) circle (0.3pt);
\draw (1.1218, 0.7770) circle (6.0pt);
\fill (1.1218, 0.7770) circle (0.3pt);
\draw (0.7808, 2.5580) circle (6.0pt);
\fill (0.7808, 2.5580) circle (0.3pt);
\draw (3.3028, 1.5346) circle (6.0pt);
\fill (3.3028, 1.5346) circle (0.3pt);
\draw (4.0574, 0.3443) circle (6.0pt);
\fill (4.0574, 0.3443) circle (0.3pt);
\draw (2.9599, 4.4791) circle (6.0pt);
\fill (2.9599, 4.4791) circle (0.3pt);
\draw (0.4810, 2.1533) circle (6.0pt);
\fill (0.4810, 2.1533) circle (0.3pt);
\draw (2.1506, 4.3503) circle (6.0pt);
\fill (2.1506, 4.3503) circle (0.3pt);
\draw (0.4715, 2.3271) circle (6.0pt);
\fill (0.4715, 2.3271) circle (0.3pt);
\draw (3.1244, 0.0686) circle (6.0pt);
\fill (3.1244, 0.0686) circle (0.3pt);
\draw (3.9625, 3.2784) circle (6.0pt);
\fill (3.9625, 3.2784) circle (0.3pt);
\draw (1.9926, 0.5429) circle (6.0pt);
\fill (1.9926, 0.5429) circle (0.3pt);
\draw (1.8580, 3.4781) circle (6.0pt);
\fill (1.8580, 3.4781) circle (0.3pt);
\draw (1.3197, 4.8543) circle (6.0pt);
\fill (1.3197, 4.8543) circle (0.3pt);
\draw (1.0908, 2.8003) circle (6.0pt);
\fill (1.0908, 2.8003) circle (0.3pt);
\draw (2.9523, 4.2649) circle (6.0pt);
\fill (2.9523, 4.2649) circle (0.3pt);
\draw (3.9781, 4.1903) circle (6.0pt);
\fill (3.9781, 4.1903) circle (0.3pt);
\draw (1.5637, 1.1690) circle (6.0pt);
\fill (1.5637, 1.1690) circle (0.3pt);
\draw (4.9887, 4.9787) circle (6.0pt);
\fill (4.9887, 4.9787) circle (0.3pt);
\draw (0.5111, 1.8644) circle (6.0pt);
\fill (0.5111, 1.8644) circle (0.3pt);
\draw (1.0420, 1.6837) circle (6.0pt);
\fill (1.0420, 1.6837) circle (0.3pt);
\draw (1.8243, 4.2292) circle (6.0pt);
\fill (1.8243, 4.2292) circle (0.3pt);
\draw (3.0690, 0.4020) circle (6.0pt);
\fill (3.0690, 0.4020) circle (0.3pt);
\draw (3.2870, 2.2939) circle (6.0pt);
\fill (3.2870, 2.2939) circle (0.3pt);
\draw (2.3730, 2.0670) circle (6.0pt);
\fill (2.3730, 2.0670) circle (0.3pt);
\draw (0.3721, 2.2662) circle (6.0pt);
\fill (0.3721, 2.2662) circle (0.3pt);
\draw (3.7440, 1.5005) circle (6.0pt);
\fill (3.7440, 1.5005) circle (0.3pt);
\draw (0.3201, 4.9745) circle (6.0pt);
\fill (0.3201, 4.9745) circle (0.3pt);
\draw (4.8990, 2.9013) circle (6.0pt);
\fill (4.8990, 2.9013) circle (0.3pt);
\draw (0.1307, 2.8998) circle (6.0pt);
\fill (0.1307, 2.8998) circle (0.3pt);
\draw (3.1697, 1.8555) circle (6.0pt);
\fill (3.1697, 1.8555) circle (0.3pt);
\draw (1.4510, 2.6482) circle (6.0pt);
\fill (1.4510, 2.6482) circle (0.3pt);
\draw (2.7933, 4.7783) circle (6.0pt);
\fill (2.7933, 4.7783) circle (0.3pt);
\draw (0.6886, 4.9768) circle (6.0pt);
\fill (0.6886, 4.9768) circle (0.3pt);
\draw (0.8517, 3.9251) circle (6.0pt);
\fill (0.8517, 3.9251) circle (0.3pt);
\draw (4.8034, 2.4744) circle (6.0pt);
\fill (4.8034, 2.4744) circle (0.3pt);
\draw (3.4956, 1.7541) circle (6.0pt);
\fill (3.4956, 1.7541) circle (0.3pt);
\draw (2.2296, 0.6791) circle (6.0pt);
\fill (2.2296, 0.6791) circle (0.3pt);
\draw (3.2492, 3.9898) circle (6.0pt);
\fill (3.2492, 3.9898) circle (0.3pt);
\draw (4.3884, 1.2632) circle (6.0pt);
\fill (4.3884, 1.2632) circle (0.3pt);
\draw (1.0393, 3.6856) circle (6.0pt);
\fill (1.0393, 3.6856) circle (0.3pt);
\draw (2.6874, 3.1938) circle (6.0pt);
\fill (2.6874, 3.1938) circle (0.3pt);
\draw (0.9999, 3.4079) circle (6.0pt);
\fill (0.9999, 3.4079) circle (0.3pt);
\draw (0.0100, 1.7079) circle (6.0pt);
\fill (0.0100, 1.7079) circle (0.3pt);
\draw (0.8285, 2.6457) circle (6.0pt);
\fill (0.8285, 2.6457) circle (0.3pt);
\end{scope};
\fill[white] (1.7,2.7) rectangle (2.4,3.05);
\coordinate[label={\small{\textcolor{blue}{$W_{i, 1}^c$}}}] (A) at (2.04, 2.55);
\coordinate[label={\small{\textcolor{green}{$Z_{i, 1}$}}}] (A) at (3.24, 4.55);
\coordinate[label={\small{\textcolor{red}{$v_{i, 1}$}}}] (A) at (2.24, 4.55);
\draw (1,2) rectangle (4,5);
\end{tikzpicture}\qquad
			\begin{tikzpicture}
\begin{scope};
\clip(1,2) rectangle (4,5);
\fill[blue!40!white] (2.5, 4.7) circle (12.0pt) ;
\fill[white] (2.79, 0) rectangle (5,5);
\begin{scope};
\clip(2.79,0) rectangle (2.96,5);
\fill[green!40!white] (2.79, 4.78) circle (12.0pt) ;
\end{scope};
\fill[blue!40!white] (0, 0) rectangle (2.3,5);
\fill[blue!40!white] (2.6041712068152414, 3.3210249874842432) circle (12.0pt) ;
\fill[blue!40!white] (2.6873888546952003, 3.193821843973229) circle (12.0pt) ;
\draw[dashed] (2.6041712068152414, 3.3210249874842432) circle (12.0pt) ;
\draw[dashed] (2.6873888546952003, 3.193821843973229) circle (12.0pt) ;
\fill[blue!40!white] (2.2, 0.0) circle (12.0pt) ;
\fill[blue!40!white] (2.2, 0.2631578947368421) circle (12.0pt) ;
\fill[blue!40!white] (2.2, 0.5263157894736842) circle (12.0pt) ;
\fill[blue!40!white] (2.2, 0.7894736842105263) circle (12.0pt) ;
\fill[blue!40!white] (2.2, 1.0526315789473684) circle (12.0pt) ;
\fill[blue!40!white] (2.2, 1.3157894736842104) circle (12.0pt) ;
\fill[blue!40!white] (2.2, 1.5789473684210527) circle (12.0pt) ;
\fill[blue!40!white] (2.2, 1.8421052631578947) circle (12.0pt) ;
\fill[blue!40!white] (2.2, 2.1052631578947367) circle (12.0pt) ;
\fill[blue!40!white] (2.2, 2.3684210526315788) circle (12.0pt) ;
\fill[blue!40!white] (2.2, 2.631578947368421) circle (12.0pt) ;
\fill[blue!40!white] (2.2, 2.894736842105263) circle (12.0pt) ;
\fill[blue!40!white] (2.2, 3.1578947368421053) circle (12.0pt) ;
\fill[blue!40!white] (2.2, 3.4210526315789473) circle (12.0pt) ;
\fill[blue!40!white] (2.2, 3.6842105263157894) circle (12.0pt) ;
\fill[blue!40!white] (2.2, 3.9473684210526314) circle (12.0pt) ;
\fill[blue!40!white] (2.2, 4.2105263157894735) circle (12.0pt) ;
\fill[blue!40!white] (2.2, 4.473684210526316) circle (12.0pt) ;
\fill[blue!40!white] (2.2, 4.7368421052631575) circle (12.0pt) ;
\fill[blue!40!white] (2.2, 5.0) circle (12.0pt) ;
\fill[red] (2.76, 4.75) rectangle (2.82, 4.81);
\draw[dashed] (2.79, 4.78) circle (12.0pt) ;
\draw (2.6042, 3.3210) circle (6.0pt);
\fill (2.6042, 3.3210) circle (0.3pt);
\draw (1.2253, 2.4525) circle (6.0pt);
\fill (1.2253, 2.4525) circle (0.3pt);
\draw (2.6770, 0.2706) circle (6.0pt);
\fill (2.6770, 0.2706) circle (0.3pt);
\draw (2.8037, 1.0638) circle (6.0pt);
\fill (2.8037, 1.0638) circle (0.3pt);
\draw (4.8615, 3.9673) circle (6.0pt);
\fill (4.8615, 3.9673) circle (0.3pt);
\draw (0.8890, 4.3527) circle (6.0pt);
\fill (0.8890, 4.3527) circle (0.3pt);
\draw (1.1218, 0.7770) circle (6.0pt);
\fill (1.1218, 0.7770) circle (0.3pt);
\draw (0.7808, 2.5580) circle (6.0pt);
\fill (0.7808, 2.5580) circle (0.3pt);
\draw (3.3028, 1.5346) circle (6.0pt);
\fill (3.3028, 1.5346) circle (0.3pt);
\draw (4.0574, 0.3443) circle (6.0pt);
\fill (4.0574, 0.3443) circle (0.3pt);
\draw (2.9599, 4.4791) circle (6.0pt);
\fill (2.9599, 4.4791) circle (0.3pt);
\draw (0.4810, 2.1533) circle (6.0pt);
\fill (0.4810, 2.1533) circle (0.3pt);
\draw (2.1506, 4.3503) circle (6.0pt);
\fill (2.1506, 4.3503) circle (0.3pt);
\draw (0.4715, 2.3271) circle (6.0pt);
\fill (0.4715, 2.3271) circle (0.3pt);
\draw (3.1244, 0.0686) circle (6.0pt);
\fill (3.1244, 0.0686) circle (0.3pt);
\draw (3.9625, 3.2784) circle (6.0pt);
\fill (3.9625, 3.2784) circle (0.3pt);
\draw (1.9926, 0.5429) circle (6.0pt);
\fill (1.9926, 0.5429) circle (0.3pt);
\draw (1.8580, 3.4781) circle (6.0pt);
\fill (1.8580, 3.4781) circle (0.3pt);
\draw (1.3197, 4.8543) circle (6.0pt);
\fill (1.3197, 4.8543) circle (0.3pt);
\draw (1.0908, 2.8003) circle (6.0pt);
\fill (1.0908, 2.8003) circle (0.3pt);
\draw (2.9523, 4.2649) circle (6.0pt);
\fill (2.9523, 4.2649) circle (0.3pt);
\draw (3.9781, 4.1903) circle (6.0pt);
\fill (3.9781, 4.1903) circle (0.3pt);
\draw (1.5637, 1.1690) circle (6.0pt);
\fill (1.5637, 1.1690) circle (0.3pt);
\draw (4.9887, 4.9787) circle (6.0pt);
\fill (4.9887, 4.9787) circle (0.3pt);
\draw (0.5111, 1.8644) circle (6.0pt);
\fill (0.5111, 1.8644) circle (0.3pt);
\draw (1.0420, 1.6837) circle (6.0pt);
\fill (1.0420, 1.6837) circle (0.3pt);
\draw (1.8243, 4.2292) circle (6.0pt);
\fill (1.8243, 4.2292) circle (0.3pt);
\draw (3.0690, 0.4020) circle (6.0pt);
\fill (3.0690, 0.4020) circle (0.3pt);
\draw (3.2870, 2.2939) circle (6.0pt);
\fill (3.2870, 2.2939) circle (0.3pt);
\draw (2.3730, 2.0670) circle (6.0pt);
\fill (2.3730, 2.0670) circle (0.3pt);
\draw (0.3721, 2.2662) circle (6.0pt);
\fill (0.3721, 2.2662) circle (0.3pt);
\draw (3.7440, 1.5005) circle (6.0pt);
\fill (3.7440, 1.5005) circle (0.3pt);
\draw (0.3201, 4.9745) circle (6.0pt);
\fill (0.3201, 4.9745) circle (0.3pt);
\draw (4.8990, 2.9013) circle (6.0pt);
\fill (4.8990, 2.9013) circle (0.3pt);
\draw (0.1307, 2.8998) circle (6.0pt);
\fill (0.1307, 2.8998) circle (0.3pt);
\draw (3.1697, 1.8555) circle (6.0pt);
\fill (3.1697, 1.8555) circle (0.3pt);
\draw (1.4510, 2.6482) circle (6.0pt);
\fill (1.4510, 2.6482) circle (0.3pt);
\draw (2.7933, 4.7783) circle (6.0pt);
\fill (2.7933, 4.7783) circle (0.3pt);
\draw (0.6886, 4.9768) circle (6.0pt);
\fill (0.6886, 4.9768) circle (0.3pt);
\draw (0.8517, 3.9251) circle (6.0pt);
\fill (0.8517, 3.9251) circle (0.3pt);
\draw (4.8034, 2.4744) circle (6.0pt);
\fill (4.8034, 2.4744) circle (0.3pt);
\draw (3.4956, 1.7541) circle (6.0pt);
\fill (3.4956, 1.7541) circle (0.3pt);
\draw (2.2296, 0.6791) circle (6.0pt);
\fill (2.2296, 0.6791) circle (0.3pt);
\draw (3.2492, 3.9898) circle (6.0pt);
\fill (3.2492, 3.9898) circle (0.3pt);
\draw (4.3884, 1.2632) circle (6.0pt);
\fill (4.3884, 1.2632) circle (0.3pt);
\draw (1.0393, 3.6856) circle (6.0pt);
\fill (1.0393, 3.6856) circle (0.3pt);
\draw (2.6874, 3.1938) circle (6.0pt);
\fill (2.6874, 3.1938) circle (0.3pt);
\draw (0.9999, 3.4079) circle (6.0pt);
\fill (0.9999, 3.4079) circle (0.3pt);
\draw (0.0100, 1.7079) circle (6.0pt);
\fill (0.0100, 1.7079) circle (0.3pt);
\draw (0.8285, 2.6457) circle (6.0pt);
\fill (0.8285, 2.6457) circle (0.3pt);
\end{scope};
\fill[white] (1.7,2.7) rectangle (2.4,3.05);
\coordinate[label={\small{\textcolor{blue}{$W_{i, 2}^c$}}}] (A) at (2.04, 2.55);
\coordinate[label={\small{\textcolor{green}{$Z_{i, 2}$}}}] (A) at (3.44, 4.55);
\coordinate[label={\small{\textcolor{red}{$v_{i, 2}$}}}] (A) at (2.34, 4.55);
\draw (1,2) rectangle (4,5);
\end{tikzpicture}
		\end{center}
		\caption{Illustration of two steps in the exploration algorithm. The solid balls represent $B_{r_0/2}(\vp)$, the blue region is the already explored region $S_{i, j}$, and the green region is the currently explored region $Z_{i, j}$. In the left figure, a new cluster is explored from a point $v_{i, 1}\in D$ (red square). In the right figure, the exploration is continued by using a point from $\vp$ as $v_{i,2}$ (red square).}
		\label{fig:2}
	\end{figure}

	We now give the precise construction of the points $v_{i, j}$.  In the first step, we let $v_{1, 1}$ be the $\preceq$-smallest point in $D$ such that  $ B_{r_0}(v_{1,1}) \cap Q \cap B_{r_0}(B ) \ne \es$  if such a point exists and otherwise just the smallest point in $D$ such that $ B_{r_0}(v_{1,1}) \cap Q \ne \es$.
	
	Next, we explain how to construct $v_{i,j+1}$ from $S_{i,j}$. Let
	$$k_{i,j+1}:=\inf \{k \in \{1,\dots,j\}:  (Q\sm S_{i,j}) \cap B_{r_0} (x_{i,k})\ne \es\}$$
	be the first index $k$ where $B_{r_0} (x_{i,k})$ is not contained in $S_{i,j}$ (with the convention $B_{r_0}(\es)=\es$). If $k_{i,j+1}<\ff$, we define $v_{i, j+1} := x_{i,k_{i,j+1} }$. Otherwise, if  $k_{i, j+1 }= \ff$, we set $N_i:=j$. In this case, $ v_{i,j+1}=v_{i+1,1} $. If $(Q\sm S_{i, j}) \cap B_{r_0}(B) \ne \es$, we define  $ v_{i+1,1} $ as the $\preceq$-smallest element of $D$ such that $B_{r_0}(v_{i+1,1})\cap (Q\sm S_{i,j}) \cap B_{r_0}(B ) \ne \es$. Else, if $B_{r_0} (B) \su S_{i ,j} $, define  $v_{i+1,1} $ as the smallest element of $D$ such that $B_{r_0}(v_{i+1,1})\cap (Q\sm S_{i,j})\ne \es$). 	
	This construction ensures that we move through the components of $B_{r_0/2} (\vp ) $ one at a time, starting with those intersecting $B_{r_0/2}(B) $ so that $\psi\cap B$ affect the thinning only on those components. 
	
	To see that the algorithm terminates after a finite number of steps (by which we mean that $\sum_{i \le I} N_i <\ff$), note that all $v_{i,1}$ are pairwise distinct points in $D$. Since $D$ is finite, $I<\ff$. Moreover, each $N_i$ is finite, since it is bounded by two times the number of points in the corresponding cluster. 
	
	Consider the case when $\iota$ is an arbitrary injective ordering and
	$D :=\de  \Z^d$, where $\de>0$ is sufficiently small so that $Q_{2\de}\su B_{r_0}$ and hence $B_{r_0}(D)  =\R^d$. Let  $\iota_\infty(x) = \sup_{1\le i\le d} |x_i|$ for $x=(x_1,\ldots,x_d)$. We order $x,y\in D$ by declaring $x<y$ if $\iota_\ff(x)<\iota_\ff(y)$ or $\iota_\ff(x)=\iota_\ff(y)$ and $x$ is lexicographically smaller than $y$. 	
	We call the resulting thinning algorithm the \emph{radial coupling} and denote the corresponding thinning operator by $\Tra_{ Q, B, \psi}$. When $B=Q^c$, we may simply write $\Tra_{ Q, \psi}$. The order in which the algorithm visits the clusters of $B_{r_0/2}(\vp )$ is given by the smallest point in $D$ within distance $r_0$ from a point in the cluster. 
\end{example}

\subsection{Homogeneity and decorrelation via disagreement coupling}
\label{ss:hd}

The fact that the radial coupling works on the clusters of the Boolean model $B_{r_0/2}(\Pd)$ allows us to relate the Gibbs process to the percolation properties of the Boolean model. Since we assume $r_0<r_c(\a_0)$, we already know that the connected components of $B_{r_0/2}(\Pd)$ are almost surely all finite. Moreover, we shall heavily rely on the following fundamental result from continuum percolation known as the \emph{sharp phase transition}. For Borel sets $A, B \su \R^d$, we let $\{A \lrsa B\}$ denote the event that there exists a connected component of the Boolean model $B_{r_0/2}(\Pd)$ intersecting both $B_{r_0/2}(A) $ and $B_{r_0/2}(B) $.  Then, for $r < r_c(\a_0)$, the connection probabilities $\P(A \lrsa B)$ decays exponentially fast in the distance between $A$ and $B$.

\bepr[Sharp phase transition, \cite{raoufi_sub}]
\label{pr:spt}
Let $r_0 < r_c(\a_0)$. Then, there are constants $c_{\ms{SPT}, 1},c_2>0$ such that
$$  \P\big(o \lrsa \pa B_s \big) < c_{\ms{SPT}, 1} \exp(-c_{\ms{SPT}, 2} s).$$
\enpr
\bep
See e.g.~\cite[Theorem 1.4]{raoufi_sub}.
\enp

Proposition \ref{pr:rad_prop} immediately has the following corollary, where $\dist(A,B)=\inf\{|x-y|\mid x\in A, y\in B\}$ denotes the distance between two Borel sets $A,B \su \R^d$.

\bec[Disagreement probabilities]
\label{cor:dperc}
Let $A\su Q \su \R^d$ be bounded Borel sets, $B\su Q^c$ and $\psi \in \Nlf_{Q^c}$. Then, 
$$\sup_{\substack{ \psi' \in \Nlf_{Q^c}\\ \psi\cap B^c=\psi'\cap B^c}}\hspace{-0.2cm} \P\big(\Tra_{Q,B, \psi}(\Pds_Q) \cap A \ne \Tra_{Q,B, \psi'}(\Pds_Q) \cap A\big) \le \P(A \lrsa B)\le c_{\ms{SPT}, 1}\big|B_{r_0}(A) \big|\exp(-c_{\ms{SPT}, 2} \dist(A,B)).$$
\enc 

\bep

It follows from Proposition \ref{pr:rad_prop} that $\Tra_{Q,B, \psi}(\Pds_Q)$ and $\Tra_{Q,B, \psi'}(\Pds_Q)$ agree on $A$ when $B_{r_0/2}(B)  $ is not connected to $B_{r_0/2}(A)$. Moreover, 
\begin{equation}\label{eq:conn_prob}
	\P({A\lrsa B})\le \E\Big[ \sum_{x\in \Pd\cap  B_{r_0}(A)} \one_{\{x \lrsa B\}}  \Big] \le \int_{B_{r_0}(A) } \P(x \lrsa B_{\dist(A,B)}(x))\d x.
\end{equation}
\enp

For the rest of this section, we highlight how the relation between Gibbs processes and the associated Boolean model provided by disagreement coupling can be used to establish homogeneity and decorrelation of  Gibbs point processes satisfying \eqref{as:(A)}.
 The homogeneity and decorrelation are captured through the total variation distance and the $\a$-mixing coefficient, respectively. Thus, we define the \emph{total variation distance} between two random variables $X, Y$ with values in a common measurable space $S$ by
\begin{align}
\label{eq:dtv}
	\dtv(X, Y) := \sup_{f\co S \to [0, 1]}\big|\E[f(X)] - \E[f(Y)]\big|,
\end{align}
where the supremum is taken over all measurable functions $f\co S \to [0, 1]$ .
Moreover, we define the \emph{$\a$-mixing coefficient} of two $\s$-algebras $\AA$ and $\BB$ by
\begin{align*}
\a(\AA, \BB) := \sup_{A\in \AA, B\in \BB} \big|\P(A \cap B) - \P(A) \P(B)\big|.
\end{align*}

%
%
\bec[Exponential decay of the total variation distance]
\label{cor:dec} 
Suppose that the PI satisfies condition \eqref{as:(A)}. For any bounded Borel sets $A\su Q \su \R^d$ and  $B\su  Q^c$, it holds that
$$
\sup_{\substack{\psi, \psi' \in \Nlf_{Q^c}\\ \psi\cap B^c=\psi'\cap B^c}} \dtv\big(\XX(Q, \psi)\cap A, \XX(Q,\psi') \cap A\big) \le c_{\ms{SPT}, 1} \big|B_{r_0}(A)\big| e^{-c_{\ms{SPT}, 2} \dist(A,B)}.
$$
\enc
\bep
 We realize $\XX(Q,\psi)$ and $\XX(Q, \psi')$ by disagreement couplings  $\Tra_{Q,B, \psi}(\Pds_Q)$ and $\Tra_{Q,B, \psi'}(\Pds_Q)$, respectively. Since
$$\dtv\big((\XX(Q, \psi) )\cap A, \XX(Q, \psi') \cap A\big) \le  \P\big(\Tra_{Q,B, \psi}(\Pds_Q) \cap A \ne \Tra_{Q,B, \psi'}(\Pds_Q) \cap A\big), $$
the results follow from Corollary \ref{cor:dperc}.
\enp

	From Corollary \ref{cor:dec}, we can derive useful homogeneity and decorrelation bounds. We note that \cite[Theorem 2]{benes} also establishes a decorrelation property via disagreement coupling, where in their case the spatial dependence is measured via the $k$-point correlation functions. However, we need $\a$-mixing for our further arguments, which is why we include the short proof of Proposition \ref{pr:lhom} below.  We also establish a form of local homogeneity in the sense that for large windows, the local distribution in a neighborhood of any space point does not depend too much on the position of the point in the window.

	\bepr[Spatial decorrelation and local homogeneity]\label{pr:lhom}
	\begin{itemize}
	\item[(i)] Let $s > 0$ and $A\su Q_n$ be  Borel. Then
	$$
			\a\big(\s\big(\XX(Q_n,\es)\cap A\big),\s\big(\XX(Q_n,\es)\sm B_s (A) \big)\big)\le c_{\ms{SPT}, 1}  \big|B_{r_0}(A)\big| e^{-c_{\ms{SPT}, 2}s},	
	$$

	\item[(ii)] Let the PI be translation-invariant. Then
	$$
\sup_{x_0 \in Q_{n - 2\sqrt n}}\dtv\big((\XX(Q_n,\es) - x_0)\cap Q_{\sqrt n}, \XX(Q_n,\es) \cap Q_{\sqrt n}\big) \le 2c_{\ms{SPT}, 1}\exp(-2c_{\ms{SPT}, 2}\sqrt n).
$$
\end{itemize}
	\enpr
	\bep
	We start with the proof of (i).
	 Let $\XX' := \XX(Q_n, \es) \cap A$, $\XX'' := \XX(Q_n, \es) \sm B_s(A) $ and $\YY := \XX(B_s(A) , \XX'') \cap A$. We note that for any measurable sets of configurations $B_1,B_2$, we have 
\begin{align*}
	\P(\XX' \in B_1, \XX'' \in B_2) &= \E\big[\P(\XX' \in B_1\ba \XX'') \one\{\XX'' \in B_2\}\big]\\
	&= \P(\YY \in B_1) \P(\XX''\in B_2) + \E\Big[\big(\P(\XX' \in B_1\ba \XX'') - \P(\YY \in B_1)\big) \one \{\XX'' \in B_2\}\Big]\\
	&\le \P(\YY \in B_1) \P(\XX''\in B_2) + \sup_{\psi \in \mathbf N_{Q_n\sm B_s (A) }} \dtv(\XX(B_s(A) ,\XX'')\cap A,\XX(B_s(A), \psi)\cap A).
\end{align*}
Noting that $\YY$ and $\XX'$ have the same distribution, we conclude the proof of part (i) by applying Corollary \ref{cor:dec} with $Q= B_s(A) $ and $B= Q_n \sm B_s(A)$. 

For the proof of (ii), note that $Q_{2\sqrt n} \su Q_{n }-x_0$ for every $x_0 \in Q_{n - 2\sqrt n}$. Therefore,
\begin{align*}
	&\dtv\big((\XX(Q_n,\es) - x_0)\cap Q_{\sqrt n}, \XX(Q_n,\es) \cap Q_{\sqrt n}\big) \\
	\quad&\le 2\sup_{\psi \in \Nlf_{Q_{2\sqrt n}^c}} \dtv\big(\XX(Q_{2\sqrt n}, \psi)\cap Q_{\sqrt n}, \XX(Q_{2\sqrt n},\es) \cap Q_{\sqrt n}\big),
\end{align*}
where we used that by  stationarity, the distribution of $(\XX(Q_n,\es)-x_0)\cap Q_{2\sqrt n}$ is $\XX(Q_{2\sqrt n },\psi)$, when conditioned on $(\XX(Q_n,\es)-x_0)\cap Q_{2\sqrt n}^c=\psi$. Hence, the claim follows from Corollary \ref{cor:dec}. 
	\enp

%
%
\subsection{Radial coupling for infinite-volume Gibbs processes}
\label{ss:emb}

So far, we have only considered disagreement coupling for Gibbs processes in bounded domains. In this section, we extend disagreement coupling to the infinite-volume Gibbs process. That is, we show that the radial coupling extends to unbounded domains, constructing the infinite-volume Gibbs process as a thinning of $\Pds$. This construction yields explicit couplings between Gibbs processes in bounded and unbounded domains via thinning of a common Poisson process allowing us to quantify the disagreement between these processes. This improves on \cite{dp}, where such couplings were only obtained qualitatively via weak convergence.  The results allow us to show mixing, and hence ergodicity, of the infinite-volume Gibbs process, which is crucial for the proof of Theorem \ref{thm:1}. For completeness, we also include  proof that the infinite-volume Gibbs process is unique, although we remark that this result is not new, having appeared in e.g.\ \cite{dereudre}, see also \cite{benes} for a proof using disagreement coupling.  
 
The main result of this section, provided by Propositions \ref{pr:xxff} and \ref{pr:uniq} below, establishes that the radial thinning on $Q_n$ converges almost surely to the infinite-volume Gibbs process when $n\tff$. 

%
%
\bepr[Convergence of radial couplings]
\label{pr:xxff}
Assume that the PI satisfies \eqref{as:(A)}. Let $U \su \R^d$  and $\psi \in \Nlf_{U^c}$ be locally finite. Then, as $n\tff$, the embeddings $\Tra_{Q_n \cap U, \psi}(\Pds_{Q_n})$ converge almost surely to  a limiting process on $U$
\begin{align}
	\label{eq:xxff}
	 T^{\ff}_{U,\psi}(\Pds_{U}):= \bigcup_{n_0 \ge 1} \bigcap_{n \ge n_0} \Tra_{Q_n \cap U, \psi}(\Pds_{Q_n\cap U})
\end{align}
such that for all bounded Borel sets $A\su \R^d$ and $n\ge 1$,
$$\P\big( T^{\ff}_{U,\psi}(\Pds_{U})  \cap A= \Tra_{Q_{m} \cap U, \psi}(\Pds_{Q_{m}\cap U})\cap A \text{ for all $m \ge n$}\big) \le \P(A \not \lrsa Q_{n-4r_0}^c ).$$
If $\k$ is translation-invariant, the process $T^{\ff}(\Pds) := T^{\ff}_{\R^d,\es}(\Pds)$ is a stationary and mixing limiting process. {Moreover, $T^{\ff}_{U,\psi}(\Pds_U) +x $ has the same distribution as $ T^{\ff}_{U+x,\psi+x}(\Pds_{U+x}) $ for all $x\in \R^d$.}
\enpr

 We stress that, although $T^{\ff}(\Pds)$ is stationary when $\k$ is translation-invariant, the radial thinning is not translation-covariant in the sense that  $T^{\ff}(\Pds+(x,0)) $ is generally not the same as  $T^{\ff}(\Pds) +x$.

%
%
\bepr[Uniqueness of infinite-volume Gibbs process]
\label{pr:uniq}
Assume that the PI satisfies \eqref{as:(A)}. The process $T^{\ff}(\Pds)$ is the distributionally unique infinite-volume Gibbs process with PI $\k$.
\enpr

The rest of this section is devoted to proving Propositions \ref{pr:xxff} and \ref{pr:uniq}. The key ingredient is the following almost-sure consistency property of the disagreement coupling in increasing domains.

%
%
\bel[Consistency of the radial thinning]
\label{lem:rad}
Let $A \su \R^d$ be bounded, $U \su \R^d$ and $\psi \in \Nlf_{U^c}$ be locally finite. Then, almost surely on the event $\{A \not \lrsa Q_{n-4r_0}^c\}$, we have for all $m \ge n$
\begin{equation}\label{e:eqonA}
	\Tra_{Q_n \cap U, \psi}(\Pds_{Q_n}) \cap A = \Tra_{Q_{m} \cap U, \psi}(\Pds_{Q_m})\cap A.
\end{equation}
\enl

%
%

\bep
Let $V_1,  \dots, V_K$ be the regions considered in the radial coupling algorithm on $Q_{m}$ and assume that $V_1,\dots, V_k$ are the ones that start from a point $v_{i,1} \notin Q_{n-2r_0}$. Then, $S_{k,N_k} = {V}_1\cup \dots \cup {V}_{{k}} $ is obtained as the union of $B_{r_0}(D\cap Q_{n-2r_0}^c)\cap (Q_m\cap U) $  and all  $B_{r_0/2}( \mc C_l)\cap (Q_m\cap U)$, $l=1,\ldots,L$, where $\mc C_l$ is a component of $B_{r_0/2}(\Pd_{ Q_{m}\cap U})$ such that $B_{r_0/2}(\mc C_l)$  contains a point in $D\cap Q_{n-2r_0}^c$.
Similarly, the disagreement coupling algorithm on $Q_n$ first considers the $\tilde{V}_1,\dots, \tilde{V}_{\tilde{k}}$ that start from a point $\tilde{v}_{i,1} \notin Q_{n-2r_0}$. Then, $\tilde{S}_{\tilde{k},N_{\tilde{k}}} = \tilde{V}_1\cup \dots \cup \tilde{V}_{\tilde{k}} $ is obtained as the union of $B_{r_0}(D\cap Q_{n-2r_0}^c) \cap (Q_n\cap U) $ and all $B_{r_0/2}(\tilde{\mc C}_{\tilde{l}})\cap (Q_n\cap U)$ where $\tilde{\mc C}_{\tilde{l}}$ is a component of $B_{r_0/2}(\Pd_{Q_n\cap U})$ such that $B_{r_0/2}(\tilde{\mc C}_{\tilde{l}})$, $\tilde{l}=1,\ldots,\tilde{L}$, contains a point in $D\cap Q_{n-2r_0}^c$.

We claim that $(Q_m\cap U) \sm S_{k,N_{k}}=(Q_n\cap U)\sm \tilde{S}_{\tilde{k},N_{\tilde{k}}}$, and hence the disagreement coupling algorithm proceeds the same from steps $k$ and $\tilde{k}$, respectively. On the event $\{A \not \lrsa Q_{n-4r_0}^c\}$, $A\su Q_m\sm S_{k,N_k}=Q_n\sm \tilde{S}_{\tilde{k},N_{\tilde{k}}}$ and hence \eqref{e:eqonA} follows. 

To see that $(Q_m\cap U) \sm S_{k,N_{k}}=(Q_n\cap U)\sm \tilde{S}_{\tilde{k},N_{\tilde{k}}}$, it is enough to show
\begin{equation} \label{eq:union}
	B_{r_0}(D\cap Q_{n-2r_0}^c)\cup \bigcup_l B_{r_0/2}(\mc C_l )  = B_{r_0}(D\cap Q_{n-2r_0}^c)\cup \bigcup_{\tilde{l}} B_{r_0/2}( \tilde{\mc C}_{\tilde{l}}),
\end{equation}
since $Q_m^c \su Q_n^c \su B_{r_0}(D\cap Q_{n-2r_0}^c)$ by the choice of $D$.
Any component $\tilde{\mc C}_{\tilde{l}}$ 
is contained in one of the components $\mc C_l$, which shows one inclusion.
Next, consider any component $\mc C_l$ containing a point in $D\cap Q_{n-2r_0}^c$ within distance $r_0/2$. Then, $\mc C_l\cap Q_n$ must be contained in a union of some of the components in $B_{r_0/2}(\Pd_{Q_n\cap U})$  that either contain a point in $D\cap Q_{n-2r_0}^c$ within distance $r_0/2$ themselves or has distance at most $r_0/2$ to $\partial Q_n$. In the latter case, it also contain a point in $D\cap Q_{n-2r_0}^c$ within distance $r_0/2$, as well as some sets of the form $B_{r_0/2}(x)$ with $x\in \Pd_{U\cap Q_m \sm Q_n}$. Sets of the latter form must satisfy $B_{r_0}(x) \su Q_{n-2r_0}^c\su B_{r_0}(D\cap Q_{n-2r_0}^c)$ by the choice of $\de$ in the definition of $D$. 
In total, this shows that $B_{r_0/2}(\mc C_l )$ is contained in the right hand side of \eqref{eq:union}.
\enp

For the rest of this section, we say that a measurable function $f \co \Nlf \to [0, 1]$ is \emph{local} if there exists $s > 0$ such that $f(\vp) = f(\vp \cap Q_s)$ holds for all $\vp \in \Nlf$. We call  $Q_s$ a locality region of $f$.

\bep[Proof of Proposition \ref{pr:xxff}]
Since $r_0 < r_c$ is sub-critical, the almost-sure convergence in \eqref{eq:xxff} follows from  Lemma \ref{lem:rad}.

For the rest of this proof, we write $\XX:=T^{\ff}(\Pds)$. To prove that $\XX$ is stationary, we must show 
$\E[f(\XX+ x)] = \E[f(\XX)]$
for every $x \in \R^d$ and local $f\co \Nlf \to [0, 1]$. Since $f$ is local and \eqref{eq:xxff} holds, this is equivalent to $\lim_{n \tff} \big(\E[f(\XX(Q_n,\es)+x )] - \E[f(\XX(Q_n,\es))]\big) = 0$. Now, let $x\in Q_{n-2\sqrt n}$, where  $n$ is so large such that $\sqrt n>s$, where $Q_s$ is the locality region of $f$. Then,
$$
|\E[f(\XX(Q_n,\es)+x )] - \E[f(\XX(Q_n,\es))]|\le \dtv((\XX(Q_n,\es)+x) \cap Q_{\sqrt n},\XX(Q_n,\es) \cap Q_{\sqrt n}),
$$
which goes to 0 by part (ii) of Proposition \ref{pr:lhom}, i.e., the spatial homogeneity.

The claim $ T^{\ff}_{U,\psi}(\Pds_{U})  +x \sim  T^{\ff}_{U+x,\psi+x}(\Pds_{U+x})  $ is shown similarly using the straightforward generalization of part (ii) of Proposition \ref{pr:lhom} to Gibbs processes with boundary conditions.

It remains to show that $\XX(\ff)$ is mixing, i.e., 
$$
\lim_{|x|\tff}\E[f(\XX + x)g(\XX)] = \E[f(\XX)]\E[g(\XX)]
$$
for every $x \in \R^d$ and local $f,g\co \Nlf \to [0, 1]$. Let $A \su \R^d$ be a fixed bounded Borel set containing the locality regions of $f$ and $g$. Then, setting $n(x) = |x|^2$, Lemma \ref{lem:rad} gives that 
$$\E\big|[f(\XX + x)g(\XX) - f(\XX(Q_{n(x)},\es) + x)g(\XX(Q_{n(x)},\es))]\big| \le 2\P(A \cup (A + x)\lrsa \pa Q_{n(x)-4r_0}).$$
for all sufficiently large $|x|$. As before, since $\Pd$ is sub-critical, the right-hand side tends to 0 as $|x| \tff$. Now, we need to show that 
$$\lim_{|x|\tff}\Big(\E\big[f(\XX(Q_{n(x)},\es) + x)g(\XX(Q_{n(x)},\es))\big] - \E\big[f(\XX(Q_{n(x)},\es))\big]\E\big[g(\XX(Q_{n(x)},\es))\big]\Big) = 0.$$
Invoking Proposition \ref{pr:lhom} (i), i.e., the spatial decorrelation, concludes the proof of the mixing property.
\enp

%
%
We now show the uniqueness result asserted in Proposition \ref{pr:uniq}.
\bep[Proof of Proposition \ref{pr:uniq}]
We must show that $\XX:=T^{\ff}(\Pds)$ satisfies the GNZ equations \eqref{eGNZ} for any $f:\R^d\times \Nlf \to [0,1]$ with support in $Q_s\times \Nlf$ and such that $f(x,\vp) =f(x,\vp\cap Q_s)$ for some $s>0$. 
The point process $\Tra_{Q_n, \es}(\Pds_{Q_n})$ is a Gibbs point process $\XX(Q_n,\es)$ on $Q_n$ with PI $\k$ and hence it satisfies the GNZ equations \eqref{eGNZ_finite} for $f$. Taking limits, we see that also $\XX(\ff)$ satisfies the GNZ equations \eqref{eGNZ} for $f$ and hence is an infinite-volume Gibbs point process with PI $\k$.  

Now, let $A\su \R^d$ be a bounded Borel set and let $\XX'$  be any infinite-volume Gibbs point processes with PI $\k$.  Then, to sample $\XX'$ in $A$, one may first generate a sample $\psi$ from $\XX'$ in $Q_n^c$ with $Q_n \supseteq A$ and then draw a sample from $\XX' \cap Q_n$ under the boundary condition $\psi$. By the DLR equations \eqref{eq:DLR2}, the conditional distribution of $\XX' \cap Q_n$ given $\psi$ can be represented as $\Tra_{Q_n, \psi}(\Pds_{Q_n})$. Hence, it suffices to show that $\lim_{n \tff}\sup_{\psi, \psi' \in \Nlf_{Q_n^c} } \P(E_n(\psi,\psi'))=0$, where 
$$E_n(\psi, \psi') := \big\{\Tra_{Q_n, \psi}(\Pds_{Q_n}) \cap A \ne \Tra_{Q_n, \psi'}(\Pds_{Q_n})\cap A\big\}.$$
Hence, we conclude by using Corollary \ref{cor:dperc}.
\enp

\begin{remark}
	While the order on $D$, and hence the order in which we thin the components of $B_{r_0/2}(\Pd)$, goes from the outside towards the origin, the proof of Lemma \ref{lem:rad} shows that the radial thinning construction of  $T^{\ff}(\Pds)$ can actually be performed by moving from the origin and outwards by considering a growing sequence of windows $Q_1,Q_2,\dots$ After the thinning of $Q_n$, we keep the thinning we had on components of $B_{r_0/2}(\Pd)$ that do not contain a point in $D \cap Q_{n-2r_0}^c$ within distance $r_0/2$. On $Q_{n+1}$, we only need to run the radial thinning until all components containing a point in $ D\cap Q_{n-2r_0}^c$ within distance $r_0/2$ have been considered, since the remaining components will not change. We continue this process for increasing $n$.
\end{remark}

{
We mention the following corollary, which is immediate from the proof of Lemma \ref{lem:rad}.
\bec[Boundary conditions for radial thinning]\label{cor:rad_bound}
Let  $U \su \R^d$ be Borel  and $\psi,\psi' \su {U^c}$ be locally finite with $Q_n\cap \psi = Q_n \cap \psi'$. Then, almost surely on the event $\{A \not \lrsa Q_{n-4r_0}^c\}$, 
\begin{equation*}
	T^{\ff}_{ U, \psi}(\Pds) \cap A = T^{\ff}_{ U, \psi'}(\Pds)\cap A. 
\end{equation*}
\enc

The following proposition, a version of the DLR-equations in unbounded domains, is used in the proof of Theorem \ref{thm:1}.

%
%
\bepr[DLR-equations in unbounded domains]
\label{pr:chimera}
Let $U \su \R^d$ be Borel. Let $\X=T^{\ff}(\Pds)$ be an infinite-volume Gibbs process on $\R^d$ and let $\XX'=(\X\cap U^c )\cup T^{\ff}_{U,\X\cap U^c}(\tilde{\PP}^{*}_{U})$, where $\tilde{\PP}^{*}$ is an independent copy of $\Pds$. Then, $\XX'$ is again an infinite-volume Gibbs process on $\R^d$.
\enpr

\bep
Let $\X_n=\Tra_{Q_{4n},\es}(\Pds_{Q_{4n}})$. 
Then, by Proposition \ref{pr:xxff},
$$\P(\XX \cap Q_{2n} \ne  \XX_n\cap Q_{2n}) \le \P(Q_{2n} \lrsa Q_{4n}).$$
Moreover, by Corollary \ref{cor:rad_bound}, 
$$\P(T^{\ff}_{U,\X_n \cap U^c}(\tilde{\PP}^{*}_{U}) \cap Q_n \ne T^{\ff}_{U,\X \cap U^c}(\tilde{\PP}^{*}_{U}) \cap Q_n)\le \P(Q_n \lrsa Q_{2n-2r_0})+\P(Q_{2n} \lrsa Q_{4n}).$$
Finally, 
$$\P\Big(T^{\ff}_{U,\X_n \cap U^c}(\tilde{\PP}^{*}_{U}) \cap Q_n \ne \Tra_{U\cap Q_{4n},\X_n \cap U^c}(\tilde{\PP}^*_{U\cap Q_{4n}})\cap Q_n \Big)\le \P(Q_n \lrsa Q_{4n-2r_0}).$$
By the DLR-equations \eqref{eq:DLR1},  $\XX'_n=(\X_n\cap U^c )\cup T^{\ms{rad}}_{U\cap Q_{4n},\X_n \cap U^c}(\tilde{\PP}^{*}_{U\cap Q_{4n}})$ has the distribution of the Gibbs process $\XX(Q_{4n},\es)$. If $f:\Nlf \to [0,1]$ is a bounded local function with locality region $A$, the above shows that $\lim_{n\to \infty} f(\XX'_n) = f(\XX')$. Moreover, Proposition \ref{pr:xxff} implies that $\XX'_n\cap A$ converges to $\XX\cap A$  in the total variation distance. Hence,  dominated convergence yields
$$\E[f(\XX')] = \lim_{n\to \infty} \E[f(\XX'_n)] =\E[f(\XX)],$$
which finishes the proof.
\enp

}

%
%
\subsection{Disagreement coupling for perturbed PIs}
\label{ss:rfin}

The goal of this section is to control the disagreement probabilities for Gibbs processes with PIs that differ locally. More precisely, let $Q \su \R^d$ be Borel, $\psi, \psi'\su Q^c$ be locally finite, and $\k, \k'$ be PIs satisfying the assumption \eqref{as:(A)}. Assume that $\k$ and $k'$ are small perturbations of each other, i.e., there is a (usually large) set $P\su Q$ such that for all $x \in P$ and $\mu \in \mathbf N$,	
\begin{equation}\label{eq:k,k'}
\k(x,\mu)=\k'(x,\mu). 	
\end{equation}
In the following, let $\XX(Q,\psi)$ and $\XX'(Q,\psi')$ be Gibbs processes with PI $\k$ and $\k'$ and boundary conditions $\psi,\psi'\in \Nlf_{Q^c}$, respectively. Throughout this section, we use a "prime" to denote quantities associated with $\kappa'$. We know already from Section \ref{ss:hd} that the total variation distance between these processes on a set $A\su P$ is small.

\bec[Total variation distance for different PIs]\label{cor:kappa'}
Suppose that $\k, \k'$ are PIs satisfying condition \eqref{as:(A)}.  Assume that $A\su P\su Q$ and $\k,\k'$ satisfy \eqref{eq:k,k'}.  Let $\psi,\psi'\in \Nlf_{Q^c}$ and $B\su Q^c$ such that $\psi\cap B^c = \psi'\cap B^c$. Then,
\begin{equation}\label{eq:dtv_simple}
\dtv(\XX(Q,\psi)\cap A\ne \XX'(Q,\psi')\cap A) \le \P(A\lrsa  B\cup (Q\sm P)).	
\end{equation}
\enc

\bep
First, consider any realizations $\XX(Q,\psi)\sm P =\xi$ and $\XX'(Q,\psi')\sm P=\xi'$. Extend this over $Q'$ by using disagreement couplings $\Tra_{P, B\cup (Q\sm P) ,\xi \cup \psi}(\Pds_Q)$ and $\Tra_{P, B\cup (Q\sm P) ,\xi' \cup \psi'}(\Pds_Q)$, which have the distributions of $\XX(Q,\psi)\cap P$ and $\XX'(Q,\psi')\cap P$ by the DLR-equation \eqref{eq:DLR1}. Then,
\begin{equation}\label{eq:disprob}
	\P(\Tra_{P,B\cup (Q\sm P) , \xi \cup \psi}(\Pds_Q)\cap A\ne \Tra_{P,B\cup (Q\sm P) , \xi' \cup \psi'}(\Pds_Q) \cap A) \le \P\big(A\lrsa B\cup (Q\sm P)\big).	
\end{equation}
Hence, \eqref{eq:dtv_simple} follows from Corollary \ref{cor:dperc} since  
\enp

However, in the proof of Theorem \ref{thm:2}, we will have a whole family of PIs and we would like to be able to simultaneously construct the corresponding Gibbs  processes and still have control over all the pairwise disagreement probabilities rather than just the total variation distance. In particular, the set $P$ may depend on which pair we compare, so the strategy for constructing the Gibbs processes in the proof of Corollary \ref{cor:kappa'} does not work. {Moreover, in proof of Theorem \ref{thm:1}, we need to control how local changes in boundary conditions affect the radial thinning in unbounded domains.} The theorem below, which is the main result of this section, allows us to control the disagreement probabilities in these settings. We remark that neither the theorem nor its proof provides precise information about how the difference between the Gibbs processes is related to the components of $B_{r_0/2}$.

%
%
\bet[Disagreement probabilities for differing PIs] \label{th:disdis}
There is a constant $c_{\ms{DP}}>0$ with the following property.
Let $Q \su \R^d$ be bounded Borel and $\iota: Q\to[0,1)$ be injective.
Suppose that $\k, \k'$ are PIs satisfying condition \eqref{as:(A)} and $\psi,\psi'\in \Nlf_{Q^c}$ \cb{ and $\psi\cap B^c = \psi'\cap B^c$ for some $B\su Q^c$}.  Assume that $s > 0$ and  $A\su P\su Q$ are such that \cb{$\dist(A,(Q_m\sm P) \cup B) \ge 3s$} and $\k, \k'$ satisfy \eqref{eq:k,k'}. Then,
\begin{align}
	\P\big(T^{{\sf{rad}}}_{Q,\psi}(\Pds_Q)\cap A\ne T^{'\sf{rad}}_{Q,\psi'}(\Pds_Q) \cap A\big) \le c_{\ms{DP}} \big|B_{2s}(A)\big| e^{-c_{\ms{SPT}, 2}s}.\label{eq:XXbou}
\end{align}
\ent

\cb{
As a corollary, we obtain the analogous statement for Gibbs processes in unbounded domains.

\bec\label{cor:ff_disag}
There are  $c_{\ms{DPP}},c_{\ms{DPP}}'>0$ such that the following holds.
Let $U\su \R^d$, $A\su P\su U$ and $B\su U^c$  such that $A$,$B$, and $U\sm P$ are bounded. Let $\psi,\psi'\in \Nlf_{U^c}$ be such that $\psi \cap B^c = \psi' \cap B^c $. Suppose that the PIs $\kappa,\kappa'$ satisfy \eqref{as:(A)} and \eqref{eq:k,k'}. Let $s=\dist(A,(U\sm P) \cup B)$. Then,
\begin{equation}\label{eq:ff-disdis}
	\P( T^{\ff}_{U,\psi} (\PP^*_{U}) \cap A \ne T^{\ff,\prime}_{U,\psi'} (\PP^*_{U})\cap A ) \le 
	c_{\ms{DPP}}|B_{2s/3}(A)|\exp(-c_{\ms{SPT}} s). 
\end{equation}
Moreover, for all $n \ge 1$,
$$\P(T^{\ff}_{U,\psi} (\PP^*_{U})  \sm Q_n \ne T^{\ff,\prime}_{U,\psi'} (\PP^*_{U}) \sm Q_n) \le c'_{\ms{DPP}}\exp(-c_\ms{SPT} n). $$
\enc}

Before giving the proof of Theorem \ref{th:disdis}, we mention the following result, which will be a corollary of Lemma~\ref{lem:distv} below. While not sufficient for the proof of our main theorems, we think the result is of independent interest. It  shows that we do not need to use disagreement coupling to bound the disagreement probabilities \eqref{eq:disprob}. In fact, the standard Poisson embedding has a similar property, at least if the ordering $\iota$ is such that $\iota(x)<\iota(y)$ whenever $x\in A$ and $y\in Q\sm A$ (we will, however, not be able to make such a choice of $\iota$ in the proof of our main theorems). This is somewhat surprising, since the standard Poisson embedding was not designed to control the pairwise disagreement.
 
\bec[Disagreement probabilities for standard Poisson embedding]\label{cor:poiss_emb}
Suppose that $\k, \k'$ are PIs satisfying condition \eqref{as:(A)}. Assume that $A\su P\su Q$ are such that $\k,\k'$ satisfy \eqref{eq:k,k'}. Let $\iota:Q\to [0,1)$ be an injective map such that $\iota(x)<\iota(y)$ whenever $x\in A$ and $y\in Q\sm A$. Let $\psi,\psi'\in \Nlf_{Q^c}$ and $B\su Q^c$ such that $\psi\cap B^c = \psi'\cap B^c$. Then,
\begin{equation*}
	\P\big(T_{Q,\psi}(\Pds_Q)\cap A\ne T_{Q,\psi'}(\Pds_Q) \cap A\big) \le \a_0(2+\a_0 |A|)|A| \P\big(A \lrsa B\cup  (Q\sm P) \big).
\end{equation*}
\enc

The key to proving  Theorem \ref{th:disdis} and Corollary \ref{cor:poiss_emb} is the following lemma which considers the standard Poisson embedding. In the lemma, we only need to choose the ordering $\iota$ locally, which makes the lemma applicable in the proofs of our main theorems.

%
%
\bel [Disagreement probability of first point]
\label{lem:distv}
Let $A\su P \su Q$ be bounded Borel, $\psi,\psi'\su Q^c$ be locally finite, and $\iota:Q\to [0,1)$ be injective such that $\iota(x)<\iota(y)$ whenever $x\in A$ and $y\in Q\sm A$. Suppose that $\k, \k'$ are PIs satisfying \eqref{as:(A)} and \eqref{eq:k,k'}. 
Then,
$$
	\P(\inf{}_{\iota} (T_{Q,\psi}(\Pds_Q) \cap A) \ne \inf{}_{\iota} (T'_{Q,\psi'}(\Pds_Q)\cap A))\le (2+\a_0 |A|)\,\dtv(\XX(Q,\psi)\cap A,\XX'(Q,\psi')\cap A),
$$
where 
 $\inf_\iota$ is denotes the smallest point with respect to the ordering $\le_\iota$.
\enl

\bep
In the following, we write $p(x, \psi) := p(x, Q, \psi)$.  
$$q(x) :=\P\big(T_{Q,\psi}(\Pds_Q)\cap Q_{(-\ff,x)}=\es\big)= \P\Big(U_i > p(X_i, \psi)\text{ for all }(X_i, U_i) \in \Pds_{Q_{(-\ff,x)}}\Big),$$ 
denote the probability that all points of $\Pds_Q$ before $x$ are rejected in the thinning. Similarly, we define $q'(x)$ and $q^\vee(x)$ by replacing $p(X_i,\psi)$ by $p'(X_i,\psi')$ and $p(X_i,\psi)\vee p'(X_i,\psi')$, respectively. 
By construction of the Poisson embedding, and the multivariate Mecke equation \cite[Theorem 4.4]{LP}, we obtain that the probability that the first points of $T_{Q,\psi}(\Pds_Q)$ and $T'_{Q,\psi'}(\Pds_Q)$ disagree on $A$ is given by
\begin{align}
	&\P\big(\inf{}_{\iota} \big(T_{Q,\psi}(\Pds_Q)\cap A\big) \ne \inf{}_{\iota} \big(T'_{Q,\psi'}(\Pds_Q)\cap A\big)\big)= \int_{A} |p(x,\psi) - p'(x,\psi')| q^\vee(x)\,\mathrm dx \nonumber\\
	&\le \int_{A} |p(x,\psi) - p'(x,\psi')| q(x)\,\mathrm dx \nonumber\\
	&\le \int_{A} |p(x,\psi)q(x) - p'(x,\psi')q'(x)|\,\mathrm dx+\int_{A}p'(x,\psi')|q'(x) - q(x)|\,\mathrm dx.\label{dtvbou}
\end{align}
Hence, since  $\P(\inf_\iota(T_{Q,\psi}(\Pds_Q)\cap A) \in B) =\int_{A \cap B} p(x,\psi)q(x)    \d x$, the definition of total variation distance yields,
\begin{align*}
\int_{A} |p(x,\psi)q(x) - p'(x,\psi')q'(x)|\,\mathrm dx& \le2 \sup_{B \in \BB^d} \big| \P\big(\inf{}_\iota\big(T_{Q,\psi}(\Pds_Q)\cap A\big) \in B\big) - \P\big(\inf{}_\iota\big(T_{Q,\psi'}'(\Pds_Q)\cap A\big)\in B\big)  \big|\\ &\le 2\dtv \big(T_{Q,\psi}(\Pds_Q)\cap A,T_{Q,\psi'}'(\Pds_Q)\cap A\big),
\end{align*}
where the first inequality can be seen by letting $$B:= \begin{cases}
\{x \in A:\,p(x,\psi)q(x)\ge p'(x,\psi')q'(x)\},\quad \text{if } \quad \P(\inf_\iota T_{Q,\psi}(\Pds_Q)\in A) \ge \P\big(\inf_\iota T_{Q,\psi'}'(\Pds_Q)\in A\big),\\
	\{x \in A:\,p(x,\psi)q(x)< p'(x,\psi')q'(x)\},\quad\text{if } \quad \P(\inf_\iota T_{Q,\psi}(\Pds_Q)\in A) < \P\big(\inf_\iota T_{Q,\psi'}'(\Pds_Q)\in A\big).
\end{cases}
$$ To bound the second integral in \eqref{dtvbou}, we use that  $p'(x,\psi')\le \a_0$ and obtain 
\begin{align*}
\int_{A}p'(x,\psi')|q'(x) - q(x)|\,\mathrm dx&\le \a_0 |A| \sup_{B \in \BB^d}\big|\P\big(T_{Q,\psi}(\Pds_Q)\cap A \cap B=\es\big) - \P\big(T_{Q,\psi'}'(\Pds_Q)\cap A \cap B=\es\big)\big|\\ &\le\a_0 |A| \dtv\big(T_{Q,\psi}(\Pds_Q)\cap A,T_{Q,\psi'}'(\Pds_Q)\cap A\big),
\end{align*}
which finishes the proof.
\enp

Before proving Theorem \ref{th:disdis}, we give the proof of Corollary \ref{cor:poiss_emb}, which is a slightly simpler version of the proof of Theorem \ref{th:disdis}.

%
%
\bep[Proof of Corollary \ref{cor:poiss_emb}]
Let $y_1,\ldots,y_K$ be the points in $\Pd_A$ and let $z_i = \{y_i \}\cap T_{Q,\psi}(\Pds_Q)$ and $z_i' = \{y_i\} \cap T'_{Q,\psi'}(\Pds_Q)$, $i=1,\ldots,K$ denote the thinning of these points. Then, 
\begin{align}\nonumber
	&\P\big(T_{Q,\psi}(\Pds)\cap A \ne T_{Q,\psi'}(\Pds) \cap A\big) = \sum_{k \ge 1} \P\big( z_i=z_i',\, i< k ,\, z_k \ne z_k',\, k\le K\big)\\
	&=\sum_{k \ge 1} \E\big[\one{\{ z_i=z_i',\, i< k \}} \E\big[\one{\{  z_k \ne z_k',\, k\le K\}} \mid \Pds_{Q_{(-\ff,y_{k-1})}} \big] \big].\label{eq:ce_sum}
\end{align}
In the following, we write $S_{k}:=Q_{(-\ff,y_{k})}$ and $S_{k}^*:=Q_{(-\ff,y_{k})}\times[0,\a_0]$.
Note that  by construction of the standard Poisson embedding,  
$$
T_{Q,\psi}(\Pds_{Q}) \sm S_{k-1}= T_{Q\sm S_{k-1} ,\{z_1,\ldots,z_{k-1}\}\cup \psi}(\Pds_{Q\sm S_{k-1}}),$$
and the similar statement holds for $T'_{Q,\psi'}(\Pds_{Q})$. Moreover, since $S_{k-1}^*$ is a stopping set, conditionally on $\Pds_{S_{k-1}}$,  $\Pds_{Q\sm S_{k-1}}$ has the distribution of a Poisson process on $Q^*\sm S_{k-1}^*$. Thus, we are in the situation of Lemma \ref{lem:distv} with $Q$ replaced by $Q \sm S_{k-1} $ and $A$ replaced by $A\sm S_{k-1}$. Therefore, we may bound the conditional expectation in \eqref{eq:ce_sum} by
$$(2+\a_0 |A|)  \E\big[ \dtv\big(\XX(Q \sm S_{k-1}, \{z_1,\ldots,z_{k-1}\}\cup \psi)\cap A, \XX'(Q \sm S_{k-1}, \{z_1',\ldots,z_{k-1}'\}\cup \psi' )\cap A\big) \mid \Pds_{S_{k-1}} \big]. $$
Moreover, on the event $\{z_i =z_i',i<k\}$, the boundary conditions differ only on $B$, so the total variation distance is bounded by $\P(A \sm S_{k-1} \lrsa B\cup Q\sm (P\cup S_{k-1} ) )$ by Corollary \ref{cor:kappa'}. Since the number of terms in \eqref{eq:ce_sum} is $K=\Pd(A)$, we obtain the bound on \eqref{eq:ce_sum}
$$(2+\a_0 |A|)\P\big(A \lrsa B\cup (Q\sm P) \big) \E[\Pd(A)] = (2+\a_0 |A|)\a_0|A| \P\big(A \lrsa B\cup (Q\sm P) \big) .$$
\enp

%
%
	\bep[Proof of Theorem \ref{th:disdis}]
Throughout this proof, we let $\XX:=T^{{\sf{rad}}}_{Q,\psi}(\Pds_Q)$ and $\XX':=T^{'\sf{rad}}_{Q,\psi'}(\Pds_Q)$. 
The idea is to bound the disagreement probability of $\XX \cap A$ and $\XX' \cap A$ by considering disagreement on the components of $B_{r_0/2}(\Pd)$ that intersect $A$. More precisely, let  $V_i$ be as in  Example \ref{ex:dc2}. Let  $V_{i_1}, V_{i_2}, \dots, V_{i_M}$ denote the subsequence of the sets $V_i$ such that the starting point $v_{i,1}\in D$ is contained in $B_s(A)$.  We introduce the events
$$ E_s :=\{A \lrsa B_s(A)^c\}^c, \qquad \tilde{E}_s:=\{B_s(A) \lrsa B_{2s}(A)^c\}^c, \qquad \bar E_s := E_s \cap \tilde E_s. $$
Then, on the event $E_s$, $\XX$ and $\XX'$ agree on $A$ if they agree on all the $V_{i_m}$. Moreover, on the event $\tilde{E}_s$, all $V_{i_m}$ are contained in $B_{2s}(A)$. 
From \eqref{eq:conn_prob}, we have
	$$
	\P(E_s^c)\hspace{-0.05cm} =\hspace{-0.05cm} \P(A  \lrsa B_s(A)^c)\hspace{-0.05cm}\le \hspace{-0.05cm} c_1\big|B_{r_0}(A) \big| e^{-c_{\ms{SPT}, 2}s},\quad \P(\tilde{E}_{s}^c)\hspace{-0.05cm}=\hspace{-0.05cm}\P( B_{s}(A)  \lrsa B_{2s}(A)^c)\hspace{-0.05cm}\le\hspace{-0.05cm} c_1\big|B_{s+r_0}(A)\big| e^{-c_{\ms{SPT}, 2}s}.
	$$
	Thus, it is enough to work under the events $E_s$ and $\tilde{E}_s$.
	
Recall the definitions of $V_{i, j}$ and  $\xi_{i,j}$ from Example \ref{ex:dc2} and let $\XX_i:=\XX \cap V_i$, $\XX'_i:=\XX' \cap V_i$ and  $\XX_{i,j}:=\xi_{i,1} \cup \dots \cup \xi_{i,j} $, $\XX_{i,j}':=\xi_{i,1}' \cup \dots \cup \xi_{i,j}'$ for $ j \le N_{i}$.  Then,
	\begin{align*}
	&\P(\bar E_s, \XX \cap A\ne \XX' \cap A)\\
	&=\sum_{m \ge 1} \P\big(\tilde{E}_{s}, (\XX_{i_1},\dots,\XX_{i_{m-1}}) = (\XX'_{i_1},\dots,\XX'_{i_{m-1}}) , \XX_{i_m}\ne \XX'_{i_m}, M \ge m\}\big)\\
	&= \sum_{m\ge 1}\E\left[\mathds 1\{(\XX_{i_1},\dots,\XX_{i_{m-1}}) = (\XX'_{i_1},\dots,\XX'_{i_{m-1}}),M\ge m\}\P\big(\tilde{E}_{s}, \XX_{i_{m}}\ne \XX'_{i_{m}} \mid \Pds_{ S_{i_{m}-1,N_{i_m-1}}}\big)\right].
	\end{align*}
	We note here that $S_{i_{m}-1,N_{i_m-1}}^*$ is again a stopping set and the event $\{M\ge m\}$ is measurable with respect to $\Pds_{ S_{i_{m}-1,N_{i_m-1}}} $.

	We now claim that 
	\begin{align}
		\label{eq:pr14}
		\P\big(\tilde{E}_{s}, \XX_{i_{m}}\ne \XX'_{i_{m}}  \mid \Pds_{ S_{i_{m}-1,N_{i_m-1}}}\big)\le  c_1e^{-c_{\ms{SPT}, 2}s} \E\Big[\sum_{j \le N_{i_m}}|V_{i_m,j}|\one\{V_{i_m,j}\su B_{2s}(A)\}\Big],
	\end{align}
	for some $c_1 > 0$. We first explain how to conclude the proof of the Proposition from \eqref{eq:pr14}. Afterwards, we establish \eqref{eq:pr14}.

	From \eqref{eq:pr14}, we get 
	\begin{align}
		&\P\big(\bar E_s, \XX \cap A\ne \XX' \cap A\big) \le c_1 e^{-c_{\ms{SPT}, 2}s}  \E\Big[\sum_{m\le M} \sum_{j\le N_{i_m}}|V_{i_m,j}|\one{\{ V_{i_m,j}\su B_{2s}(A)\}}\Big].\label{bou:XneqX}
	\end{align}
	Next, we bound the double sum in \eqref{bou:XneqX} by
	\begin{align}
		\nonumber
		\sum_{m\le M} \sum_{j\le N_{i_m}} \big|V_{i_m,j}\big|\one{\{V_{i_m,j}\su B_{2s}(A)\}}
		\le 	& \sum_{m\le M} \sum_{j\le N_{i_m}} \big|V_{i_m,j}\big| \one{\{\Pd_{ V_{i_m,j}}=\es, V_{i_m,j}\su B_{2s}(A)\}}\\
		&+ \sum_{m\le M} \sum_{j\le N_{i_m}} \big|V_{i_m,j}\big| \one{\{\Pd_{V_{i_m,j}}\ne \es, V_{i_m,j}\su B_{2s}(A)\}}.\label{bouXX1}
		\end{align}
	Note that if $\Pd_{V_{i_m,j}}=\es$, then $V_{i_m,j}=Z_{i_m,j}$. Since the $Z_{i_m,j}$ are disjoint and contained in $B_{2s}(A)$,  the expected value of the first term on the right hand side of \eqref{bouXX1} is bounded by $\big|B_{2s}(A)\big|$. For the expected value of the second term on the right hand side of  \eqref{bouXX1}, note that every point in $\Pd$ can only be contained in one  of the sets $Z_{i_m,j}\su V_{i_m,j}$. Thus the number of terms in the sum is bounded by $\PP(B_{2s}(A))$. By Proposition \ref{pr:rad_prop}, $|V_{i_m,j}|\le |B_{r_0}|$. Hence, we conclude that the second term on the right hand side of \eqref{bouXX1} is bounded by
	$$
|B_{r_0}| \E[\Pd (B_{2s}(A))]=|B_{r_0}|\a_0 \big|B_{2s}(A)\big|.
	$$
	These estimates show that for some constant $c_2>0$,  \eqref{bou:XneqX} is bounded by
	\begin{align*}
		c_1(1+\a_0\big|B_{r_0}\big|)\big|B_{2s}(A)\big| e^{-c_{\ms{SPT}, 2}s} \le c_2\big|B_{2s}(A)\big| e^{-c_{\ms{SPT}, 2}s}.
	\end{align*}

	It remains to prove \eqref{eq:pr14}. 
To that end, note that 	given $\Pds_{ S_{i_m-1,N_{i_m-1}}}$, the processes $\XX_{i_{m}}$ and $\XX'_{i_{m}}$ are restrictions of Gibbs processes on $Q\sm S_{i_m-1,N_{i_m-1}}$  with boundary condition $ \psi$ and $\psi'$, respectively, to $V_{i_m}$ (since $B_{r_0}(\Pd\cap S_{i_m-1,N_{i_m-1}})\cap Q \su S_{i_m-1,N_{i_m-1}}$). Hence, it is enough to consider the case $i_m = 1$. Here, we obtain the bound on \eqref{eq:pr14}
	\begin{align}
	&\P\Big(\{V_1\su B_{2s}(A) \}\cap \Big\{\XX \cap \bigcup_{j\le N_1}V_{1,j}\ne \XX'\cap \bigcup_{j\le N_1} V_{1,j} \Big\}\Big)\nonumber\\
	&\, \le	\sum_{j \ge 1}	\P(\XX_{1,j-1}=\XX'_{1,j-1}, \xi_{1,j}\ne \xi'_{1,j},N_1\ge j, V_{1,j}\su B_{2s}(A))\nonumber\\
		&\,=\E\big[\sum_{j \le N_1}\one \{\XX_{1,j-1}=\XX'_{1,j-1}, V_{1,j}\su B_{2s}(A)\}\, \P( \xi_{1,j}\ne  \xi'_{1,j} \mid \Pd_{S_{1,j-1})}\big],\label{eq:Y1n}
	\end{align}
	where $\{\XX_{1,0}=\XX'_{1,0}\}$ is interpreted as trivially satisfied. Here, we have used that both the set $V_{1, j}$ and the event $\{N_1\ge j\}$ are measurable with respect to  $\Pd_{ S_{1,j-1}}$. 
		Since $S_{1,j-1}^*$ is a stopping set, conditionally on $\Pds_{ S_{1,j-1}^*}$, we have that $\xi_{1,j}$ is the restriction of a Gibbs process 
	$$
	\tilde{\XX}_{1,j}:=T_{Q\sm S_{1,j-1}, \XX_{1,j-1} \cup \psi}(\Pds_{Q\sm S_{1,j-1}})
	$$
	and analogously for $\xi'_{1,j}$. Recalling that $\xi_{1, j}$ is $ \inf_{\iota_{1,j}} (\Pd \cap V_{1,j})\cap \tilde{\XX}_{1,j}   $ (with $\inf(\es)$ interpreted as $\es$), the event $\{ \xi_{1,j}\ne  \xi'_{1,j} \}$ is contained in the event $\{\inf_{\iota_{1,j}} (\tilde{\XX}_{1,j} \cap V_{1,j})\ne \inf_{\iota_{1,j}} (\tilde{\XX}_{1,j}'\cap V_{1,j}) \}$.  Thus, we may apply Lemma \ref{lem:distv} with $Q$ replaced by $Q\sm S_{1,j-1}$ and $A$ replaced by $V_{1,j}$.   Since $|V_{1,j}|\le \big|B_{r_0}\big|$,  Lemma \ref{lem:distv} shows that the conditional probability in \eqref{eq:Y1n} is bounded by
	\begin{align}\label{eq:abc}
	c_3  \dtv\big(\XX(Q\sm S_{1,j-1}, \XX_{1,j-1} \cup \psi)\cap V_{1,j},  \XX'(Q\sm S_{1,j-1},\XX'_{1,j-1} \cup \psi')\cap V_{1,j}\big),
	\end{align}
	where $c_3 := 2+\a_0\big|B_{r_0}\big|$.	Using here that $\psi, \psi'\su Q^c\su P^c$ and that we are on the event $\{\XX_{1,j-1}=\XX'_{1,j-1}\}$, Corollary \ref{cor:kappa'} with $A$ replaced by $ V_{1, j}$ \cb{and  $Q$ replaced by $Q\sm S_{1,j-1}$ gives the bound $\P(V_{1,j} \lrsa (Q\sm P)\cup B)$} on the total variation distance in \eqref{eq:abc}, which is again bounded by 
	$
	c_1\big|B_{r_0}( V_{1,j})\big| e^{-c_{\ms{SPT}, 2}s},
	$ by the argument in \eqref{eq:conn_prob} using that $s\le\dist(B_{2s}(A), P^c)$ by assumption.
	Letting $c_4 := c_3c_1$, we use this finding  in \eqref{eq:Y1n} to obtain the asserted bound in \eqref{eq:pr14} 
	\begin{equation}
		\P\Big(E_s, \XX \cap V_1\ne \XX'\cap  V_1 \Big)\le 	c_4e^{-c_{\ms{SPT}, 2}s} \E\Big[\sum_{j\ge 1}|V_{1,j}|\one{\{N_1\ge j, V_{1,j}\su B_{2s}(A)\}}\Big]. \label{eq:W1lim}
	\end{equation}
\enp

\cb{
	\bep[Proof of Corollary \ref{cor:ff_disag}] We may assume that $A, B ,U\sm P \su Q_m $ for some $m$.  
	From Proposition \ref{pr:xxff}, we have  a $c_1>0$ such that
	\begin{align*}
	\P( T^{\ff}_{U,\psi} (\PP^*_{U}) \cap A  \ne \Tra_{Q_{m+m'}\cap U,\psi } (\PP_{Q_{m+m'}\cap U}^*) \cap A )& \le \P( Q_{m} \lrsa Q_{m+m'-2r_0}^c )\\&\le 2c_1|Q_{m+r_0}|\exp(-c_\ms{SPT} (m+m')).
	 \end{align*}
	A similar statement holds for $T^{\ff,\prime}_{U,\psi'} (\PP^*_{U})$.	By Theorem \ref{th:disdis}, 
	\begin{align*}
		\P\big(T^{\ff}_{Q_{m+m'}\cap U,\psi} (\PP^*_{Q_{m+m' }\cap U}) \cap A  & \ne T^{\ff,\prime}_{Q_{m+m'}\cap U,\psi'} (\PP^*_{Q_{m+m'}\cap U}) \cap A \big)  \le c_{\ms{DP}}|B_{2s/3}(A)|\exp(-c_\ms{SPT} s).
	\end{align*}
	Letting $m'\to \infty$, we obtain \eqref{eq:ff-disdis}.
	
	Still assuming $B,U\sm P \su Q_m$, take $n\ge m$, and $A=Q_{n+1}\sm Q_n$. Then $s\ge n-m$ and \eqref{eq:ff-disdis} shows that there must be a $c_2>0$ such that
	$$\P(T^{\ff}_{U,\psi} (\PP^*_{U})  \cap (Q_{n+1}\sm Q_{n}) \ne T^{\ff,\prime}_{U,\psi'} (\PP^*_{U}) \cap (Q_{n+1}\sm Q_n)) \le c_2\exp(-c_\ms{SPT} n). $$
	Summing these probabilities, we find that 
	$$\P(T^{\ff}_{U,\psi} (\PP^*_{U})  \sm Q_n \ne T^{\ff,\prime}_{U,\psi'} (\PP^*_{U}) \sm Q_n) \le c_3\exp(-c_\ms{SPT} n), $$
	which shows the claim.
	\enp}

	\section{Proof of Theorem \ref{thm:1} -- CLT under weak stabilization}
\label{sec:weak}

In this section, we prove Theorem \ref{thm:1}, i.e., the CLT for Gibbsian functionals under weak stabilization. The main idea is to extend the martingale approach that was developed for Poisson point processes in \cite[Theorem 3.1]{penrose}. A martingale is obtained by conditioning on the point process in an increasing sequence of windows. The centered functional can then be expressed as a sum of martingale differences. Thus, general CLT results apply, if these differences can be controlled. Since the martingale differences are essentially obtained by locally changing the boundary conditions, they can be handled using the disagreement coupling results from Section \ref{sec:const}. 

The rest of the section is organized as follows. First, in Section \ref{sec:gen}, we establish a CLT under a general moment and stabilization condition. Second, in Section \ref{sec:clt}, we prove the main CLT for Gibbs processes under general conditions. In Section \ref{sss:moment}, we discuss when these conditions are satisfied. Finally, in Section \ref{sss:vp}, we demonstrate how to show variance positivity in the setting of Betti numbers, thereby completing the proof of Corollary \ref{cor:betti}. 

%
%
\subsection{CLT under general conditions}
\label{sec:gen}
In the following, we let $(a_n)$ be a sequence such that $\lim_{n\to \infty}a_n= \infty$.  Let $\YY_{a_n}$, $n\geq 1$, denote a sequence of point processes on $\R^d$. Later in this section, we always consider the specific choices $a_n=n+\sqrt n$ and $\YY_{a_n} = \XX$. However, since it causes no additional work, we formulate the proofs in the general setting. For instance, in other applications, one might want to take $\YY_{a_n} = \XX(Q_{a_n},\psi_n)$. For $n\ge1$, we consider functionals $H_n(\vp)$ defined for any $\vp\in \Nlf$. In the present subsection, we do not require translation-invariance of $H_n$ and allow $H_n$ to depend on $n$.

Let $Z_{a_n}$ be the set of lattice points $z\in \Z^d$ such that $Q_{z,1}$ intersects the window $Q_{a_n}$. Then, $Q_{a_n}\su \bigcup_{z\in Z_{a_n}} Q_{z,1}$. Order $Z_{a_n}$ lexicographically as $z_1,\dots,z_{k_n}$. Let $\FF_{0,n}$ be the trivial $\s$-algebra and $\FF_{k_n,n}=\sigma(\YY_{a_n})$. Moreover, define 
$$ \FF_{i,n}:= \s\bigg(\YY_{a_n} \cap \bigcup_{{z\in \Z^d, z\preceq z_i}} Q_{z,1}\bigg)$$
 for $0<i<k_n$, where $\preceq$ denotes lexicographical ordering on $\Z^d$. 

The key observation in the martingale approach is that 
\begin{equation*}
	H_n(\YY_{a_n}) - \E[H_n(\YY_{a_n})]= \sum_{i\le k_n} \De_{i,n}, 
\end{equation*}
where 
$$\De_{i,n} := \De_{z_i,n}:= \E[H_n(\YY_{a_n}) | \FF_{i,n}] - \E[H_n(\YY_{a_n}) | \FF_{i-1,n}]. $$
For each $n$, the $ \E[H_n(\YY_{a_n}) | \FF_{i,n}]$ define a martingale with respect to $\FF_{i,n} $. By orthogonality of martingale differences, the variance is given by 
\begin{equation}\label{eq:var_mart_diff}
	\Var(H(\YY_{a_n}) - \E[H(\YY_{a_n})])= \sum_{1\le i\le k_n}\E[  \De_{i,n}^2]. 
\end{equation}
We impose the following two conditions:
\begin{itemize}
	\item[(i)] 	$
		\sup_{n\ge 1, z\in Z_{a_n}} \E[\De_{z,n}^4] < \infty.
		$
	\item[(ii)] There exists a stationary ergodic point process $\YY$ such that for every $z\in \Z^d$, there is a random variable $\De_z = \De_z(\YY)$ such that ${\De}_ {z,n} \to \De_z$ in probability (and hence in $L^2$ due to (i)). The limit is shift invariant meaning that for any $z_0\in \Z^d$, $\De_{z+z_0} (\YY + z_0) = \De_ z(\YY)$. The convergence must be uniform in the sense that	
	\begin{equation}\label{eq:unif_conv}
		\lim_{n\to \infty} \sup_{z\in Z_{a_n-\rho(a_n)}}\lVert\De_z-\De_{z,n}\rVert_{L^2}=0,
	\end{equation}
	where $\rho:(0,\infty) \to \R$ is a function with $\rho(x)\in o(x)$ and $\lim_{x\to \infty} \rho(x) = \infty$.
\end{itemize}

The proposition below and its proof is an adaption of \cite[Thm. 3.1]{penrose} to general point processes.

\bepr[CLT for general functionals]
\label{thm:martingale}
Let $(a_n)$ be a sequence with $\lim_{n\to \infty} a_n = \infty$ and $\YY_{a_n}$ a sequence of point processes on $\R^d$. 	Under Conditions (i) and (ii), 
\begin{equation}
	a_n^{-d}\Var(H_n(\YY_{a_n} )) \to \s^2\ge 0, \label{eq:lim_var}
\end{equation}
where $\s^2 = \E[\De_0^2]$, and 
\begin{equation*}
	a_n^{-d/2}(H_n(\YY_{a_n} ) - \E[H_n(\YY_{a_n} )] ) \to N(0,\s^2)
\end{equation*}
in distribution.
\enpr

\bep
Since $\lim_{n\to \infty}a_n^d/k_n=1$, Proposition \ref{thm:martingale} follows from the martingale CLT of \cite{mcleish} if the following three conditions are satisfied: 
\begin{itemize}
	\item[1.] $\sup_{n \ge 1} k_n^{-1}\E\big[ \max_{i \le k_n} \De_{i,n}^2\big] < \ff$;
	\item[2.] $k_n^{-1/2}\max_{i\le k_n}|\De_{i,n}| \to 0$ in probability; and
	\item[3.]$k_n^{-1}\sum_{i \le k_n} \De_{i,n}^2 \to \s^2$ in $L^1$.
\end{itemize}
In particular, \eqref{eq:lim_var} follows from \eqref{eq:var_mart_diff} and 3. The first two conditions follow from Condition (i) exactly as in \cite[proof of Thm. 3.1]{penrose}. To make our presentation self-contained, we briefly reproduce the argument.

\begin{itemize}
	\item[1.] Condition (i) implies
	\begin{equation*}
		\sup_{ n\ge1} k_n^{-1}\E \big[ \max_i \De_{i,n}^2\big] \le \sup_ n k_n^{-1 } \sum_{i\le k_n} \E [ \De_{i,n}^2]<\infty.
	\end{equation*}
	
	\item[2.] This follows from Markov's inequality and (i) because
	\begin{equation*}
		\P\big(k_n^{-1/2}\max_{i\le k_n}|\De_{i,n}| \ge \eps\big) \le \sum_{i\le k_n}\P\big(k_n^{-1/2}|\De_{i,n}| \ge \eps\big) \le k_n \frac{\max_{i\le k_n}\E[\De_{i,n}^4]}{k_n^ 2\eps^2}, 
	\end{equation*}
	which tends to 0 as $n\to \infty$.
	
	\item[3.] We start by noting that $\De_z$ is in $L^2$, since $\De_{z,n}$ converges to $\De_z$ in probability, and (i) implies uniform integrability.
	
	By Assumption (ii), ${\De_z(\YY)}_{z\in \Z^d}$ forms a multidimensional ergodic sequence, meaning that the sum 
	\begin{equation*}
		|Z_{a_n-\rho(a_n)}|^{-1}\sum_{z\in Z_{a_n-\rho(a_n)}} \De_ z^2 \xrightarrow{L^1} \s^2,
	\end{equation*}
	see \cite[Thm. 10.12]{kallenberg}. 
	Moreover, since $|Z_{a_n-\rho(a_n)}| = a_n^d + o(a_n^{d}) = k_n+o(k_n)$ and $\De_z\in L^2$,
	\begin{equation*}
		(k_n^{-1}-|Z_{a_n-\rho(a_n)}|^{-1})\sum_{z\in Z_{a_n-\rho(a_n)} } \De_ z^2 \xrightarrow{L^1} 0,
		\quad \text{	and }\quad
		k_n^{-1}\sum_{z\in Z_{a_n}\sm Z_{a_n-\rho(a_n)} } \De_ z^2 \xrightarrow{L^1} 0.
	\end{equation*}
	
	It remains to show that 
	\begin{equation*}
		k_n^{-1}\sum_{z\in Z_{a_n}} |\De_ z^2 - \De_ {z,n}^2| \xrightarrow{L^1} 0.
	\end{equation*}
	For this, note that
	\begin{equation*}
		\E \big[|\De_ z^2 - \De_{z,n}^2|\big] \le \E \big[(\De_ z + \De_{z,n})^2\big]^{1/2} \E\big[ (\De_ z - \De_{z,n})^2\big]^{1/2}. 
	\end{equation*}
	The first term is uniformly bounded by (i) and the fact that $\De_ z \in L^2$.
	From Assumption (ii) we have that
	\begin{equation*}
		k_n^{-1} \sum_{z\in Z_{a_n-\rho(a_n)}} \E\big[ (\De_ z - \De_{z,n})^2\big] \to 0 
	\end{equation*}
	as $n \to 0$. Finally,
	\begin{equation*}
		k_n^{-1} \sum_{z\in Z_{a_n}\sm Z_{a_n-\rho(a_n)}} \E \big[(\De_ z - \De_{z,n})^2\big] \to 0, 
	\end{equation*}
		again by (i) and the fact that $\De_z\in L^2$. This shows 3.
\end{itemize}

\enp

%
%
\subsection{A CLT for translation-invariant functionals of Gibbs point processes}
\label{sec:clt}
We now return to Gibbs point processes. 
The present subsection is devoted to connecting the conditions (i) and (ii) stated in the general CLT in Proposition \ref{thm:martingale} with the conditions \eqref{eq:moment_cond} and \eqref{eq:as_stab} from Section \ref{sec:mod} in the Gibbs setting. 
We consider a functional $H_n$ of the form 
\begin{equation}\label{eq:H_n_Gibbs}
	H_n(\vp) = H(\vp \cap Q_n) 
\end{equation}
where $\vp \in \Nlf$ and $H$ is a translation-invariant functional on $\Nlf_0$.

Throughout this section we let $\X=T^\ff(\Pds)$ denote the infinite-volume Gibbs process on $\R^d$. 
We define sets
\begin{equation*}
V_z = \bigcup_{w\preceq z} Q_{w,1},\quad U_z = \bigcup_{w \succ z} Q_{w,1},
\end{equation*}
where $\preceq$ refers to the lexicographic ordering on $\Z^d$, and define $V_{z,n}:=V_z \cap Q_n$ and $U_{z,n}:=U_z\cap Q_n$. \begin{wrapfigure}{r}{0.515\textwidth}
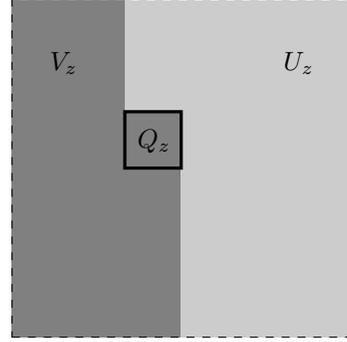

	\centering
	\betp[scale=1.5]

\fill[black!50] (-1.5, -2.0) rectangle (0, 1);
\fill[black!20] (-.5, 0) rectangle (1.5, 1);
\fill[black!20] (0, -2.0) rectangle (1.5, 1);

\draw[dashed] (-1.5, -2) rectangle (1.5, 1);

\coordinate[label={$Q_z$}] (A) at (-.24, -0.45);
\coordinate[label={$V_z$}] (A) at (-1.04, 0.25);
\coordinate[label={$U_z$}] (A) at (1.04, 0.25);

\draw[very thick] (-.5, -.5) rectangle (0, 0);
\entp
	\caption{Illustration of the sets $U_z$ and $L_z$.}
	\label{fig:fig1}
	\vspace{-.5cm}
\end{wrapfigure}
We refer the reader to Figure \ref{fig:fig1} for an illustration of these sets.

	We aim at proving the following CLT under the conditions (i) and \eqref{eq:as_stab}. This will be the most general version of the CLT. We discuss further in which situations the conditions are satisfied in Subsection \ref{sss:moment}. In the following, let $z\in \Z^d$ and let $z_{-}$ denote the point in $\Z^d$ that comes just before $z$ in the lexicographic order.

\bet[CLT for translation-invariant functionals]
\label{thm:gibbs_CLT}
Let $\XX$ be an infinite-volume Gibbs point process on $\R^d$
with translation-invariant PI satisfying \eqref{as:(A)}. Let $H$ be a translation-invariant functional of finite point patterns such that \eqref{eq:as_stab} holds. Let $H_n$ be as in \eqref{eq:H_n_Gibbs}. Suppose, \begin{equation}\label{eq:moment(i)}
\sup_{n\ge 1, z_i\in Z_{n+\sqrt n}} \E[(\E[H_n(\XX) | \FF_{i,n}] - \E[H_n(\XX) |\FF_{i-1,n}] )^4] < \infty.
\end{equation}
Then, the limit
\begin{equation*}
	\De_0 = \lim_{n\to \infty } \left(\E[ H_n(\X )| \X\cap V_ 0 ] - \E[ H_n(\X ) | \X\cap V_{0_-}]\right)
\end{equation*}
exists in $L^2$. Let $\s^2:= \E[\De_0^2]$. Then,
$
|Q_n|^{-1}\Var(H_n(\X) ) \to \s^2
$
and
\begin{equation*}
	|Q_n|^{-1/2}\left(H_n(\X) - \E[H_n (\X)] \right) \to N(0,\s^2).
\end{equation*}
\ent

For the proof, we make the following constructions. Define
$$\eta_z:=T^{\ff} (\Pds)\cap V_z, \quad \XX_z' = \eta_z \cup T^{\ff,z}_{U_z,\eta_z}(\tilde{\PP}^*_{U_z}),$$
where $\tilde{\PP}^*$ is an independent copy of $\Pds$ and $T^{\ff,z}$ denotes the radial coupling centered at $z$, i.e.\ $T^{\ff,z}(\vp)=T^{\ff}(\vp-z)+z$. Then $\XX_z'\sim \XX$ by Proposition \ref{pr:chimera}. Note that
$$ \E[H_n(\XX) | \eta_{z_-} ]= \E[H_n(\XX_{z_-}') | \eta_{z_-} ] = \E[H_n(\XX_{z_-}') | \eta_z ].$$
Thus, we can write the martingale differences as 
$$\Delta_{z,n} = \E[H_n(\XX_z') \mid \eta_z ]-\E[H_n(\XX_{z_{-}}') \mid \eta_{z_-}] = \E[H_n(\XX_z') - H_n(\XX_{z_{-}}') \mid \eta_z ].$$
Note that $\XX_z'$ and $\XX_{z_-}'$ both equal $\eta_{z_-}$ on $V_{z_-}$. They are extended from $\eta_z$ and $\eta_{z_-}$ on $U_z$ and $U_{z_-}$, respectively. Lemma \ref{lem:ff_disag} below will allow us to control the disagreement this causes in regions far from $z$. 

\bel\label{lem:ff_disag}
Suppose that the PI satisfies \eqref{as:(A)}. There exists a random variable $N_z \ge 1$ with $\P(N_z\ge n)\le c_{\ms{DPP}}'\exp(-c_\ms{SPT}n)$ such that almost surely 
$$ \XX'_z\cap (U_z\sm Q_{z,N_z}) = \XX'_{z_-}\cap (U_{z_-}\sm Q_{z,N_z}). $$
\enl

\cb{
\bep
It is enough to consider the case $z=0$. 
Note that $T^{\ff}_{U_0,\eta_0}(\PP^*_{U_0})=T^{\ff,\prime}_{U_{0_-},\eta_{0_-}}(\PP^*_{U_0}) $,
where the thinning operator $T^{\ff,\prime}$ makes use of the PI
$$\kappa'(x,\mu ) = \begin{cases} \kappa(x,\mu \cup (\eta_0\sm \eta_{0_-})),& x\in U_{0},\\
	0, & x\in Q_{0,1}.
\end{cases}$$
By Corollary \ref{cor:ff_disag}, 
\begin{equation*}
	\P\big( T^{\ff}_{U_{0_-},\eta_{0_-}}(\PP^*_{U_{0_-}})\sm Q_n\ne T^{\ff,\prime}_{U_{0_-},\eta_{0_-}}(\PP^*_{U_0})\sm Q_n \big) \le c_{\ms{DPP}}'\exp(-c_\ms{SPT}n).
\end{equation*}
Letting $N$ be the minimal $n$ such that $ T^{\ff}_{U_{0_-},\eta_{0_-}}(\PP^*_{U_{0_-}})\sm Q_n\ne T^{\ff,\prime}_{U_{0_-},\eta_{0_-}}(\PP^*_{U_0})\sm Q_n$ shows the claim.
\enp
}

\bep[Proof of Theorem \ref{thm:gibbs_CLT}]
We show Condition (i) and (ii) of Proposition \ref{thm:martingale} with $a_n=n+\sqrt n$, $\YY_{a_n}=\XX$ and $\rho(a_n)=2\sqrt n$. Assumption \eqref{eq:moment(i)} immediately implies (i). It remains to show (ii).
Fix $z\in \Z^d$. By Lemma \ref{lem:ff_disag}, there is a.s.\ an $N_z$ such that $\XX_z'\sm Q_{z,N_z}=\XX_{z_{-}}'\sm Q_{z,N_z}$. Hence,
\begin{align*}
	 H(\XX_z') - H(\XX_{z_{-}}') &:=\lim_{n\to \infty} H_n(\XX_z') - H_n(\XX_{z_{-}}')\\
	& =\lim_{n\to \infty} H_n(\XX_z') - H_n(\XX_z'\sm Q_{z,N_z}) - H_n(\XX_{z_-}') + H_n(\XX_{z_{-}}'\sm Q_{z,N_z})
\end{align*}
exists almost surely by \eqref{eq:as_stab}. Note that
$$\E[\E[H(\XX_z') - H(\XX_{z_{-}}')|\eta_z]^4] \le M:=\sup_{n\ge 1, z\in\Z^d\cap Q_n}\E[\Delta_{z,n}^4]$$
by the conditional Fatou's Lemma. 

 We now show that 
$$\E[H_n(\XX_z') - H_n(\XX_{z_{-}}') \mid \eta_z ] \to \E[H(\XX_z') - H(\XX_{z_{-}}')\mid \eta_z] =:\Delta_z \text{ in } L^2 .$$

First, note that \eqref{eq:as_stab} implies that for any $z\in \Z^d$, $l\ge 1$, and $\eps>0$, there exists a random variable $R(z,l,\eps)$ such that whenever $Q_{z,R(z,l,\eps)} \su Q_{w,m}$,
$$|(H(\XX)-H(\XX\sm Q_{z,l})) -( H(\XX\cap Q_{w,m})-H(\XX\cap Q_{w,m} \sm Q_{z,l}))|<\eps.$$
Indeed, otherwise, for $z=0$, there would be a $w_n\in \Z^d, m_n\in \N $ with $Q_n\su Q_{w_n,m_n}$ for any $n\in \N$, such that
$$|(H(\XX)-H(\XX\sm Q_{z,l})) -( H(\XX\cap Q_{w_n,m_n})-H(\XX\cap Q_{w_n,m_n}\sm Q_{z,l}))|\ge \eps.$$
Since $\bigcup_{k\ge 1} \bigcap_{n\ge k}Q_{w_n,m_n}=\R^d$, this contradicts \eqref{eq:as_stab}. The claim follows for any other $z$ by stationarity.) Clearly, the distribution of $R(z,l,\eps)$ is independent of $z$.

Let $\eps>0$ be given and fix $z\in \Z^d$. For $l,n\in \N$, we define events
$$E_{z,l} :=\{\XX_z'\sm Q_{z,l} =\XX_{z_{-}}'\sm Q_{z,l} \} , \quad F_{z,l,n} :=\{R(z,l,\eps)\le \sqrt n-l-1\} $$
Choose $l$ such that $\P(N_z>l)< \eps$ where $N_z$ is the random variable from Lemma \ref{lem:ff_disag}. Then
$\P(E_{z,l})> 1-\eps$. Moreover, let $n \ge l$ be such that
 $\P(F_{z,n,l}) \ge 1-\eps$. Then, for all $m \ge n$, the conditional Jensen inequality and the Cauchy-Schwarz inequality yield
\begin{align*}
	&\| \E[H_m(\XX_z') - H_m(\XX_{z_{-}}') \mid \eta_z ]- \E[H(\XX_z') - H(\XX_{z_{-}}') \mid \eta_z ]\|_{L^2}^2\\ &\le 
	\| \one_{E_{z,l} \cap F_{z,l,n} \cap F_{z_-,l,n}}((H_m(\XX_z') - H_m(\XX_{z_{-}}')) -( H(\XX_z') - H(\XX_{z_{-}}') ))\|^2_{L^2}\\ &\quad +4\E[\Delta_{z,m}^4 + \Delta_z^4]^{1/2} (\P(E_{z,l}^c) + \P(F_{z,l,n}^c)+\P(F_{z_-,l,n}^c))^{1/2}\\
	&\le 4\|\one_{F_{z,n,l}} (H_m(\XX_z') - H_m(\XX_z'\sm Q_{z,l}) - H(\XX_z') + H(\XX_z'\sm Q_{z,l}) )\|^2_{L^2} \\
	&\quad+4\|\one_{F_{z_-,n,l}}( H_m(\XX_{z_{-}}') -H_m(\XX_{z_{-}}'\sm Q_{z,l}) - H(\XX_{z_{-}}') + H(\XX_{z_{-}}'\sm Q_{z,l}) )\|^2_{L^2} +4\sqrt{6} M^{1/2} \eps^{1/2}\\
	&\le 8\eps^2 +4\sqrt{6} M^{1/2} \eps^{1/2}
\end{align*}
Since the distributions of $N_z$ and $R(z,\eps,l)$ are independent of $z$, the bound is uniform for all $z\in Q_{m-\sqrt m}$. Since $\eps$ was arbitrary, this shows \eqref{eq:unif_conv}. The fact that $\Delta_{z+z_0}(\XX+z_0) = \Delta_z(\XX) $ follows immediately by construction, and $\XX$ was stationary and ergodic by Proposition \ref{pr:xxff}.
\enp

\subsubsection{More on the conditions of Theorem \ref{thm:gibbs_CLT}}\label{sss:moment}

In this section, we give some conditions under which the assumptions \eqref{eq:as_stab} and \eqref{eq:moment(i)} in Proposition~ \ref{thm:gibbs_CLT} are satisfied. The first lemma together with Theorem \ref{thm:gibbs_CLT} implies Theorem \ref{thm:1}. 

\bel[Moment condition]
	\label{lem:(i)'}
	Let $\XX$ be an infinite-volume Gibbs point process
	with translation-invariant PI satisfying \eqref{as:(A)}.
	Consider a translation-invariant functional $H$ and let $H_n$ be defined as in \eqref{eq:H_n_Gibbs}. Then, \eqref{eq:moment_cond} implies \eqref{eq:moment(i)}. 
\enl

%
%
\bep[Proof of Lemma \ref{lem:(i)'}]
First consider the case $z_{i-1}=z_{i-}$. Write $z=z_i$ for simplicity.
We may write
\begin{equation}\label{eq:(A)_m}
	H_n(\XX_z') - H_n(\XX_{z_-}') = \sum_{m\ge1} \one_{A_{z,m}} (H_n(\XX_z') - H_n(\XX_{z_{-}}')), 
\end{equation}
where $A_{z,m}$ is the event $\{\lceil N_z \rceil = m\}$ with $N_z$ as in Lemma \ref{lem:ff_disag}.
 Then, on the event $A_{ z,m}$, we have that $\X_z'\sm Q_{z,m}=\X_{z_-}'\sm Q_{z,m}$ and hence	
\begin{align}\label{eq:minusWindow}	\one_{A_{z,m}} (H_n(\X_z') - H_n(\X_{z_-}') )\nonumber
		= \one_{A_{z,m}}\big( (H_n(\X_z') - H_n(\X_z'\sm Q_{z,m})) - (H_n(\X_{z_-}')- H_n(\X_{z_-}'\sm Q_{z,m})) \big).
\end{align}
By \eqref{eq:(A)_m} and the conditional Jensen inequality, it is enough to show that 
\begin{equation}\label{eq:inequality1}
	\sup_{n,z} \E \Big[\Big(\sum_{m\geq 1} \one_{A_{z,m}}(H_n( \X_z') - H_n( \X_z'\sm Q_{z,m}))\Big)^4\Big] < \ff.
\end{equation}
We compute
\begin{align*}
&\E \Big[\Big(\sum_{m\geq 1} \one_{A_{z,m}}(H_n( \X_z') - H_n( \X_z'\sm Q_{z,m}))\Big)^4\Big]= 	 \E \Big[\sum_{m\geq 1} \one_{A_{z,m}}\big(H_n( \X_z') - H_n( \X_z'\sm Q_{z,m})\big)^4\Big] \\
&\le \sum_{m\geq 1} \P(A_{z,m})^{1/5} \E\Big[(H_n( \X_z') - H_n( \X_z'\sm Q_{z,m}))^5\Big]^{4/5}\\
	&\le \sum_{m\geq 1} c' \exp(-c_2(m-1)/5) c_{\ms M}m^{p_{\ms M}},
\end{align*}
which is finite. Here, we used that the $A_{z,m}$ are disjoint, the Hölder inequality, the exponential tails of $N_z$ from Lemma \ref{lem:ff_disag}, and \eqref{eq:moment_cond}.

For the case $z_{i-1}\ne z_{i-}$, we note that this implies that $z_i$ is at the boundary of $Q_{a_n}$. Hence, $V_{z_i}\sm V_{z_{i-1}}$ is the union of $Q_{z_i,1}$ and a subset of $Q_{a_n}^c$
Define $B_n$ to be the event that $Q_{n} {\not\lrsa} Q_{a_n-2-4r_0}^c $ in ${\tilde{\Pd}}$.
We now construct $\XX_{z_i}''=T^\ff_{U_{z_i},\eta_{z_i}}(\tilde{\PP}^*_{U_{z_i}})$ and $\XX_{z_{i-1}}''=T^\ff_{U_{z_{i-1}},\eta_{z_{i-1}}}(\tilde{\PP}^*_{U_{z_{i-1}}})$ using radial thinning centered at the origin (rather than at $z_i$). Then,
$$ \De_{i,n} = \E\big[ H_n(\X_{z_i}'') - H_n(\X_{z_{i-1}}'') | \eta_{z_i} \big]. $$
 On $B_n$, $\XX_{z_i}''\cap Q_n=\XX_{z_{i-1}}''\cap Q_n$ by Lemma \ref{lem:rad}. Hence, by the conditional Jensen and Hölder inequalities,
\begin{align*}
		\E\big[\De_{i,n}^4\big] &\le \E\big[1_{B_n^c} ( H_n(\X_{z_i}'') -H_n(\X_{z_{i-1}}'') )^4\big] \\
	& \le \P(B_n^c)^{1/5} \big(\E\big[ (H_n(\X_{z_i}'')- H_n(\XX_{z_{i-1}}''))^5\big]\big)^{4/5}\\
	&\le 2^4 \P(B_n^c)^{1/5} \big(\E\big[ H_n(\X_{z_i}'')^5\big]^{4/5}+ \E\big[H_n(\XX_{z_{i-1}}'')^5\big]\big)^{4/5}.
\end{align*}
Applying \eqref{eq:moment_cond} with $Q_{z,m}=Q_n$, we have 
\begin{align*}
	\E\big[ H_n(\X_{z_i}'')^5\big] \le c_{\ms M}n^{p_{\ms M}}.
\end{align*}
Hence, there are constants $c_1,c_2>0$ such that
\begin{align*}
	& \E\big[\De_{i,n}^4\big] \le C_1 n^{p}\exp(-C_2\sqrt n),
\end{align*}
which is bounded in $n$, thereby concluding the proof.
\enp

The following proposition gives a deterministic criterion for \eqref{eq:moment_cond} to be satisfied.

\bel\label{lem:growth}
Let $\XX$ be an infinite-volume Gibbs point process
with translation-invariant PI satisfying \eqref{as:(A)}.
Consider a translation-invariant functional $H$ and let $H_n$ be defined as in \eqref{eq:H_n_Gibbs}. Then, the condition \eqref{eq:m1} implies \eqref{eq:moment_cond}.
\enl

\bep
Note that
\begin{align}\nonumber
	\E \big[\big(H_n( \X) - H_n( \X\sm Q_{z,m})\big)^5\big]
	&\le c\E \Big[\Big(\sum_{x\in 
		\X \cap Q_{z,m}} \exp( c\X ( B_R(x)))\Big)^5\Big] \\
		&\le c\E \Big[\Big(\sum_{x\in 
		{\Pd} \cap Q_{z,m}} \exp(c{\Pd} ( B_R(x)))\Big)^5\Big].\label{eq:exp_sum}
\end{align}
The first inequality is obtained by adding points in $Q_n\cap \X \cap Q_{z,m}$ to $Q_n\cap \X \sm Q_{z,m}$ one at a time and repeatedly use \eqref{eq:m1}. For the second inequality, we realize $\XX$ as a thinning $\XX = T^\ff(\Pd)$.

Multiplying out $ \big(\sum_{x\in {\PP} \cap Q_{z,m}} \exp(c\Pd( B_R(x)) )\big)^5$, we get that \eqref{eq:exp_sum} is a finite linear combination of terms of the form 
\begin{align*}
	&\E\Big[\sum_{(x_1,\dots ,x_k)\in ({\Pd}_{ Q_{z,m}})_{\neq}^k} \prod_{i=1}^k \exp\big(s_i c {\Pd}\big( B_R(x_i)\big)\big)\Big]\\
	& \le \int_{Q_ m^k}\a_0^k \E_{x_1,\dots,x_k} \Big[\exp\Big(5c{\Pd}\Big( \bigcup_i B_R(x_i)\Big)\Big)\Big]\d x_1\dots \d x_k\\
	&\le C' m^{5d},
\end{align*}
where $k\le 5$, $ s_1+ \dotsm +s_k = 5$ and $C'>0$ is a constant independent of $m$. Here, $(\vp)_{\neq}^k$ denotes the set of $k$-tuples of pairwise distinct points in the point pattern $\vp\in \Nlf$ and $\E_{x_1,\dots,x_k}$ denotes Palm expectation. Inserting this in \eqref{eq:exp_sum} proves the claim.
\enp

The following proposition gives a deterministic criterion for $H$ that ensures that \eqref{eq:as_stab} holds.

\bel\label{lem:s_stab}
Condition \eqref{eq:s1} implies \eqref{eq:as_stab}. 
\enl

\bep
Let $l$ be given. For $n$ sufficiently large, $Q_l\su Q_{w_n,m_n}$ and hence $\vp\cap Q_{w_n,m_n}$ is obtained from $\vp \cap Q_{w_n,m_n}\sm Q_{l}$ by adding the points of $\vp \cap Q_l = \{x_1,\dots,x_k\}$ one add a time. Thus,
\begin{align*}
	&H(\vp \cap Q_{w_n,m_n}) - H(\vp \cap (Q_{w_n,m_n} \sm Q_l)) \\
&= \sum_{i=1}^k ( H((\vp \cap (Q_{w_n,m_n}\sm Q_l)) \cup \{x_1,\dots,x_i\} ) - H((\vp \cap (Q_{w_n,m_n} \sm Q_l))\cup \{x_1,\dots,x_{i-1}\} ). 
\end{align*}
For each term, the limit exists by the assumption \eqref{eq:s1} and translation invariance of $H$. 
\enp

	\subsubsection{Variance positivity}\label{sss:vp}

We now complete the proof of Corollary \ref	{cor:betti} by showing the positivity of the limiting variance.
	\bep[Proof of variance positivity in Corollary \ref{cor:betti}]
		 In the proofs of Theorem \ref{thm:martingale} and \ref{thm:gibbs_CLT}, instead of using unit cubes, we could have used cubes of the form $az+Q_a$ for $z\in \Z$ and any $a>0$. In the following we choose $a>2s+2(s\vee r_0)$. Defining $U_{z,a}$, $\eta_{z,a}$ and $\XX_{z,a}'$ analogously to $U_z$, $\eta_z$ and $\XX_z'$, the limiting variance is $\s^2=a^{-d}\E\De_{0,a}^2$, where $\Delta_{0,a}$ is the $L^2$-limit of 
		\begin{align*}
			\E [H_n( \X_{0,a}') 
			- H_n( \X_{0_-,a}') \mid \eta_{0,a}]. 
		\end{align*}
		Thus, it is enough to show $\liminf_n \E[\E [H_n( \X_{0,a}') - H_n( \X_{0_-,a}') \mid \eta_{0,a} ]^2]>0$. 
		
		For the proof, we construct an auxiliary process ${\X}_{0,a}'' = \eta_{0_-,a}\cup T^{\ff}_{U_{0,a},\eta_{0_-,a}}(\tilde{\PP}_{U_{0,a}})$. Note that $\XX_{0,a}''$ depends only on $\XX$ via $\eta_{0_-,a}$ and $\XX_{0,a}''\cap Q_a = \es$.
		For a point process $\mc Y$, we define disjoint events
		\begin{align*}
			E_1{}& = \{ \mc Y \cap Q_a = \es\} \\
			E_2{}& = \{ \mc Y \cap (Q_a\sm Q_{a-2(r_0\vee s)}) = \es, \beta_q^{r,s}(\mc Y \cap Q_{a-2(r_0\vee s)})=1 \}.
		\end{align*}
		Both events have positive probabilities when $\mc Y$ is a Poisson process, see \cite[Example 1.8]{shirai}. Then,
		\begin{align*}
			& \E\Big[\E \big[H_n( \X_{0,a}') - H_n( \X_{0_-,a}')  \mid \eta_{0,a} \big]^2 \Big]\\
			&\quad \ge \E\Big[\E \big[(H_n( \X_{0,a}') - H_n( \X_{0_-,a}') )(\one\{{\X}_{0,a}'\in E_1\}+\one\{{\X}_{0,a}'\in E_2\}) \mid \eta_{0,a} \big]^2 \Big] \\
			&\quad = \E\Big[\E \big[(H_n({\X}_{0,a}'\sm Q_a) - H_n( \X_{0_-,a}') )\one\{{\X}_{0,a}'\in E_1\} \\
			&\quad \quad+ (H_n( {\X}_{0,a}'\sm Q_a) - H_n( \X_{0_-,a}') + 1 )\one\{ \X_{0,a}'\in E_2 \} \mid \eta_{0,a} \big]^2 \Big]\\
			& \quad = \E\Big[\E \big[H_n( {\X}_{0,a}'') - H_n( \X_{0_-,a}') \mid\eta_{0,a} \big]^2 \one\{{\X}_{0,a}'\in E_1\} \\
			&\quad \quad + \E\big[H_n( {\X}_{0,a}'') - H_n( \X_{0,a}') + 1 \mid \eta_{0,a} \big]^2 \one\{\X_{0,a}'\in E_2\} \Big]. 
		\end{align*}
		
		Here, we used disjointness of $E_1$ and $E_2$, measurability of $ \one\{{\X}_{0,a}'\in E_1\} $ and $ \one\{{\X}_{0,a}'\in E_2\} $ with respect to $\eta_{0,a}$, and the fact that on $E_1$ and $ E_2$, $\X_{0,a}'$ agrees with ${\X}_{0,a}''$ on $\R^d$ and $\R^d\sm Q_a$, respectively. Using the conditional Jensen inequality and the fact that both ${\X}_{0,a}''$ and $ \X_{0_-,a}'$ depend on $\eta_{0,a}$ only via $\eta_{0_-,a}$, we obtain the lower bound
				\begin{align*}
			& \E\Big[\E \big[H_n( {\X}_{0,a}'') - H_n( \X_{0_-,a}') \mid \eta_{0_-,a} \big]^2 \E\big[\one\{{\X}_{0,a}'\in E_1\} \mid\eta_{0_-,a} \big] \\
			&\quad\quad + \E\big[H_n( {\X}_{0,a}'') - H_n( \X_{0_-,a}') + 1 \mid \eta_{0_-,a} \big]^2\E\big[\one\{{\X}_{0,a}'\in E_2\} \mid\eta_{0_-,a} \big] \Big] \\
			&\quad\ge \E\Big[\max\big\{\E \big[H_n( {\X}_{0,a}'') - H_n( \X_{0_-,a}') \mid \eta_{0_-,a} \big]^2,\, \big(\E \big[H_n( {\X}_{0,a}'') - H_n( \X_{0_-,a}') \mid\eta_{0_-,a} \big]+1\big)^2 \big\} \\
			&\quad\quad \times \min_{i=1,2} \P({\X}_{0,a}'\in E_i\mid\eta_{0_-,a}) \Big]\\
			&\quad\ge \tfrac 1 4\E\Big[ \min_{i=1,2} \P({\X}_{0,a}'\in E_i\mid\eta_{0_-,a}) \Big],
		\end{align*}	
		 where the inequality $a_1b_1 + a_2b_2\geq (a_1\vee a_2) \cdot (b_1\wedge b_2)$ for $a_1,a_2,b_1,b_2\geq 0$ was used to obtain the second inequality, and the third inequality used that $x^2 \vee (x+1)^2\ge 1/4$ for all $x\in \R$.
		 
For $\P( {\X}_{0,a}'\in E_1\mid\eta_{0_-,a} )$, we have the deterministic lower bound 
$$\P( {\X}_{0,a}'\in E_1\mid \eta_{0_-,a} ) \ge \P( \tilde{\PP}\in E_1\mid \eta_{0_-,a} ) = \P(\tilde{\Pd}\cap Q_a = \es).$$ 
It thus remains to bound $\P( {\X}_{0,a}'\in E_2) $. We have 
		\begin{align*}
			&\P( {\X}_{0,a}'\in E_2)\\
			&\quad = \P(\X_{0,a}' \cap (Q_a\sm Q_{a-2(r_0\vee s)}) = \es )
			 \P(\beta_q^{r,s}(\X_{0,a}\cap Q_{a-2(r_0\vee s)})=1 \mid \X_{0,a} \cap (Q_a\sm Q_{a-2(r_0\vee s)}) = \es)\\
			 &\quad \ge \P(\tilde{\PP} \cap (Q_a\sm Q_{a-2(r_0\vee s)})=\es)\P(\beta_q^{r,s}(\X( Q_{a-2(r_0\vee s)}, \es))=1 ),
		\end{align*}
		where the last inequality used the DLR equations \eqref{eq:DLR1}. The probability $\P(\beta_q^{r,s}(\X( Q_{a-2(r_0\vee s)}, \es))=1)$ is non-zero because the corresponding probability for a Poisson process $\P(\beta_q^{r,s}({\PP} \cap Q_{a-2(r_0\vee s)})=1 )$ is positive and $\X( Q_{a-2(r_0\vee s)}, \es)$ has a positive density with respect to $\Pd\cap Q_{a-2(r_0\vee s)}$ by the assumption $\kappa>0$. 
	\end{proof}

	\section{Proof of Theorem \ref{thm:2} -- normal approximation for Gibbsian score sums}
\label{sec:gibbs}
Before elaborating on the technical details of the proof of Theorem \ref{thm:2}, we first give a broad overview of the main idea. 
{We consider the general framework for a random measure $\Xi$ on \cb{$\R^d$}. In the proof of Theorem~\ref{thm:2}, this will be applied to the random measure }
	\begin{align}
	\Xi := \Xi[\XX]& :=\sum_{x \in \XX \cap Q} g(x, \XX)\,\delta_x, \label{eqn:poixidefXi}
\end{align}
noting that $|\Xi|:=\Xi(\R^d)=H(\XX)$ is the score sum defined at \eqref{eqn:poixidef}. 

The concept of Palm theory is of central importance for our approach.
We recall that a collection of random measures $\{\Xi_x\}_{x \in \R^d}$ is a \emph{Palm version} of $\Xi$ if it satisfies 
\begin{align}
	\E\Big[ \int_{\R^d} f(x,\Xi) \,\Xi(\d x)\Big]=\int_{\R^d} \E [f(x,\Xi_x)] \,\La(\d x)  \label{def:palm1}
\end{align}
for all non-negative measurable $f$, where $\La := \E[\Xi]$ denotes the intensity measure. Moreover, we set $\s^2:=\Var (|\Xi|) := \Var(\Xi(\R^d))$. For an introduction to Palm theory and in particular for the definition of Palm processes with respect to {random measures}, we refer to \cite[Section 6]{randMeas}.

Assume that $\Xi$ and its Palm version $\{\Xi_x\}_{x \in \R^d}$ are defined on the same probability space and let
\begin{align}
	Y_x :=|\Xi_x|-|\Xi|, \quad  \De_x :=\f{Y_x}{{\s}}=\f{|\Xi_x|-|\Xi|}{{\s}},\quad x \in \R^d. \label{def:Ws}
\end{align}
{The proof of Theorem \ref{thm:2} relies on the following theorem from \cite{chen}, which follows by an application of Stein's method.}
%
%
\bet [Theorem 3.1 in \cite{chen}] \label{thm:rollin}
Let $\Xi$, $\{\Xi_x\}_{x \in \R^d}$, $	\{Y_x\}_{x \in \R^d} $, and $\{\De_x\}_{x \in \R^d}$ be as in \eqref{def:Ws}. Then,
\begin{align*}
	\dk\Big(\f{|\Xi|-\E |\Xi|}{{\s}}, N(0,1)\Big) \le 2E_1+5.5E_2+5E_3+10E_4+7E_5
\end{align*}
where
\begin{align*}
	E_1&:=\f1{\s^2} \E\Big|\int_{\R^d} (Y_x\, \one \{|Y_x|\le \s\}-\E\big[Y_x \, \one \{|Y_x|\le \s\}\big] )\,\La(\d x)\Big|,\\
	E_2&:=\f1{\s^3}\int_{\R^d} \E \big[Y_x^2\,\one \{|Y_x| \le \s\}\big] \,\La(\d x),\\
	E_3&:=\f1{\s^2}\int_{\R^d} \E \big[|Y_x|\,\one\{|Y_x| > \s\}\big] \,\La(\d x),\\
	E_4&:=\f1{\s^2} \int_{-1}^1 \int_{\R^d} \int_{\R^d} \Cov\big(\phi_{x}(t), \phi_{y}(t)\big) \, \La(\d x)\,\La(\d y)\,\d t,\\
	E_5&:=\f1{\s}\Big(\int_{-1}^1 \int_{\R^d} \int_{\R^d} |t|\,\Cov\big(\phi_{x}(t), \phi_{y}(t)\big) \, \La(\d x)\,\La(\d y)\,\d t\Big)^{1/2},
\end{align*}
where to simplify the notation, we write
\begin{align*}
	\phi_{x}(t)=\begin{cases}
		\one\{1 \ge \De_x>t>0\},\quad &t>0,\\
		\one\{-1 \le \De_x<t<0\},\quad &t<0.
	\end{cases}
\end{align*}
\ent
Typically, the most delicate expressions are the terms $E_1$, $E_4$, and $E_5$. This is because, loosely speaking, those terms encode a bound on the deviation of the difference $Y_x$ when averaged over space. In contrast, the terms $E_2, E_3$ only involve moment bounds at individual space points.

%
%
To prove of Theorem \ref{thm:2}, we will specify the processes $\Xi$ and $\{\Xi_x\}_{x\in \R^d}$ in Theorem \ref{thm:rollin} and show that the terms $E_1,\dots,E_5$ can be bounded by the right-hand side asserted in Theorem \ref{thm:2}. We write $\{\XX_x^\Xi\}_{x\in Q}$ for a Palm version of $\XX$ with respect to $\Xi$ in $Q$. That is, for all measurable $f:\X \ti \mathbf N \to [0,\ff)$,
\begin{align}
\E \int_{Q} f(x,\XX) \,\Xi(\d x)= \int_{Q} \E f(x,\XX_x^\Xi) \, \La(\d x). \label{def:palm2}
\end{align}
We note that  assumption \eqref{eq:m2} ensures the $\s$-finiteness of $\La$, and thereby the existence of a Palm version $\{\XX_x^\Xi\}_{x\in Q}$, see \cite{kallenberg}. 
{Taking $f(x,\XX) = \tilde{f}(x,\Xi[\XX])$ in \eqref{def:palm2} shows that $\{\Xi[\XX_x^\Xi ] \}_{x\in Q}$ is a Palm version $\{\Xi_x\}_{x\in Q}$ of $\Xi$ satisfying \eqref{def:palm1}.}
Lemma \ref{lem:palmgibbs} below shows that the reduced Palm process of $\XX$ with respect to $\Xi$ at some point $x \in Q$ is a Gibbs process. This property is used in the proof of Theorem \ref{thm:2} to justify that we can construct a coupling of $\Xi$ and its Palm version via disagreement coupling.

%
%
\bel
\label{lem:palmgibbs} For a.a.\  $x \in Q$, the reduced process $\XX_x^{!,\Xi}:=\XX_x^\Xi\sm \{x\}$ is a Gibbs process with PI $\k_x$ given by
	$$\k_x(y,\om):=\k(y,\om \cup \{x\}) \f{g(x,\om \cup \{x,y\})}{g(x,\om \cup \{x\})}$$
	where $0/0:=0$.
	\enl

Note that for $\|y-x\| > R(x,\omega)$, we have that $\k_x(y,\om)=\k(y,\om)$. This will be the case in our application of Lemma \ref{lem:palmgibbs} in the proof of Theorem \ref{thm:2}.

	%
	%
	\bep
	We must verify \eqref{eGNZ} for $\XX_x^{!,\Xi}$. Let $f:\R^d \ti \mathbf N \to [0,\ff)$ be measurable. We note that the Palm versions are only defined uniquely up to a \cb{Lebesgue measure zero set} of points $x \in Q$. Therefore, it suffices to show that for all measurable $h:\R^d \to [0,\ff)$,
	\begin{align}
		\int_{Q}h(x) \E\Big[ \sum_{Y_i \in \XX_x^{!,\Xi}}  f(Y_i, \XX_x^{!,\Xi})\Big] \,\La(\d x)=\int_{Q}h(x) \int\E  \Big[f(y,\XX_x^{!,\Xi} \cup \{y\}) \k_x(y,\XX_x^{!,\Xi}) \Big]\d y \,\La(\d x). \label{eqn:palmlem}
	\end{align}
Here, the right-hand side is, by the definitions  of $\XX_x^{!,\Xi}$ in \eqref{def:palm2} and of $\k_x(y,\om)$, given by
\begin{align*}
	&\E \Big[\int_{Q}h(x) \int_{Q}f(y,(\XX \sm \{x\})\cup \{y\}) \k_x(y,\XX\sm \{x\}) \,\d y \,\Xi(\d x)\Big] \\
	&\quad =\E\Big[ \sum_{X_i \in \XX} h(X_i) \int_{Q}f(y,(\XX \sm \{X_i\})\cup \{y\}) \k(y,\XX) g(X_i,\XX \cup \{y\}) \,\d y\Big].
\end{align*}
We use Fubini's theorem and apply \eqref{eGNZ} to the integral over $y$ to get
\begin{align*}
	\E \Big[\sum_{X_i \ne X_j \in \XX} h(X_i) f(X_j,\XX \sm \{X_i\}) g(X_i,\XX) \Big] 	=\E\Big[	\int_{Q}h(x)  \sum_{X_j\in \XX \sm \{x\}} f(y,\XX \sm \{x\})  \,\Xi(\d x)\Big].
\end{align*}
It follows from the definition \eqref{def:palm2} of $\XX_x^{!,\Xi}$ that the above coincides with the left-hand side in \eqref{eqn:palmlem}.
\enp

The proof of Theorem \ref{thm:2} is based on the following coupling construction. Let 
\begin{align*}
	 r:=4s := 4\max\{r_0,c_{\ms s} \log |Q|\}\quad \text{ with }\quad c_{\ms s} := {120}\max(c_{\ms{SPT},2}^{-1}, c_{\ms{es}}^{-1}),
\end{align*}
 where $c_{\ms{SPT},2}>0$ is the constant from Corollary \ref{cor:dec} and $c_{\ms{es}}>0$ is introduced in \eqref{eq:s2}. Here, the choice of the prefactor 120 guarantees that $r$ and $s$ are large enough such that the bounds asserted in Lemma \ref{lem:Yintbou} and Lemma \ref{lem:thm2cov} hold. Let $\tilde \XX$ be a Gibbs process with PI $\k$ and let $\tilde \XX_x^\Xi$, $x \in Q$, be independent reduced Palm processes of $\XX$ with respect to $\Xi$ at $x$, which is possible by the existence of uncountable product measures, see  \cite[Corollary 6.18]{kallenberg}. We specify $Q^+:=B_{2r}(Q)$, $\Psi:=\tilde \XX \cap (Q^+)^c$, $\Psi_x=\tilde \XX_x^\Xi \cap ((Q^+)^c \cup B_{s}(x)) $ and ${\k'_x}(z,\vp):=\k_x(z,\vp_{B_{s}(x)^c} \cup (\Psi_x \cap B_{s}(x))) \one\{z \in B_{s}(x)^c\}$, with $\k_x$  from Lemma \ref{lem:palmgibbs}. Hence, when conditioned on $\Psi$ and $\{\Psi_{x}\}_{x\in Q}$, we can carry out the disagreement coupling of $\tilde \XX$ and $\tilde \XX_{x}^\Xi$, i.e., 
\begin{align}
	\XX&:=T^{\mathsf{rad}}_{Q^+,(Q^+)^c,\Psi}(\Pds_{Q^+}) \cup \Psi,\qquad	\XX_{x}^\Xi:= \begin{cases}
		T^{\mathsf{rad},x}_{Q^+,(Q^+)^c,\Psi_x \cap (Q^+)^c}(\Pds_{Q^+})\cup \Psi_x & \text{ if $R(x, \tilde \XX_x^\Xi)\le s$,}\\
		\tilde \XX_x^\Xi\cup \Psi_x & \text{ if $R(x,\tilde \XX_x^\Xi)> s$,}
	\end{cases} 
	\label{eqn:coupling}
\end{align} 
{where $T^{\mathsf{rad},x}$ denotes radial thinning with respect to $\k_x'$.} It is important to note that despite the more complicated definition, the process $\XX_x^\Xi$ has also the correct distribution as a Palm process. This is because of the DLR equation, which is applied by first conditioning on $\Psi_x$, and then sampling the remaining Gibbs process with one of the two options depending on the value $R(x,\tilde \XX_x^\Xi)$.
Also note that while in the random measure $\Xi$ we only sum over points in $Q$, the dependence through the score functions means that $\Xi$ also depends on the configuration of $\XX$ outside $Q$. Therefore, we use the enlarged window $Q^+ = B_{2r}(Q)$ in the disagreement couplings in \eqref{eqn:coupling}.

Given $\Psi_x$, it follows from Proposition \ref{th:dc_prop} that ${T^{\mathsf{rad},x}}_{Q^+,(Q^+)^c,\Psi_x \cap (Q^+)^c}({\Pds_{Q^+}})$ is a Gibbs process on $Q^+$ with PI ${\k'_x}$ and boundary condition $\Psi_x$. On the other hand, $ \XX_x^\Xi \cap Q^+$ conditioned on $\Psi_x$ is by \eqref{eq:DLR2} also a Gibbs process on $Q^+$ with PI $\k'_x$ and boundary condition $\Psi_x$. Since the distribution of a Gibbs process is unique on compact domains, we conclude that {conditionally on $\Psi_x$,} $${T^{\mathsf{rad},x}}_{Q^+,(Q^+)^c,\Psi_x \cap (Q^+)^c}({\Pds_{Q^+}}) \stackrel{d}{=} \tilde \XX_x^\Xi \cap Q^+ x,$$
{and similarly $T^{\mathsf{rad}}_{Q^+,(Q^+)^c,\Psi}(\Pds_{Q^+})  \stackrel{d}{=} \tilde \XX \cap Q^+ $ conditionally on $\Psi$.}
Therefore, $ \XX_{x}^\Xi \stackrel d= \tilde \XX_x^\Xi$ and similarly $\XX \stackrel d=\tilde \XX$ {by the DLR-equations \eqref{eq:DLR2}.} In particular,  $ \{\XX_{x}^\Xi\}_{x \in Q}$ is a Palm version of $ \XX$ with respect to $\Xi$ on $Q$. It follows that $\{\Xi[\XX_{x}^\Xi]\}_{x \in Q}$ is a Palm version of $\Xi$ on $Q$. 

The moment bound provided by the following lemma will be used in the proof of Theorem \ref{thm:2}.

%
%
\bel[Moment bound for $Y_{x}$]
\label{lem:Yintbou}
Let $Q \subset \R^d$ be a bounded Borel set. Assume that the PI $\k$ satisfies \eqref{as:(A)} and that \eqref{eq:m2} and \eqref{eq:s2} hold. For $x \in Q$ let $\Xi:=\Xi[\XX]$ and $\Xi_{x}:=\Xi[\XX_{x}^\Xi]$, where $\XX$ and $\XX_{x}^\Xi$ are given at \eqref{eqn:coupling}. For  $m \le 4$ there is some $q_0=q_0(d, \a_0, r_0,c_{\ms{es}}, c_{\ms{SPT},2})>0$ such that for  $|Q|>q_0$ we have
$$
\int_{Q} \E[|Y_{x}|^m]\, \La (\d x) \le c_{\ms{Y}} |Q| s^{dm},
$$
where $c_{\ms{Y}}=c_{\ms{Y}}({d},\a_0, r_0, m,c_{\ms{m}})>0$.
\enl 

\bel[Covariance bound] \label{lem:thm2cov}
Let $Q \subset \R^d$ be a bounded Borel set. Assume that the PI $\k$ satisfies \eqref{as:(A)} and that \eqref{eq:m2} and \eqref{eq:s2} hold. Let $p > 0$  and $F:\mathbf N \to [0,\ff)$ be measurable. Furthermore, assume that there is some $q_0=q_0(d, \a_0, r_0, c_{\ms{es}}, c_{\ms{SPT},2}, p)>0$ such that $F(\vp)\le  |Q|^p \#\vp$ for $|Q|>q_0$ and $\vp \in \Nlf$. For those $Q$ we have
$$
\iint_{Q^2 \cap \{|x-y|>4r\}} \Cov(F_{x}, F_{y})\, \La^2(\d x, \d y) \le |Q|^{-1}, 
$$
where $F_{x} := F(\XX_{x, r})$ with $\XX_{x, r}:=\XX_{x} \cap B_{r}(x)$  for any measurable $F:\mathbf N \to [0,\ff)$. 
\enl

\ber \label{rem:cov}
The bound in Lemma \ref{lem:thm2cov} is not at all sharp but sufficient for our purposes.
As the proof reveals, Lemma \ref{lem:thm2cov} also holds if $F_{x}$ or $F_{y}$ is replaced by $F(\XX \cap B_{r}(x))$ or $F(\XX \cap B_{r}(y))$, respectively. In this case, the first step of the proof (where $F_{z}$ is truncated) can be avoided.
\enr

From Lemma \ref{lem:Yintbou} we can now directly conclude that
\begin{align*}
	\s^3(E_2+E_3) \le   \int_{Q} \E[Y_{x}^2]\,\La(\d x)\le c_{\ms{Y}} |Q| s^{2d} .
\end{align*}
For the remaining error terms, we will apply Lemma \ref{lem:Yintbou} and Lemma \ref{lem:thm2cov} to deduce that
\begin{align}
E_1\le c_{E_1} s^{2d} \s^{-3} |Q|,\qquad E_4 \le c_{E_4} s^{2d} \s^{-3} |Q|,\qquad E_5 \le c_{E_5} s^{2d} \s^{-3} |Q| \label{E145}
\end{align}
for some constants $c_{E_i}=c_{E_i}(d,\a_0,c_{\ms{m}})>0,\,i=1,4,5$. This implies the assertion of Theorem \ref{thm:2} with the constant $c_{\ms{norm}}:=c_{\ms{s}}^{2d} (c_{\ms{Y}}+c_{E_1}+c_{E_4}+c_{E_5})$.

We now proceed with the proofs of the two lemmas and thereafter establish the bounds in \eqref{E145}.

\begin{proof}[Proof of Lemma \ref{lem:Yintbou}]
Set  $\Xi^>:=\Xi- \Xi^\le$, $\Xi_{x}^>:=\Xi_{x} - \Xi_{x}^\le$ and $\XX \De\XX_{x}^\Xi:=(\XX \setminus \XX_{x}^\Xi) \cup (\XX_{x}^\Xi \setminus \XX)$, where 
	\begin{align}
		\Xi_{x}^\le:= \sum_{y \in \XX_{x}^\Xi \cap Q} g(y,\XX_{x}^\Xi) \mathds 1\{R(y,\XX_{x}^\Xi)\le s\}\de_y,\qquad \Xi^\le:= \sum_{y \in \XX \cap Q} g(y,\XX) \mathds 1\{R(y,\XX)\le s\}\de_y.\label{def:Xile}
	\end{align}
	Furthermore, we also put $\Xi_{x, r} :=\Xi_{x}^\le \big(B_{r}(x)\big)$.
We decompose $Y_{x}$ as
	\begin{align*}
	\big(|\Xi_{x}^> |-|\Xi^>|\big) + \big(\Xi_{x}^\le(B_r(x))-\Xi^\le (B_r(x)) \big) +\big(\Xi_{x}^\le ( B_{r}(x)^c)-\Xi^\le (  B_{r}(x)^c)\big) \one \{\XX \Delta \XX_{x}^\Xi \not \su B_{3s}(x)\}.
	\end{align*}
From the Hölder inequality applied to both the expectation and the integral with respect to $\La$, we get
	\begin{align}
		\int_{Q} \E[ | Y_{x}|^m] \, \La(\d x) \le \hspace{-0.4cm} \sum_{\substack{(i_1,\dots,i_6):\\ i_1+\cdots+i_6=m}} &\Bigg[\Big(\int_{Q} \E [\Xi_{x}^\le\big(B_{r}(x)\big)^m] \, \La(\d x)\Big)^{\f{i_1}{m}}  \Big(\int_{Q} \E [\Xi^\le\big(B_{r}(x)\big)^m] \, \La(\d x)\Big)^{\f{i_2}{m}}\nonumber\\
		&\times \Big(\int_{Q} \E [|\Xi_{x}^>|^m] \, \La(\d x)\Big)^{\f{i_3}{m}}  \Big(\int_{Q} \E [|\Xi^>|^m] \, \La(\d x)\Big)^{\f{i_4}{m}}\nonumber\\
		&\times \Big(\int_{Q} \E [|\Xi_{x}|^m  \one \{\XX \De\XX_{x}^\Xi \not\su B_{3s}(x)\}] \, \La(\d x)\Big)^{\f{i_5}{m}}\nonumber\\
		&\times \Big(\int_{Q} \E [|\Xi|^m  \one \{\XX \De\XX_{x}^\Xi \not\su B_{3s}(x)\}] \, \La(\d x)\Big)^{\f{i_6}{m}} \Bigg].\label{eqn:minpalmloc}
	\end{align}
We now bound the six integrals on the right-hand side in \eqref{eqn:minpalmloc} separately. For the first integral, we obtain by an expansion of the $m$th power and iteratively applying \eqref{eGNZ},
\begin{align}\nonumber
		&\int_{Q} \E \big[\Xi_{x}^\le \big(B_r(x)\big)^m\big] \La(\d x)\\ \nonumber
		&\le	\int_{Q}\E \Big[ \Big(\sum_{X_j \in \XX_{x}^\Xi \cap B_{r}(x)}   g(X_j, \XX_{x}^\Xi)\Big)^m\Big]\, \La(\d x)\\ \nonumber
		&= \E \Big[\sum_{X_i \in \XX\cap Q}g(X_i, \XX) \Big(\sum_{X_j \in  \XX \cap B_{r}(X_i)}   g(X_j, \XX)\Big)^m \Big]\\ \nonumber
		&=\int_{Q}  \sum_{k=0}^{m} \sum_{\substack{(i_1, \dots, i_{k + 1}):\\ i_{1}+\cdots+i_{k + 1}=m}} \E\Big[\k(x,\XX) \, g(x,\XX\cup \{x\})^{i_{k + 1}+1}\sum_{(X_1, \dots, X_k) \in (\XX  \cap B_{r}(x))_{\neq}^k}  \prod_{j=1}^{k} g(X_j,\XX \cup \{x\})^{i_j}  \Big]\,\d x\\ 
		&=\int_{Q}   \sum_{k=0}^{m} \sum_{\substack{(i_1, \dots, i_{k + 1}):\\ i_{1}+\cdots+i_{k + 1}=m}} \int_{ B_{r}(x)^k} \E\Big[\k(\{x,\ww\}, \XX)\, g(x,\XX \cup \{x,\ww\})^{i_{k + 1}+1}\prod_{j=1}^{k} g(w_j,\XX \cup \{x,\ww\})^{i_j}\Big] \,\d\ww\, \d x, \label{eqn:firstfactor}
\end{align}
	where $\k(\{x,\ww\},\XX):=\k(x,\XX)\k(w_1,\XX \cup \{x\})\cdots \k(w_k,\XX \cup \{x,w_1,\dots,w_{k-1}\})$ for $\ww = (w_1,\ldots,w_k)\in (\R^d)^k$.
	Using the moment condition \eqref{eq:m2}, we find that the above is bounded by
	\begin{align*}
		&|Q| \sum_{k\le m} \sum_{\substack{(i_1, \dots, i_{k + 1}):\\ i_{1}+\cdots+i_{k + 1}=m}} \a_0^{k+1} \sup_{\xx \in Q^{k+1}} {\E [}g(x_1,\XX \cup \{\xx\})^{k+1}]|B_{r}(x)|^k \le c_1  |Q| r^{dm}
	\end{align*}
	for some constant $c_1=c_1(d,\a_0,m, c_{\ms{m}})$. The same bound can be established for the second integral on the right-hand side in \eqref{eqn:minpalmloc}.
	
	For the third integral in \eqref{eqn:minpalmloc}, we find analogously to above that 
	\begin{align}
		&\int_{Q} \E\Big[\Big|\sum_{X_j \in \XX_{x}^\Xi} g(X_j,\XX_{x}^\Xi) \mathds 1\{R(X_j,\XX_{x}^\Xi)>s\}\Big|^m\Big] \La(\d x)\nonumber\\
		&\le \int_{Q}   \sum_{k=0}^{ m} \sum_{\substack{(i_1, \dots, i_{k + 1}):\\ i_{1}+\cdots+i_{k + 1}=m}} \Big\{\int_{ Q^k} \E\Big[\k(\{x,\ww\}, \XX)\, g(x,\XX \cup \{x,\ww\})^{i_{k + 1}+1}\nonumber\\
		&\qquad \qquad \qquad \qquad \qquad  \times \prod_{j=1}^{k} g(w_j,\XX \cup \{x,\ww\})^{i_j} \one\{R(w_j,\XX \cup \{x,\ww\})>s\} \Big] \,\d\ww\, \d x\Big\}\nonumber\\
		&\le   \sum_{k=0}^{m} \sum_{\substack{(i_1, \dots, i_{k + 1}):\\ i_{1}+\cdots+i_{k + 1}=m}} \Big\{ \a_0^{k+1} |Q|^{k+1} \nonumber\\
		&\qquad \qquad \qquad \qquad \qquad  \times\sup_{\xx \in Q^{k+1}}\E\big[  g(x_{k+1},\XX \cup \{\xx\})^{i_{k + 1}+1} \prod_{j=1}^{k}g(x_j,\XX \cup \{\xx\})^{i_j} \one\{R(x_j,\XX \cup \xx ) > s\}\big]\Big\}.\nonumber
		\end{align}
	Here, we bound every indicator in the expectation except for the one with $j=1$ by $1$. Then, we apply the Hölder inequality with $m+2$ factors (counted with multiplicities). This gives by \eqref{eq:m2} and \eqref{eq:s2} the bound
	\begin{align}
	& \sum_{k\le m} \sum_{\substack{(i_1, \dots, i_{k + 1}):\\ i_{1}+\cdots+i_{k + 1}=m}} \a_0^{k+1} |Q|^{k+1} \sup_{\xx \in Q^{k+1}}\P(R(x_1,\XX \cup \{\xx\}) > s )^{{\f1{m+2}}}\sup_{\xx \in Q^{k+1}}\E[g(x_1,\XX \cup \{\xx\})^{m+2}]^{\f{m+1}{m+2}}\nonumber\\
	&\quad \le |Q|^{m+1} e^{-\f{c_{\ms{es}}s}{2(m+2)}},
	\label{eqn:Yint3}
	\end{align}
where the last inequality holds for $|Q|$ large enough. The integral of  $\E [|\Xi^>|^m]$ in \eqref{eqn:minpalmloc} can be bounded similarly.
	
For the fifth integral in \eqref{eqn:minpalmloc}, we find from the Hölder inequality the bound
	\begin{align}
		&\Big(\int_{Q}\E \Big[{|\Xi_{x}|^{m+1}}\Big]\,\La(\d x)\Big)^{\f{m}{m+1}} \Big( \int_{Q} \P\big(\XX \De\XX_{x}^\Xi \not\su B_{3s}(x)  \big) \, \La(\d x)\Big)^{\f1{m+1}}. \label{eqn:Yint5}
	\end{align}
	Now, it follows from the bounded moments condition \eqref{eq:m2} analogously to the bounds of the terms in \eqref{eqn:firstfactor} that the first integral is bounded by $c_2 |Q|^{m+2}$ for some $c_2=c_2(\a_0,m, c_{\ms{m}})>0$, whereas we bound the probability in the second integral by
	\begin{align}
		\P(R(x,\XX_{x}^{\Xi})>s)+ 	\P\big(\XX \De\XX_{x}^\Xi \not\su B_{3s}(x), \,  R(x,\XX_{x}^{\Xi})\le s\big). \label{eq:5int2}
	\end{align}
By definition of $\XX_{x}^{\Xi}$, the $\La$-integral of the first probability is given by
\begin{align}
\E \Big[ \sum_{x \in Q} \one\{R(x,\XX)>s\} g(x,\XX)\Big]&=\int_{Q} \E[\k(x,\XX) \one\{R(x,\XX \cup \{x\})>s\} g(x,\XX)]\,\d x\nonumber\\
&\le \a_0 \int_{Q}  \P(R(x,\XX \cup \{x\})>s)^{1/2} \E [g(x,\XX \cup \{x\})^2]^{1/2}\,\d x\nonumber\\
&\le e^{-c_{\ms{es}}s /4},\label{eq:Rpalm}
\end{align}
where the last inequality holds for $|Q|$ large enough and we have used \eqref{as:(A)} and the Cauchy-Schwarz inequality in the last line. 
To treat the second probability in \eqref{eq:5int2}, we apply Theorem \ref{th:disdis}. 
To that end, we note that for $R(x,\XX_{x}^{\Xi})\le s$ we have by the stopping property of $R$ that conditioned on $\Psi_{x,0}:=\Psi_x \cap B_s(x)$,
$$
{\k'_x}(y,\vp)=\k_x(y,\vp)=\k(y,\vp),\quad y \in Q \setminus B_{2s}(x).
$$
Hence, $\k$ and ${\k'_x}$ satisfy assumption \eqref{as:(A)} and we obtain from Theorem \ref{th:disdis} with $B:=B_s(x)$ and $A=P:=Q \setminus B_{3s}(x)$ that 
	\begin{align}
	\P\big(\XX \De\XX_{x}^\Xi \not\su B_{3s}(x), \,  R(x,\XX_{x}^{\Xi})\le s \big) \le c_{\ms{DP}}|Q| \exp(-c_{\ms{SPT},2} s/3), \label{eqn:Yint5b}
	\end{align}
	where $c_{\ms{DP}}$ and $c_{\ms{SPT},2}$ are the constants from Theorem \ref{th:disdis}. 
	The sixth integral in \eqref{eqn:minpalmloc} is treated analogously. 
\end{proof}

%
%

%
%
\bep[Proof of Lemma \ref{lem:thm2cov}]
For $z \in \{x,y\}$ we set 	$R_{z} := R(z,\XX_{z})$ and write 
	$$\tilde F_{z}:= F(\XX_{z} \cap B_{r}(z)) \one\{R_{z}\le s\}.$$
Then,
	\begin{align}
	\Cov(F_{x}, F_{y})
	&=\Cov(\tilde F_{x}, \tilde F_{y}) + \Cov(\tilde F_{x},  F_{y}\one\{R_{y}> s\}) +  \Cov(F_{x} \one\{R_{x}> s\}, F_{y} )\label{lemco1bou2}
\end{align}
As the steps for the second and the third covariance are similar, we only discuss the second in detail.  Here, using that $F$ is nonnegative,  the Cauchy-Schwarz inequality gives the upper bound
\begin{align}
	\E[F_{x}^2]^{1/2} \E[F_{y}^2 \one\{R_{y}>s\}]^{1/2}. \label{EFxy}
\end{align}
An upper bound for the negative covariance can be derived similarly. By Jensen's inequality and \eqref{as:(A)}, the $\La^2$-integral of \eqref{EFxy} is bounded by
\begin{align*}
	&\La(Q)\Big(\iint \E[F_{x}^2]  \E[F_{y}^2 \one\{R_{y}>s\}]\La^2(\d x, \d y)\Big)^{\f 12}\\
	&\le   \La(Q) \Big(\E\Big[\sum_{x \in \XX} F(\XX \cap B_r(x))^2 g(x,\XX) \Big] \E\Big[\sum_{y \in \XX} F(\XX \cap B_r(y))^2 g(y,\XX) \one\{R(y,\XX)>s\} \Big]\Big)^{\f 12}.
\end{align*}	
We set $N_{z,r}:=\#((\XX \cup \{z\}) \cap B_{r}(z))$ for $z \in \{x,y\}$, $q_{y}:=\P(R(y,\XX \cup \{y\})>s)$ and find from the GNZ equation {\eqref{eGNZ}}, Hölder's inequality and \eqref{eq:s2} that the above is for $|Q|>q_0$ bounded by
\begin{align}\nonumber
	&\, \le {\a_0} |Q|^{2p} \La(Q) \Big(\iint \E[N_{x,r}^2 g(x,\XX \cup \{x\})] \E[N_{y,r}^2 g(y,\XX \cup\{y\})\one\{R(y,\XX \cup \{y\})>s\}] \,\d x \,\d y)\Big)^{\f 12}\\ \nonumber
	&\, \le {\a_0} |Q|^{2p} \La(Q) \Big(\iint \E[N_{x,r}^4]^{\f 12} \E[g(x,\XX \cup \{x\})^2]^{\f 12} \E[N_{y,r}^4]^{\f 14} \E[g(y,\XX \cup\{y\})^4]^{\f 14} \sqrt{q_y}\,\d x \,\d y)\Big)^{\f 12}	\\ \label{eq:bound_cov_1}
	&\, \le {\a_0} |Q|^{2p} \La(Q)^2 \sup_{y \in Q} \E[g(y,\XX \cup \{y\})^4]^{1/4}  \sup_{y \in Q} \E[N_{y,r}^4]^{3/8}q_y^{1/4}.
\end{align}
Here we use \eqref{eq:m2} and that $\sup_{y \in Q} \E[N_{y,r}^4] \le c_1 r^{4d}$ and thus find that the above is bounded by $c_2 r^{3d/2} |Q|^{2(p+1)}q_y^{1/4}$, where $c_1>0$ and $c_2>0$ only depend on $\a_0$ and $c_{\ms{m}}$.

For the first covariance on the right-hand side in \eqref{lemco1bou2}, we first work conditioned on $\Psi_x$ and $\Psi_y$
Then, again, we want to apply Theorem \ref{th:disdis} with $A:=B_{2s}(z)$ and $P= B_{2r}(x)$, $B=B_s(x) $ for $z=y$ and $P= B_{2r}(y)$, $B=B_s(y) $ for $z=x$. Setting $\Psi_{z, 0} : = \Psi_z \cap B_{s}(z)$, we consider the PI
$${\k_{1,z}(w, \vp) :=\k(w,\vp_{B_{s}(z)^c}  \cup \Psi_{z})\one\{w\in  Q^+\sm B_{s}(z)\},}$$
 and  the interpolation PI
\begin{align*}
	\k_1'(w,\vp):=\begin{cases}
		\k(w,\vp_{B_{s}(x)^c} \cup  \Psi_{x, 0}) & w \in B_{2r}(x) \sm B_{s}(x),\\
			\k(w,\vp_{B_{s}(y)^c} \cup  \Psi_{y, 0}) &w \in B_{2r}(y) \sm B_{s}(y),\\
			0 &\text{otherwise}.
	\end{cases}
\end{align*}
Then, as before conditioned on $\Psi_{x}$ and $\Psi_{y}$, we  note that on the event $\{R(z,\Psi_{z, 0}) \le s\}$, the process  $\XX_{z}^\Xi$ is a Gibbs process on $Q^+\sm B_{s}(z)$ with PI $\k_{1,z}(w, \vp)$ and boundary conditions $\Psi_{z,0}$.
Hence, when defining 
$$
\XX':=T^{'_1\ms{rad}}_{Q^+,(Q^+)^c,\Psi_{x,0} \cup \Psi_{y,0}}(\PP^*_{Q^+}) \cup \Psi_{x, 0}  \cup  \Psi_{y, 0},$$
where $\k'_1$ is used in the definition of $T^{'_1\ms{rad}}$, Theorem \ref{th:disdis} gives for $S_{z}:=\{\XX_{z}^{\Xi} \cap B_{2s}(z)=\XX' \cap B_{2s}(z)\}$ that 
\begin{align}
\P(S_{z}^c) \le c_{\ms{DP}} |B_{4s}| e^{-c_{\ms{SPT},2}s}.\label{est:Szn}
\end{align}
Similarly as before, we bound the first covariance on the right-hand side in \eqref{lemco1bou2} by
\begin{align*}
	&\Cov(\tilde F_{x} \one_{S_{x}}, \tilde F_{y}\one_{S_{y}}) + \Cov(\tilde F_{x} \one_{S_{x}^c}, \tilde F_{y}\one_{S_{y}}) + \Cov(\tilde F_{x} , \tilde F_{y}\one_{S_{y}^c}).
\end{align*}
Here, we use \eqref{est:Szn} to bound the last two covariances similarly to {\eqref{eq:bound_cov_1}}. Setting $\tilde F_{x}':=\tilde F(\XX'_{x, r})$ {where $\XX'_{x, r}:=\XX'\cap B_{r}(x)$,} the first covariance is further decomposed into
\begin{align}
	\Cov(\tilde F_{x}' , \tilde F_{y}')-	\Cov(\tilde F_{x}' \one_{S_{x}}, \tilde F_{y}'\one_{S_{y}^c}) 
	&-\Cov(\tilde F_{x}' \one_{S_{x}^c}, \tilde F_{y}'\one_{S_{y}}) -	\Cov(\tilde F_{x}' \one_{S_{x}^c}, \tilde F_{y}'\one_{S_{y}^c}). \label{lem:co1co}
\end{align}
Note that the computations for the last three covariances are similar to the ones carried out above. Hence,  we only consider the first covariance in detail.
It is given by
\begin{align}
&\E \big[\Cov(\tilde F_{x}' , \tilde F_{y}'\mid \Psi_x, \Psi_y)\big] + \Cov\big(\E\big[\tilde F_{x}'\mid \Psi_x,\Psi_y\big], \E\big[\tilde F_{y}'\mid \Psi_x,\Psi_y\big] \big).\label{lem:co_cosplit}
\end{align}
To bound the first term, we first note that by the tower property of conditional expectation, 
\begin{align*}
\Cov(\tilde F_{x}' , \tilde F_{y}'\mid \Psi_x, \Psi_y)=\E\big[\tilde F_{y}' 	\big(\E\big[\tilde F_{x}' \mid \Psi_x, \XX'_{y, r} \big]-\E[\tilde F_{x}'\mid \Psi_x, \Psi_y]\big)\mid \Psi_x, \Psi_y \big]. 
\end{align*}	
To bound the difference of the expectations in the round brackets, we apply Corollary \ref{cor:dec} with $A := B_{r}(x)$ and the Gibbs process $\XX'$.  Using here that $\tilde F(\XX'_{x, r}) \le |Q|^p \#(\XX_{x, r}') \le |Q|^p \Pd(B_{r}(x))$ almost surely, we obtain
\begin{align*}
	&\big|	\E\big[\tilde F_{x}' \mid \Psi_x, \XX'_{y, r}  \big]-\E[\tilde F_{x}'\mid \Psi_x, \Psi_y] \big|\\
	&\quad \le |Q|^p\E\big[ \Pd (B_{r}(x)) \dtv\big(\mc L(\XX'_{x,  r}\hspace{-.15cm}\mid \Psi_x, \XX'_{y, r}) ,\mc L(\XX'_{x, r}\hspace{-.15cm}\mid \Psi_x, \Psi_y )\big)\big]\\
	&\quad \le |Q|^p c_{\ms{SPT},1} |B_{r+r_0}(x)| e^{-2c_{\ms{SPT},2}s}\E[ \Pd (B_{r}(x))] \le c_3 c_{\ms{SPT},1} e^{-c_{\ms{SPT},2} s} r^{2d} |Q|^p.
\end{align*}
In particular,
$
\Cov(\tilde F_{x}' , \tilde F_{y}'\mid \Psi_x, \Psi_y)\le c_4 c_{\ms{SPT},1}e^{-c_{\ms{SPT},2} s}  r^{3d} |Q|^{2p}
$
for some $c_3=c_3(\a_0,r_0,d)>0$ and $c_4=c_4(\a_0,r_0,d)>0$. We split the second term from \eqref{lem:co_cosplit} into
\begin{align}
&\Cov\big(\E\big[\tilde F_{x}'\mid \Psi_x\big], \E\big[\tilde F_{y}'\mid \Psi_y\big] \big)+\Cov\big(\E\big[\tilde F_{x}'\mid \Psi_x,\Psi_y\big]- \E\big[\tilde F_{x}'\mid \Psi_x\big], \E\big[\tilde F_{y}'\mid \Psi_y\big] \big)\nonumber\\
&\quad +\Cov\big(\E\big[\tilde F_{x}'\mid \Psi_x,\Psi_y\big], \E\big[\tilde F_{y}'\mid \Psi_x,\Psi_y\big] -  \E\big[\tilde F_{y}'\mid \Psi_y\big]\big).\label{lem:co_cond}
\end{align}
Due to the independence of $\Psi_x$ and $\Psi_y$, the first covariance vanishes. The second and the third covariances can be treated very similarly, so we only consider the second one. Here, we use that $\tilde F(\XX_{x, r}') \le |Q|^p \#(\XX_{x, r}') \le |Q|^p \Pd(B_{r}(x))$ and obtain the bound
\begin{align*}
	&\E\big[\big(\E\big[\tilde F_{x}'\mid \Psi_x,\Psi_y\big]- \E\big[\tilde F_{x}'\mid \Psi_x\big]\big)^2\big]^{1/2}{ \E\big[(\tilde F_{y}')^2 \big]^{1/2}} \\
	&\le |Q|^{2p} \E\big[\Pd(B_{r}(x))^2 \dtv\big(\mc L(\XX'_{x, r}\mid \Psi_x,\Psi_y),\mc L(\XX'_{x, r}\mid \Psi_x ) \big)\big]^{1/2}  \E [ \Pd(B_{r }(x))^2]^{1/2}.
\end{align*}
By Corollary \ref{cor:dec}, the total variation distance is    bounded by $c_{\ms{SPT},1} |B_{r+r_0}(x)|e^{-2c_{\ms{SPT},2} s}$ uniformly in $\Psi_y$. Therefore, we obtain the bound
$
c_5 c_{\ms{SPT},1} |Q|^{2p} r^{4d} e^{-2c_{\ms{SPT},2}s}
$
for some $c_5=c_5(\a_0,r_0)$.
\end{proof}

%
%
\bep[Proof of Theorem \ref{thm:2}] 
It remains to establish the bounds from \eqref{E145}.
 \medskip
 
 %
%
\ni{$\bs E_1$.} 
First, we have that 
\begin{align}
	\s^{4}E_1^2&\le \E\Big[\Big(\int_{Q} (Y_{x}\, \one \{|Y_{x}|\le \s\}-\E\big[Y_{x} \, \one \{|Y_{x}|\le \s\}\big] )\,\La(\d x)\Big)^2\Big]\nonumber\\
	&= \iint_{Q^2}  \Cov(Y_{x}  \one \{|Y_{x}|\le \s\},Y_{y} \one \{|Y_{y}|\le \s\})\,\La(\d x) \La(\d y).\label{eqn:E1bou0}
\end{align}
It suffices to show that the right-hand side is bounded by $c_{\ms{Y}} \s^{-2} |Q|^2 r^{4d}$. We write $Y_{x}^\le:=Y_{x}  \one \{|Y_{x}|\le \s\}$ and $Y_{x}^>:=Y_{x}  \one \{|Y_{x}|> \s\}$, $x \in \R^d$, and reformulate the covariance in \eqref{eqn:E1bou0} by
\begin{align}
  \Cov(Y_{x},Y_{y})-
  \Cov(Y^\le_{x},Y^>_{y})-   \Cov(Y^>_{x},Y_{y}).\label{eqn:E1bouu}
\end{align}
The second and third covariances can be treated similarly. Hence, we only discuss the integral over the first and the second covariance. Here, we use the Hölder inequality to bound the covariance in the second integral by
\begin{align*}
	\Cov(Y^\le_{x},Y^>_{y})&\le \E[|Y^\le_{x}-\E Y^\le_{x}|^4]^{1/4}\E \Big[\Big((|Y^>_{y}-\E Y^>_{y}|)^{1/3}\Big)^4\Big]^{3/4}\\
	&\le \E[(|Y^\le_{x}|+|\E Y^\le_{x}|)^4]^{1/4}\E\Big[\Big(|Y^>_{y}|^{1/3}+|\E Y^>_{y}|^{1/3}\Big)^4\Big]^{3/4}.
\end{align*}
Now, we evaluate the fourth powers under the expectations and apply Jensen's inequality with the convex mappings $z\mapsto z^4$ and $z \mapsto z^{4/3}$ to each of the resulting terms. This yields the bound
\begin{align*}
(16 \E [|Y^\le_{x}|^4])^{1/4} (16\E [|Y^>_{y}|^{4/3}])^{3/4} &\le 16 \E[Y_{x}^4]^{1/4}\E[|Y_{y}^>|^{4/3}]^{3/4}\\
& \le 16\E[Y_{x}^4]^{1/4}\E[Y_{y}^4]^{1/4} \P(Y_{y}^4>\s^4)^{1/2}.
\end{align*}
Hence, bounding the last probability by Markov's inequality, we find 
from Lemma \ref{lem:Yintbou} and  from Jensen's inequality applied to the normalized $\La$-integrals, that the integral over the second covariance in \eqref{eqn:E1bouu} is bounded by
\begin{align*}
	{16}{\s^{-2}}	 \iint_{Q^2}  \E[Y_{x}^4]^{1/4}\E[Y_{y}^4]^{3/4}\,\La(\d x) \La(\d y)\le {16\La(Q)}{\s^{-2}} \int_{Q}  \E[Y_{x}^4]\,\La(\d x) \le 16 \a_0 c_{\ms{Y}} \s^{-2}|Q|^2 s^{4d}.
\end{align*}
It remains to bound the integral over the first covariance in \eqref{eqn:E1bouu}, which we split by
 \begin{align}
	 \int_{Q^2} \one\{|x-y|\le 4r\} \Cov(Y_{x},Y_{y})\,\La(\d x)\La(\d y)+ \int_{Q}\int_{Q \sm B_{4r}(y)}  \Cov(Y_{x},Y_{y})\,\La(\d x) \La(\d y).\label{eqn:E1Ysplit}
 \end{align}
Here, the Cauchy-Schwarz inequality gives $\Cov(Y_{x},Y_{y})\le \sqrt{\E [Y_{x}^2]} \sqrt{\E [Y_{y}^2]}\le \E [Y_{x}^2]+\E [Y_{y}^2]$. Hence, we can bound the first integral in \eqref{eqn:E1Ysplit} using Lemma \ref{lem:Yintbou} by
$$
2|B_{4r}(y)| \int_{Q} \E[Y_{x}^2]\, \La(\d x) \le 2 \k_d 16^d c_{\ms{Y}} |Q| s^{3d}.
$$
To deal with the second integral in \eqref{eqn:E1Ysplit}, we set $\Xi_{x, r}^{\le} :=\Xi_{x}^\le \big(B_{r}(x)\big)$, and split $Y_{x}$ as 
\begin{align*}
Y_{x} =Y_{x}^{>,*}+Y_{x}^{\le, \ms{in}}+Y_{x}^{\le, \ms{out}} &:=\Big(|\Xi_{x}^> |-|\Xi^>|\Big) + \Big(\Xi_{x, r}^\le-\Xi^\le \big( B_{r}(x)\Big) \\
&\quad +\Big(\Xi_{x}^\le \big( B_{r}(x)^c\big)-\Xi^\le \big(  B_{r}(x)^c\big)\Big) \one \{\XX \Delta \XX_{x}^\Xi \not \su B_{3s}(x)\},
\end{align*}
and bound the second integral in \eqref{eqn:E1Ysplit} by 
\begin{align}
&\int_{Q^2} \bigg(\Cov(Y_{x}^{>,*}, Y_{y}^{>,*}) + 2 \Cov(Y_{x}^{>,*}, Y_{y}^{\le, \ms{in}}) + 2 \Cov(Y_{x}^{>,*}, Y_{y}^{\le, \ms{out}}) +  \Cov(Y_{x}^{\le, \ms{out}}, Y_{y}^{\le, \ms{out}}) \nonumber\\
&\quad + 2 \Cov(Y_{x}^{\le, \ms{out}}, Y_{y}^{\le, \ms{in}})\bigg) \,\La(\d x) \La(\d y) +  \int_{Q}\int_{Q \sm B_{4r}(y)} \Cov(Y_{x}^{\le, \ms{in}}, Y_{y}^{\le, \ms{in}}) \,\La(\d x) \La(\d y).
\end{align}
Now, the Cauchy-Schwarz inequality together with the bounds \eqref{eqn:Yint3}, \eqref{eqn:Yint5} and \eqref{eqn:Yint5b} in the proof of Lemma \ref{lem:Yintbou} give that the first integral is bounded by $c_{\ms{Y}} |Q|s^{4d}$.  For the covariance in the second integral, we  rely on the decomposition
 \begin{align*}
	 \Cov(Y^{\le, \ms{in}}_{x} ,Y^{\le, \ms{in}}_{y})&= \Cov\big(\Xi_{x, r}^\le, \Xi_{y, r}^\le\big) -\Cov\big(\Xi_{x, r}^\le, \Xi^\le \big( B_{r}(y)\big)\big) \\
	 &\quad -\Cov\big(\Xi^\le \big( B_{r}(x)\big), \Xi_{y,r}^\le\big)+\Cov\big(\Xi^\le \big( B_{r}(x)\big), \Xi^\le \big( B_{r}(y)\big)\big).
 \end{align*}
 Since the arguments for the other terms work analogously, we only
consider the first term on the right-hand side. As a first step, we let $M:=|Q|^{10}$ and write
\begin{align*}
\Xi_{x,r}^\le &= \sum_{w \in \X_x} \one\{w \in B_r(x)\}g(w,\XX_{x}) \one\{R(w,\XX_{x})\le s, g(w,\XX_{x})\le M\}  \\
&\quad + \sum_{w \in \X_x} \one\{w \in B_r(x)\} g(w,\XX_{x}) \one\{R(w,\XX_{x})\le s, g(w,\XX_{x})>M\}=:\Xi_{x,r}^{\le M} + \Xi_{x,r}^{>M}.
\end{align*}
This gives
\begin{align}
&	\Cov\big(\Xi_{x, r}^\le,  \Xi_{y, r}^\le\big) =\Cov\big(\Xi_{x, r}^{\le M}, \Xi_{y,  r}^{\le M} )+\Cov\big(\Xi_{x,r}^{\le M},\Xi_{y,  r}^{> M}) + \Cov\big(\Xi_{x, r}^{> M} , \Xi_{y,r} ).\label{thm2:E1tr}
\end{align}
The $\La^2$-integral of the first covariance is by Lemma \ref{lem:thm2cov} and Remark \ref{rem:cov} bounded by $|Q|^{-1}$ for $|Q|$ large enough. The integrated second covariance term is by the Jensen inequality bounded by
\begin{align}
	\int_{Q^2}	\sqrt{\Var(\Xi_{x, r}^{\le M}) 	\Var(\Xi_{y, r}^{> M})}  \La^2(\d x, \d y) \le\La(Q) 	\sqrt{\int_Q \E\big[\big(\Xi_{x,r}^{\le M}\big)^2\big] \La(\d x){\int_Q \E\Big[\big(\Xi_{y,r}^{> M}\big)^2\Big] \La(\d y)}}. \label{eqn:thm2intM}
\end{align}
Here we find by a similar computation as in \eqref{eqn:firstfactor} that the first integral is bounded by $c_1  |Q| r^{2d}$ for some $c_1=c_1(d,\a_0,c_{\ms{m}})$. 
For the second integral, we obtain similarly to \eqref{eqn:firstfactor} the bound
\begin{align*}
	&\int_{Q} \E\Big[\big(\Xi_{y, r}^{> M}\big)^2\Big] \La(\d y)\\
	& \quad  \le |Q| \sum_{k=0}^{2} \sum_{\substack{(i_1,\dots,i_{k + 1}):\\ i_1+\cdots+i_{k + 1}=2}} \hspace{-0.3cm} \Big\{\a_0^{k+1} |B_{r}(y)|^{k} \\
	&\quad \quad \qquad \qquad \qquad \qquad\times \sup_{\xx \in Q^{k+1}}\E\Big[  g(x_{k+1},\XX \cup \{\xx\})^{i_{k + 1}} \prod_{j=1}^{k}g(x_j,\XX \cup \{\xx\})^{i_j} \one\{g(x_j,\XX \cup \{\xx\}) > M\}\Big]\Big\}.
\end{align*}
Here, we apply the Hölder inequality and the Markov inequality to the expectations above and obtain from \eqref{eq:m2} that the second integral from \eqref{eqn:thm2intM} is bounded by $c_2 r^{2d} |Q|^{-3}$ for some $c_2=c_2(d,\a_0,c_{\ms{m}})$. 
The same bound can be established for the third covariance on the right-hand side in \eqref{thm2:E1tr}.  
Therefore, for some $c_3=c_3(d,\a_0,c_{\ms{m}})$,
$$
\int_{Q}\int_{Q \sm B_{4r}(y)}\Cov(Y_{x}^{\le, \ms{in}}, Y_{y}^{\le, \ms{in}}) \,\La(\d x) \La(\d y) \le c_3 r^{2d}.
$$
Since $r=4s$ and $\s^2=\Var(\Xi) \le \E[\Xi^2]\le c_4 |Q|$ for some $c_4=c_4(\a_0,c_{\ms{m}})$, this shows the first bond in \eqref{E145}.

\bigskip

\ni{$\bs E_4$.}
We perform the outer integration in $E_4$ only from $0$ to $1$. The remaining integral from $-1$ to $0$ is bounded analogously. We deal separately with the cases $|x-y| \le 4 r$ and $|x-y| > 4 r$. First, consider the case  $|x - y|\le 4r$.  
By the bound $\Cov(\phi_x(t),\phi_y(t))\le \E[\phi_x(t)]$, \eqref{eq:m2} and Lemma \ref{lem:Yintbou},
\begin{align*}
\s^{-2}\int_0^1\int_{Q}\int_{B_{4r}(x)} \Cov\big(\phi_{x}(t),\phi_{y}(t))\, \La(\d y)\,\La(\d x) \d t &\le \s^{-2} \int_0^1 \int_Q \La(B_{4r}(x)) \E[\phi_x(t)] \La(\d x) \d t\\
&\le c_{\ms{m}} (4r)^d \k _d \s^{-3}\int_Q \E[|Y_x|]\, \La(\d x) \le c_5 r^{2d} \s^{-3}|Q|
\end{align*}
for some $c_5=c_5(d,\a_0,c_{\ms{m}})>0$.

Finally, we consider the case $|x - y|> 4r$. Writing $U_{x}:=\{Y_{x}^{>,*} = Y_{x}^{\le, \ms{out}}=0\}$,
we decompose the covariance  $ \Cov(\phi_{x}(t),\phi_{y}(t))$ as 
\begin{align}
\Cov(\phi_{x}(t) \one\{U_{x}^c\},\phi_{y}(t))
	+	  \Cov(\phi_{x}(t) \one\{U_{x}\},\phi_{y}(t)\one\{U_{y}^c\})
 + 	  \Cov(\phi_{x}(t) \one\{U_{x}\},\phi_{y}(t)\one\{U_{y}\}).\label{eq:thm2E4}
\end{align}
Using that 
\begin{align*}
\Cov(\phi_{x}(t) \one\{U_{x}^c\},\phi_{y}(t)) \le \P(U_{x}^c) &\le \E\Big[\sum_{y \in \XX_x^{\Xi}} \one\{R(y,\XX_x^{\Xi})>s \text{ or } \XX \Delta \XX_x^{\Xi} \not \subseteq B_{3s}(x)\} \Big]\\
&\le \E\Big[\sum_{y \in \XX_x^{\Xi}} \one\{R(y,\XX_x^{\Xi})>s\}\Big] + \E\Big[\sum_{y \in \XX_x^{\Xi}} \one\{ \XX \Delta \XX_x^{\Xi} \not \subseteq B_{3s}(x)\} \Big],
\end{align*}
we find analogously to the bounds of the third and fifth term in \eqref{eqn:minpalmloc} with $g \equiv 1$ that
$$
	\int_0^1 \int_Q \int_{Q \setminus B_{4r}(y)}   |\Cov(\phi_{x}(t) \one\{U_n^c\},\phi_{y}(t))| \La(\d x) \La (\d y)\d t \le c_6 |Q| s^{2d}.
$$
for some constant $c_6=c_6(d,\a_0,c_{\ms{m}})>0$. A similar estimate holds for the second covariance term. Finally, we use that on $U_x$ it holds that $\De_x=\De_{x}^\le:=(|\Xi_{x}^\le|-|\Xi^\le|)/{{\sigma}}$, and hence
\begin{align*}
\phi_x(t)=\phi_{x}^\le(t)=\begin{cases}
		\one\{1 \ge \De_{x}^\le>t>0\}&t>0,\\
		\one\{-1 \le \De_{x}^\le<t<0\}&t<0,
	\end{cases}
\end{align*}
and decompose the third covariance as
\begin{align*}
  \Cov(\phi_{x}^\le(t) ,\phi_{y}^\le(t)) -   \Cov(\phi_{x}^\le(t) \one\{U_{x}^c\},\phi_{y}^\le(t)\one\{U_{y}\})- \Cov(\phi_{x}^\le(t),\phi_{y}^\le(t)\one\{U_{y}^c\}).
\end{align*}
Here, the first term is by Lemma \ref{lem:thm2cov} with $p=0$  for $|Q|$ large enough bounded by $|Q|^{-1}$. Note that this result indeed applies since $(|\Xi_{x}^\le|-|\Xi^\le|)/{{\sigma}}$ is a function of $\XX_x \cap B_r(x)$. The second and the third covariance\cb{s} are bounded similarly to the first and second covariance\cb{s} in \eqref{eq:thm2E4}.
This gives the bound on $E_4$ that was asserted in \eqref{E145}.
\medskip

%
%
\ni{$\bs E_5$.} The steps are very similar to the arguments for $E_4$.  We deal separately with the cases $|x-y| \le 4 r$ and $|x-y| > 4 r$. First, consider the case where $|x - y|\le 4r$. Then,
$$\Big|\int_0^1 t\Cov(\phi_{x}(t),\phi_{y}(t)) \d t\Big| \le \int_0^1 t\E[\phi_{x}(t)]\d t
	\le  0.5\s^{-2} \E[Y_{x}^2].$$
Therefore, by Lemma \ref{lem:Yintbou},
\begin{align*}
\s^2\int_0^1\int_{Q}\int_{B_{4r}(x)} t\,\Cov(\phi_{x}(t),\phi_{y}(t))\, \La(\d y)\,\La(\d x) \d t&\le \f 12 \int_Q \La(B_{4r}(x)) \E[Y_x^2] \La(\d x) \le c_7 |Q| s^{3d}
\end{align*}
for some $c_7=c_7(d, \a_0, c_{\ms{m}})$, and an analogous calculation holds for $t \in [-1, 0]$.

Hence, it remains to consider the case where $|x - y|> 4r$. We use the same bound 
for $\Cov(\phi_{x}(t),\phi_{y}(t))$ as in the case $E_4$. Since $\s^2 \le c_4 |Q|$, the asserted bound on $E_5$ follows.
\end{proof}

	\appendix

\section{Relation to Gibbs processes obtained from perfect simulation}
\label{sec:EDD} 
There is a substantial literature available \cite{CX22,gibbs_limit,gibbsCLT} on limit theorems for a class of Gibbs point processes, usually denoted $\Psi^*$ and satisfying the additional condition \cite[(3.7)]{gibbs_limit}, obtained from a perfect simulation technique originally suggested in \cite{ferrari}.   A precise definition of $\Psi^*$ is given in \cite{gibbs_limit}. Whenever we speak of the class $\Psi^*$ below, it is implicitly assumed that \cite[(3.7)]{gibbs_limit} holds. As observed in Section \ref{ssec:gibbs}, for finite-range Gibbs processes, our Assumption \eqref{as:(A)} is more general than the class $\Psi^*$. 

The key achievement of \cite{CX22} is to provide a quantitative CLT with respect to the Wasserstein distance for point processes that are not necessarily of Gibbsian type. The main assumption here is the exponential decay of dependence (EDD). In \cite[Section 3.1]{CX22}, it is explained that the Gibbs point processes in $\Psi^*$  are EDD. Although not stated explicitly, this means that the condition \cite[Equation (1.4)]{gibbs_limit} is also assumed. 
In Corollary \ref{cor:EDD} below, we show using disagreement couplings that also the Gibbs processes satisfying Assumption \eqref{as:(A)} satisfy the EDD, meaning that the results of \cite{CX22} also apply in this case. To make this precise, for disjoint bounded Borel sets $A, B\su \R^d$ we define the $\beta$-mixing coefficient 
$$\beta_{A, B}(\X) := \dtv\big(\XX_{A \cup B}, \XX_A \cup \XX_B'\big),$$
where $\XX_B'$ is an independent copy of $\XX_B$.  Then, the EDD says that there exists a constant $\th_0 > 0$ with the following property. Whenever, $A, B \su \R^d$ are disjoint bounded Borel sets with  $\dist(A, B) \ge \th_0 \log(\diam(A) \vee \diam(B)\vee \th_0)$, then
$$\beta_{A, B}(\X) \le \th_0 (\diam(A) \vee \diam(B) \vee 1)^{\th_0} e^{-\dist(A, B)/\th_0}.$$

%
%
\bec[Assumption \eqref{as:(A)} implies EDD]
        \label{cor:EDD}
        Let $\k$ be a PI satisfying Assumption \eqref{as:(A)}. Then, the associated infinite-volume Gibbs process $\XX$ satisfies the EDD. 
\enc

\bep
	Let $\X^1$ and $\X^2$ be two independent copies of $\X$. Let $\psi_i=\X^i \cap A$ for $i=1,2$. Define $\tilde{\X}^i:=\psi_i\cup T^\ff_{A^c,\psi_i}(\Pds)$, where $\Pds$ is independent of $\X^1$ and $\X^2$. Then $\tilde{\X}_i \sim \X$ by Proposition \ref{pr:chimera} and $\tilde{\X}^2$ is independent of $\tilde{\X}^1_A=\X^1_A$. Hence
	$$
	\dtv(\X_{A\cup B},\X_{A}\cup \X_{B}' ) = \dtv(\tilde{\X}^1_{A\cup B},{\X}_{A}^1\cup \tilde{\X}_{B}^2 )\le \P(\tilde{\X}^1_{ B}\ne \tilde{\X}_{B}^2).  $$
	 By Corollary \ref{cor:ff_disag} (with the roles of $A$ and $B$ switched), 
	$$\P\big(\tilde{\X}^1\cap B \ne \tilde{\X}^2\cap B \mid \psi_1,\psi_2) \le c_{\ms{DPP}}|B_{2\dist(A,B)/3}(B)|\exp(-c_{\ms{SPT}}\dist(A,B)).$$
	This is enough to show the claim since $|B_s(B)| \le c_1(2s+\diam(B))^d\le c_2(s^d\vee \diam(B)^d) $ for suitable $c_1,c_2>0$. This follows because $\diam(B_s(B)) \le \diam(B)+2s$. Using that $|B_s(B)| \le |B_1(o)| \diam(B_s(B))^d$ concludes the proof.
	\enp

As mentioned above, for finite-range interactions, the class of Gibbs processes satisfying Assumption \eqref{as:(A)} is larger than $\Psi^*$. Conversely, our paper requires the interaction range to be finite, whereas $\Psi^*$ also allows some infinite-range Gibbs processes. Currently, the exploration algorithm in the proof of Theorem \ref{th:disdis} heavily relies on this finiteness property. Therefore, while an extension of our arguments to infinite ranges does not seem to be completely inaccessible, it would make the arguments dramatically more complicated so that we leave it for future work. Taking this into account, it is interesting to discuss to what extent our techniques can be extended to infinite-range Gibbs processes  in $\Psi^*$. In particular, we consider the key Theorem \ref{th:disdis}.

Loosely speaking, the idea of the perfect simulation algorithm is that to construct a Gibbs point process in a Borel set $B \su \R^d$, it is only necessary to know a its \emph{ancestor clan} $\bs A_B \su \R^d \times \R$. Since the definition of the ancestor clan would require us to repeat the entire construction of the perfect simulation scheme, we refer the reader to \cite{gibbsCLT} for details. However, we shall only use the following  property: We denote by $\diam(\bs A_B)$ the diameter of the projection of $\bs A_B\su \R^d \times \R$ to the $\R^d$ component.
 Then, with high probability, the diameter $\diam(\bs A_B)$  is not much larger than that of the set $B$. More precisely, \cite[Equation (3.6)]{gibbs_limit} shows that for all Gibbs processes in $\Psi^*$ satisfying the crucial assumption  \cite[Equation (3.7)]{gibbs_limit} there exists some $c_{\ms PS} > 0$ such that for all $r \ge 1$, we have 
\begin{align}
	\label{eq:36}
	\P\big(\diam (\bs A_B) \ge r + \diam(B)\big) \le c_{\ms PS} (|B| \vee 1)\exp(-r/c_{\ms PS}).
\end{align}
Equipped with this information, we have the following version of Theorem \ref{th:disdis}.

\bec \label{cor:disag_ps}
	Consider two PIs $\kappa,\kappa'$ from the class $\Psi^*$. Let $A\su P \su Q\su \R^d$ be bounded Borel sets and $\psi,\psi'\in \Nlf_{Q^c}$ that differ only on $B\su Q^c$. Suppose $\kappa$ and $\kappa'$ satisfy \eqref{eq:k,k'}. Let $\XX(Q, \psi)$ and $\XX'(Q,\psi')$ be the Gibbs processes obtained from the perfect simulation technique using $\kappa$ and $\kappa'$, respectively. Then,
	$$\P\big(\XX(Q, \psi)\cap A \ne \XX'(Q,\psi') \cap A\big) \le  C |B_{r_0}(A)|\exp(-\dist(A,B\cup (Q\sm P))/C).$$
\enc

%
%
\bep[Proof of Corollary \ref{cor:dperc} for perfect-simulation processes]
When coupling $\XX(Q, \psi)\cap A$ and $\XX(Q, \psi')\cap A$ through the perfect simulation, we see that these two processes differ only if the ancestor clan $\bs A_A$ hits the set $(B\cup (Q\sm P))\times \R^d$. Hence, the coupling characterization of the total variation distance gives that
\begin{align*}
\P\big(\XX(Q, \psi)\cap A \ne \XX'(Q,\psi') \cap A\big) &\le \P\big(\bs A_A \cap (B\cup (Q\sm P))\times \R^d \ne \es\big) \\
&\le \P\big(\diam(\bs A_A) \ge \diam (A) +   \dist(A, B\cup (Q\sm P)\big).
\end{align*}
 Now, for any fixed $r_0 >0$, we have $|A| \vee 1\le c_1|B_{r_0}(A)|$ for some $c_1 > 0$. Hence, an application of \eqref{eq:36} concludes the proof.
\enp

The results of \cite{gibbs_limit,gibbsCLT} also required the Gibbs process to be Poisson-like.
That is, (i) $\XX$ is stochastically dominated by a homogeneous Poisson point process, and (ii) there exists $c_{\ms{PL}}, r_1 > 0$ such that for all $r \ge r_1$ and $\vp \in \Nlf_{B_r(x)^c}$, it holds that 
$$\P(\XX_{B_r(x)} = \es \ba \XX_{B_r(x)^c} = \vp) \le  e^{-c_{\ms{PL}}r^d}.$$
We show below that a Gibbs point process $\XX$ satisfying Assumption \ref{as:(A)} is Poisson-like. 

\bepr[Assumption \eqref{as:(A)} implies Poisson-likeness]
\label{pr:like}
Consider a PI $\k$ bounded by $\a_0$ and having $\inf_{y\in \R^d}\kappa(y,\es)>0$. Then, the associated infinite-volume Gibbs process $\XX$ is Poisson-like. 
\enpr
\bep
Part (i) follows immediately from Proposition \ref{pr:uniq}, so we concentrate on part (ii). 

We follow the idea of \cite[Lemma 3.3]{gibbs_limit} and start by choosing $r_1>r_0$. Then, for some $p>0$
$$
\sup_{ \psi\in \Nlf_{B_{r_1}(x)^c} } \P(\XX_{B_{r_1}(x)} = \es \ba \XX_{B_{r_1}(x)^c} = \psi )\le 1-\frac{\P(\PP(B_{r_1}(x))=1)}{\E\big[\alpha_0^{\PP(B_{r_1}(x))}\big]} \int_{B_{r_1-r_0}(x)}\kappa(y,\es) \d y \le p < 1,$$
where we used  $Z_{B_{r_1}(x)}(\psi)\le \E[\alpha_0^{\PP(B_{r_1}(x))}]$.
For $r>r_1$, it is possible to place $\Omega(r^d)$ disjoint balls of radius $r_1$ inside $B_r(x)$. Call them $B_{r_1}(y_1),\ldots, B_{r_1}(y_{k(r)})$. Define events $F=\{\XX_{B_r(x)^c} = \vp\}$ and  $F_{i}=\{\XX_{B_{r_1}(y_1)}=\es \}$. Then, by the DLR-equations \eqref{eq:DLR2},
$$\P(\XX_{B_r(x)} = \es \ba \XX_{B_r(x)^c} = \vp)\le \P\Big(\bigcap_i F_i\mid F\Big)\le \prod_i \P\Big(F_i| F\cup \bigcup_{j<i} F_j \Big)\le p^{k(r)}.$$
The result follows because $k(r)$ is $\Omega(r^d)$. 
\enp

\cb{
	\section{Proof of Proposition \ref{pr:iota}}
	\label{sec:iota}
}

\begin{proof}\cb{[Proof of Proposition \ref{pr:iota}]}
	We define a partial ordering \cb{$\le_\iota$} of $\R^d$ by \cb{$x<_\iota y$} if $\iota(x) <\iota (y)$ or $x=y$. If $\lambda( \iota^{-1}(r) )= 0$ for all $r\in \R$, \cb{then almost surely,}
	$\iota$ induces a total ordering on the points of $\PP^d$. Moreover, the window $\cb{Q_{(x,\ff)}} = \{y\in Q\mid i(x)< i(y)\}$ is again a measurable set. In fact, \cb{letting $\mu\in \Nlf$ be locally finite,} $(x,\mu)\mapsto \mu{\cap Q_{\cb{(x,\ff)}}} $ is measurable and hence the thinning probabilities $p(x,\psi)$ are measurable maps. Thus, using  $\iota$ for the \cb{standard Poisson} embedding is well-defined.
	
	\cb{
		
		We may assume $\iota(Q)\su (0,1]$. For each $n=1,2,\ldots $, we construct an \cb{injective} measurable map $\iota^{n}: Q \to \R$  which  provides almost the same ordering as  $\iota$. Let \cb{$\tilde{\iota}:Q \to (0,1]$ be any injective measurable} map. For $x\in Q$, there is a unique $m\in \{0,\ldots,2^{n-1}\}$ such that $x\in\iota^{-1}(m2^{-n+1}-2^{-n},m2^{-n+1}+2^{-n}]$, and we define $\iota^{n}(x) = m + \phi(x) \in (m,m+1]$. Letting $A_m=\iota^{-1}(m2^{-n+1}-2^{-n},m2^{-n+1}+2^{-n}]$, we have that $x<_{\iota^n}y$ if either $x\in A_{m_1}$ and $y\in A_{m_2}$ with $m_1<m_2$, or $x,y\in A_m$, and $\tilde{\iota}(x)<\tilde{\iota}(y)$.

		For each $n$, the ordering induced by $\iota^n$ may be used in the Poisson embedding to obtain $\XX^{n}(Q,\psi):=T_{Q,\psi,\iota^n}(\Pd_Q)$, which we know has the correct distribution $\XX(Q,\psi)$. We must show three things:
		\begin{itemize}
			\item[(i)] $\XX^{n}(Q,\psi)$ converges almost surely to a point process $\XX^{\ff}(Q,\psi)$.
			\item[(ii)] The limit  $\XX^{\ff}(Q,\psi)$ is the same point process $T_{Q,\psi,\iota}(\Pd_Q)$ we would get if we had used the ordering $\iota$. 
			\item[(iii)] The limit  $\XX^{\ff}(Q,\psi)$ again has the correct distribution $\XX(Q,\psi)$.
		\end{itemize}

		Let $Q_{(-\ff,x)}=\{y\in Q\mid i(y)< i(x)\}$, $Q_{(x,\ff)}^n=\{y\in Q\mid\iota^n(y) >\iota^n(x) \}$, and $Q_{(-\ff,x)}^{n}=\{y\in Q\mid\iota^n(y) <\iota^n(x) \}$. We claim that for any fixed $x\in Q$, 
		\begin{equation}\label{eq:symdif}
			\lim_{n \to \ff }|Q_{(x,\ff)}\Delta Q_{(x,\ff)}^{n}|=0, \qquad \lim_{n \to \ff }|Q_{(-\ff,x)}\Delta Q_{(-\ff,x)}^{n}|=0.
		\end{equation}
		Note first that $Q_{(x,\ff)}\Delta Q_{(x,\ff)}^{n}=Q_{(-\ff,x)}\Delta Q_{(-\ff,x)}^{n}$.
		Suppose $i(x)\in (m_{n,x}2^{-n+1}-2^{-n},m_{n,x}2^{-n+1}+2^{-n}]$. Then, 
		\begin{align*}
			Q_{(x,\ff)} \sm Q_{(x,\ff)}^n &\su \iota^{-1}(i(x) , m_{n,x}2^{-n+1}+2^{-n}]\\ 
			Q_{(x,\ff)}^n \sm Q_{(x,\ff)} &\su \iota^{-1}(m_{n,x}2^{-n+1}-2^{-n},i(x)].
		\end{align*}
		In total,
		\begin{align}\label{eq:symdif2}
			Q_{(x,\ff)}^n \Delta Q_{(x,\ff)}  \su \iota^{-1}((m_{n,x}2^{-n+1}-2^{-n},m_{n,x}2^{-n+1}+2^{-n}]). 
		\end{align}
		The claim now follows because the sets 
		\begin{align*}
			\iota^{-1}((m_{n,x}2^{-n+1}-2^{-n},m_{n,x}2^{-n+1}+2^{-n}])
		\end{align*}
		decrease towards $\iota^{-1}(\iota(x))$, which has Lebesgue measure zero by assumption. The claim for $Q_{(-\ff,x)}\Delta Q_{(-\ff,x)}^{n}$ is shown similarly. 
		
		Now fix $(x,\psi)$ with $\psi \in \Nlf_{Q_{(-\ff,x)}}$. By  \eqref{eq:symdif2}, there is an $n_0$ such that also $\psi \in  \Nlf_{Q_{(-\ff,x)}^n} $ for $n> n_0$. For such $n$, we have  
		\begin{align*}
			p^n(x, \psi) := \k(x,\psi) \frac{\E[e^{-J(\mu \cap {Q_{(x,\ff)}^n},\psi + \delta_x)}]}{\E[e^{-J(\mu\cap {Q_{(x,\ff)}^n},\psi )}]} \to \k(x,\psi) \frac{\E[e^{-J(\mu\cap {Q_{(x,\ff)}},\psi + \delta_x)}]}{\E[e^{-J(\mu\cap {Q_{(x,\ff)}},\psi )}]}= p(x,\psi) 
		\end{align*}
		as $n \to \infty$. Indeed, looking at the denominator first, we have
		\begin{align*}
			\E\left[e^{-J(\mu\cap {Q_{(x,\ff)}^n},\psi )}\right] = {}&\E\left[e^{-J(\mu \cap {Q_{(x,\ff)}},\psi )}\right] + \E\left[e^{-J(\mu \cap {Q_{(x,\ff)}^n },\psi )}1_{\mu \cap (Q_{(x,\ff)} \Delta Q_{(x,\ff)}^{n}) \ne \es}\right] \\
			&- \E\left[e^{-J(\mu_\cap {Q_{(x,\ff)}},\psi )}1_{\mu\cap (Q_{(x,\ff)} \Delta Q_{(x,\ff)}^n)\ne \es}\right] .
		\end{align*}
		The middle term goes to zero by an application of the Cauchy-Schwartz inequality, since the stability assumption ensures that $\E[e^{-2J(\mu \cap {Q_{(x,\ff)}^n},\psi )}]<\ff$ and from \eqref{eq:symdif} we have that
		$$P(\mu\cap (Q_{(x,\ff)} \Delta Q_{(x,\ff)}^{n})\ne \es) = 1-e^{-|Q_{(x,\ff)} \Delta Q_{(x,\ff)}^{n}|},$$ which goes to zero when $n\to \ff$. The last term goes to zero by the same reasoning. The numerator is treated similarly.
		
		We now show (i). For almost all $\omega$, $\iota$ and $\iota^{n}$ will define the same ordering of $\PP^d_Q(\omega)$  whenever $n>N(\omega) $ according to \eqref{eq:symdif2}. Moreover, since the thinning probabilities converge pointwise, for almost all $\omega$ there is an $N'(\omega)\ge N(\omega) $ such that the thinned process does not change whenever  $n \ge N'(\omega)$. We define $\XX^\ff(Q,\psi)$ to be this limiting point process and note that it is a.s.\ the same as one would have obtained by using the ordering $\iota$ for the thinning, which also shows (ii). 
		
		It remains  to show (iii), that $\XX^\ff(Q,\psi)$ has the correct Gibbs distribution. For this, let $f: \Nlf_Q \to \R$ be a bounded measurable function. Then 
		\begin{align*}
			&|\E[ f( \XX^n (Q,\psi))] - \E[f( \XX^\ff (Q,\psi))] |  \le
			2|f|_{\ff}P( \XX^n (Q,\psi) \ne  \XX^\ff (Q,\psi)).  
		\end{align*}
		The right hand side goes to $0$ and the left hand side is constant since all $\XX^n(Q,\psi)$ have the same distribution $\XX(Q,\psi)$, so it follows that also $\XX^\ff(Q,\psi)$ must have the same distribution.
	}
\end{proof}

	\bibliographystyle{plainnat}
	\bibliography{./lit}
\end{document}